# ADDRESSING MATHEMATICAL INCONSISTENCY: CANTOR AND GÖDEL REFUTED

J. A. PEREZ

ABSTRACT. This article undertakes a critical reappraisal of arguments in support of Cantor's theory of transfinite numbers. The following results are reported:
– Cantor's proofs of nondenumerability are refuted by analyzing the logical inconsistencies in implementation of the *reductio* method of proof and by identifying errors. Particular attention is given to the diagonalization argument and to the interpretation of the axiom of infinity.
– Three constructive proofs have been designed that support the denumerability of the power set of the natural numbers, $\mathcal{P}(\mathbb{N})$, thus implying the denumerability of the set of the real numbers $\mathbb{R}$. These results lead to a Theorem of the Continuum that supersedes Cantor's Continuum Hypothesis and establishes the countable nature of the real number line, suggesting that all infinite sets are denumerable.
Some immediate implications of denumerability are discussed:
– Valid proofs should not include inconceivable statements, defined as statements that can be found to be false and always lead to contradiction. This is formalized in a Principle of Conceivable Proof.
– Substantial simplification of the axiomatic principles of set theory can be achieved by excluding transfinite numbers. To facilitate the comparison of sets, infinite as well as finite, the concept of relative cardinality is introduced.
– Proofs of incompleteness that use diagonal arguments (e.g. those used in Gödel's Theorems) are refuted. A constructive proof, based on the denumerability of $\mathcal{P}(\mathbb{N})$, is presented to demonstrate the existence of a theory of first-order arithmetic that is consistent, sound, negation-complete, decidable and (assumed p.r. adequate) able to prove its own consistency. Such a result reinstates Hilbert's Programme and brings arithmetic completeness to the forefront of mathematics.

## 1. THE INCONSISTENCY OF ZFC

The consistency of classical set theory, more precisely the combination of the Zermelo-Fraenkel axioms and the Axiom of Choice (ZFC) [35,39,56], has recently been challenged by the construction of a proof of Goodstein's Theorem in first-order arithmetic [55]. Goodstein's Theorem is a statement about basic arithmetic sequences and, as such, has a purely number-theoretic character [36,56,59]. It was originally proved, not by means of first-order arithmetic, but by using the well-ordered properties of transfinite ordinals [36,56]. It was later shown by Paris and Kirby to be unprovable in Peano Arithmetic (PA), with a proof that needs both transfinite induction and Gödel's Second Incompleteness Theorem [46,59]. A previous result was required, i.e. the proof of the consistency of PA, also achieved by transfinite induction [29,30]. Consequently, if PA were capable of proving Goodstein's Theorem, PA would also be able to prove its own consistency, thus violating Gödel's Second Theorem [46]. Hence, the claim was made that PA cannot prove Goodstein's Theorem. Since this claim has now been disproved [55], it is important to establish the cause of the inconsistency of ZFC. While the problem may lie with the formulation of the axioms of ZFC, it is also possible that it lies in one or more of the basic concepts underlying the claim of Paris and Kirby [46]. This article undertakes a systematic review of these key concepts, namely transfinite number theory and mathematical incompleteness.

The article starts with a review of Cantorian mathematics (Section 2). Section 3 deals with Cantor's proofs of nondenumerability and provides a methodical refutation of all the main rationales by showing logical inconsistencies and by identifying errors. Section 4 describes constructive arguments that prove the countable nature of the power set of the natural numbers $\mathcal{P}(\mathbb{N})$ and, consequently, of the set of real numbers $\mathbb{R}$.

The refutation of Cantor's proofs and the demonstration of denumerability have profound implications for many areas of mathematics. Section 5 proposes a simple and practical principle to avoid flawed implementation of the *reductio* method of proof. Section 6 introduces a number of conjectures, principles and definitions that may lead to the simplification of the axiomatic principles of set theory, and provide the means to compare quantitatively the size of infinite sets. Cantor's diagonalization argument was widely adopted as a method of proof in the field of logic [20,37] and its implementation is central to important results on mathematical incompleteness by Gödel and others [31,59]. The refutation of the diagonalization argument is extended in Section 7 to the various proofs of arithmetic incompleteness. The denumerability of $\mathcal{P}(\mathbb{N})$ is used to construct a proof for the completeness of first-order arithmetic that confirms the above refutations and reinstates Hilbert's Programme [41,59].

The results reported in this article, combined with the elementary proof of Goodstein's Theorem already cited [55], achieve the mathematical consistency that Cantorian mathematics failed to deliver.

## 2. Cantorian Mathematics

When Georg Cantor established the foundations of set theory, his primary objective was to provide rigorous operational tools with which to characterise the evasive and largely paradoxical concept of infinity [5,23]. He built on the work of Dedekind who, by elaborating on the observations made by Galileo and Bolzano on the nature of infinite systems [5], defined a set $S$ as infinite when it is possible to put all its elements into one-to-one correspondence with the elements of at least one of its proper subsets; that is, to show that both set and subset have the same cardinality [5,26,27]. This property (not shared by finite sets) provided the formal definition of infinity which Cantor used as the starting point for his development of the subject.

It is important to emphasize that the critical issue at the heart of Cantor's treatment of infinity is the ability (or inability) to find a one-to-one correspondence (bijection) between two infinite sets, thus determining whether or not the cardinality of one set is greater than the cardinality of the other [2,5,35]. In this way, Cantor was able to prove, for the first time, that the set of rational numbers $\mathbb{Q}$ has the same cardinality as the set of natural numbers $\mathbb{N}$, despite having a much higher density on the number line [2,6,35]. The same could be said of the set of algebraic real numbers [6]. Accordingly, Cantor coined the term 'denumerable' for all those infinite sets which could be put into one-to-one correspondence with $\mathbb{N}$ [2,6,44].

Cantor's correspondence with Dedekind [28] shows how the failure of his repeated attempts to establish the denumerability of the set $\mathbb{R}$ of all real numbers (transcendental as well as algebraic) eventually drove him to contemplate the possibility of the opposite being true, that is, that the cardinality of $\mathbb{R}$ (which Cantor named $\boldsymbol{c}$, the power of the continuum) might be greater than the cardinality of $\mathbb{N}$ (which he named $\aleph_0$, the first 'aleph'). Cantor published the first *reductio* proof supporting nondenumerability in 1874, using the construction of a sequence of nested closed bounded intervals on the real number line, which allowed him to claim that $\boldsymbol{c} > \aleph_0$ [6].

The subsequent development of Cantor's theory of transfinite numbers (ordinal and cardinal) depended entirely on his ability to prove conclusively the existence of infinite sets of ever growing cardinality. His first assumption was that infinities larger than $\boldsymbol{c}$ would be found when moving from the real number line to the two-dimensional plane and, from there, to higher dimensions. However, Cantor was very surprised to find this was not the case [7,28], as he confided to Dedekind in 1877 when he famously wrote "*Je le vois, mais je ne le crois pas*" [28].

Cantor's failure to justify greater infinities (beyond $\boldsymbol{c}$) was finally overcome in 1891 by the publication of his famous diagonalization argument [8,28]. This new proof succeeded by extending his previous result of 1874 to ever higher infinities via the cardinality of the power set [8]. Cantor's Theorem, as it was coined in 1908 by Zermelo when describing his own more set-theoretical proof [62], states that the cardinality of the power set $\mathcal{P}(A)$ of any set $A$ is larger than the cardinality of set $A$ itself, i.e. $|A| < |\mathcal{P}(A)|$. When applied to $\mathbb{N}$, this theorem translates into $|\mathbb{N}| < |\mathcal{P}(\mathbb{N})|$ (i.e. $\aleph_0 < 2^{\aleph_0}$). Since it is possible to prove that $|\mathbb{R}| = |\mathcal{P}(\mathbb{N})|$ [35,49], the implication is that $\boldsymbol{c} = 2^{\aleph_0} > \aleph_0$, which brings again Cantor's previous result. The importance of Cantor's Theorem rests on the fact that the application of the same rule to subsequent sets (starting with $|\mathcal{P}(\mathbb{R})| = |\mathcal{P}(\mathcal{P}(\mathbb{N}))| > |\mathcal{P}(\mathbb{N})|$, and so on) generates a never ending sequence of larger and larger cardinalities. These were the higher infinities that Cantor had failed to identify in dimensions beyond the real number line [7]. Cantor's Theorem makes the existence of power sets an essential ingredient of the foundational principles of set theory, as highlighted by its formulation among the Zermelo-Fraenkel (ZF) axioms [35,56].

Cantor first proposed his Continuum Hypothesis (CH) in 1878, and spent

most of his working life trying to prove it without success [7,23]. His theory of transfinite ordinals generates a well-ordered set of cardinals (alephs) $\aleph_\lambda$ ($\lambda$ denotes a given infinite ordinal) which, accepting the Axiom of Choice (AC), completely describe the cardinalities of all possible infinite sets [35,39]. Therefore, it made sense to conjecture that $\boldsymbol{c}$ is in fact $\aleph_1$, i.e. the next biggest cardinal after $\aleph_0$, and that $2^{\aleph_0} = \aleph_1$, as stipulated by CH [7,35]. A natural extension to CH was to suggest that the same pattern holds for all cardinals $\aleph_\lambda$, i.e. $2^{\aleph_\lambda} = \aleph_{\lambda+}$ (where $\aleph_{\lambda+}$ denotes the next biggest cardinal after $\aleph_\lambda$), and this constitutes the Generalized Continuum Hypothesis (GCH) [35,56].

The relevance of CH to progress in many areas of pure mathematics was widely publicized in 1900 when the David Hilbert chose its proof as first on his list of the 23 most important unresolved problems challenging 20th Century mathematicians [5,40]. For Hilbert, who gave an enthusiastic endorsement to Cantor's theory of transfinite numbers [5,39,40], the resolution of CH (alongside the Well-Order Theorem, WOT) was central to many of the underpinnings of mathematical theory in which set theory was already playing such a key role [39,40]. Although an answer to WOT came as early as 1904 at the hands of Zermelo's Axiom of Choice (AC) [39,61], Hilbert could not have anticipated that, in future years, the solution to CH would be left unchallenged and, to a certain degree, unchallengeable. First, Gödel proved in 1936 the consistency of AC and GCH with the axioms of ZF, implying (because of Gödel's own incompleteness theorems) that neither AC nor GCH can be disproved from ZF [32,33,35]. Later, in 1963, Cohen completed Gödel's result by proving that the consistency of ZF is also conserved by *i*) the axiomatic addition of the negation of AC, or *ii*) the addition of AC and the negation of CH, implying that neither AC nor CH (and therefore GCH) can be proved from ZF [17,18,35]. These findings led both Gödel and Cohen to believe that CH is false and suspect that the cardinality of the continuum ($\boldsymbol{c}$) is even greater than the second aleph ($\aleph_1$), an opinion that is largely shared by modern set theorists [19,34,39].

The integration of CH into the axiomatic fabric of set theory generates technical difficulties which have not yet been fully resolved [39]. The Cantorian treatment of infinity was ultimately responsible for the extension of set theory beyond its more basic (naïve) formulation, needed to resolve many of the paradoxes that famously troubled the discipline in its early days [35,52,56]. However, the relative success of such an extension is still tarnished by the problems encountered in dealing with AC [39].

Most of the current complexity of set theory stems, therefore, from the need to incorporate transfinite numbers into the fabric of mathematics. The implication is that the need for complexity could readily disappear if the arguments supporting transfinite theory were found to be flawed. For well over a century, Cantor's transfinite edifice has rested heavily on three key results, namely his first proof for the nondenumerability of $\mathbb{R}$ [6], the diagonalization argument [8], and the set-theoretical proof of Cantor's Theorem [62]. These proofs rely on *reductio* arguments which, as far as the author is aware, have not been superseded by later results.

## 3. Refutations of Cantor's Proofs of Nondenumerability

None of the *reductio* arguments that Cantor used to prove the uncountable nature of $\mathbb{R}$ and $\mathcal{P}(\mathbb{N})$ can be considered technically demanding. In all cases, the proofs do not extend beyond a few lines and a short number of mathematical statements. This is particularly true of Cantor's diagonalization argument [8,28]. Such a simplicity, combined with the counter-intuitive nature of the conclusion, may be the reason why the diagonalization argument has been repeatedly questioned, though to date unsuccessfully [43]. Inevitably, any new attempts to refute Cantor's argument are likely to be treated with sceptism.

Nevertheless, a number of new approaches that refute Cantor's proofs of nondenumerability are presented here, preceded by an initial discussion on the formulation of proofs by contradiction.

3.1. **Proofs by contradiction and mathematical consistency.** As an indirect proof, the method of *reductio* (*ad absurdum*), i.e. proof by contradiction, was used by the ancient Greeks and formalised by Leibniz [22,52]. The great power of the *reductio* method lies in the fact that the truth of a mathematical statement $P$ can be established even when a proof for it has not been constructed. For the purpose of the analysis that follows, it is convenient to describe this method as a suitable extension to the construction of direct proofs.

3.1.1. *Direct proofs.* In principle, the truth of a given statement $P$ is established whenever it can be derived from another statement $A$, already known to be true, via the conditional $A \Rightarrow P$ [2,52]. By the application of conditional elimination, also known as *modus ponens* [37,52], the truth of $P$ will be a direct consequence of the truth of $A$. However, the construction of the proof will normally entail the incorporation of a string of interconnecting statements $Q_1, Q_2, \cdots, Q_n$ such that

$$A \Rightarrow Q_1 \Rightarrow Q_2 \Rightarrow \cdots \Rightarrow Q_n \Rightarrow P \,, \tag{3.1}$$

where the logical connectives between statements can be single conditionals '$\Rightarrow$' (i.e. 'if ... then'), as shown in (3.1), or biconditionals '$\Leftrightarrow$' (i.e. 'if and only if ... then', iff ), without affecting the proof. Provided that all intermediate statements $Q_1, Q_2, \cdots, Q_n$ are true, the rule of hypothetical syllogism [2,52] will apply to the sequence (3.1) to derive $A \Rightarrow P$, and the truth of $P$ will then follow.

3.1.2. *Proofs by external (or conventional) contradiction.* In the most common implementation of the *reductio* method, the proof works by first assuming the truth of the opposite statement (i.e. the negation of $P$, $\neg P$), and then allowing the implementation (as before) of standard rules of inference to proceed through a string of interconnecting statements $Q_1, Q_2, \cdots, Q_n$ until a final *contradiction* is reached (i.e. a statement which is always false [22,37,52]). It is commonplace to describe such a contradiction as a composite statement $(R \wedge \neg R)$, hence reinforcing the fact that a contradiction is a statement that can never be true [37]. Therefore,

$$\neg P \Rightarrow Q_1 \Rightarrow Q_2 \Rightarrow \cdots \Rightarrow Q_n \Rightarrow (R \wedge \neg R) \,. \tag{3.2}$$

As with direct proofs, the truth of the intermediate statements $Q_1, Q_2, \cdots, Q_n$ is a fundamental requirement of the proof, which makes it possible to implement the rule of hypothetical syllogism to (3.2) and derive

$$\neg P \Rightarrow (R \wedge \neg R)\ . \tag{3.3}$$

The contradiction $(R \wedge \neg R)$ is obviously false. Consequently, the rules of *modus tollens* and double negation [47,52] allow the following to be written

$$\neg(R \wedge \neg R) \Rightarrow \neg(\neg P) \Rightarrow P\ , \tag{3.4}$$

completing the proof.

It is important to stress that the sequence (3.2)-(3.4) does not equate to writing $\neg P \Rightarrow P$. In this type of proof, the truth of $\neg P$ does not directly convey the truth of $P$, in the way that the truth of $A$ implies the truth of $P$ in the direct proof represented by (3.1). Instead, such a conclusion is achieved on the basis of the falsehood of $(R \wedge \neg R)$ by resorting to (3.4). Furthermore, the falsehood of this eventual contradiction $(R \wedge \neg R)$ is completely independent of $\neg P$ being true or false.

The purpose of the preceding observation is not only to emphasize the primary difference between a direct proof and a proof by external contradiction, but also to consider the various situations in which a contradiction can be found. Clearly, for the same reasons that $(R \wedge \neg R)$ can never be true, the logical coexistence of $\neg P$ and $P$ (i.e. any scenario where these two statements are simultaneously true or false) would also have to be judged contradictory. Therefore, if the truth of $\neg P$ were to be legitimately associated with the truth of $P$, it is conceivable that a claim for proof could be made.

3.1.3. *Proofs by internal (or self-referential) contradiction.* This supplementary classification can be introduced to account for constructions where the initial (assumed true) statement is $\neg P$ but, after a chain of interconnecting statements $Q_1, Q_2, \cdots, Q_n$, the concluding statement is $P$ itself, instead of the falsehood $(R \wedge \neg R)$ obtained with proofs by external contradiction, Hence,

$$\neg P \Rightarrow Q_1 \Rightarrow Q_2 \Rightarrow \cdots \Rightarrow Q_n \Rightarrow P\ . \tag{3.5}$$

Provided the statements $Q_1, Q_2, \cdots, Q_n$ are true, the rule of hypothetical syllogism can be applied to (3.5) to conclude that $\neg P \Rightarrow P$. The repeat rule [52] means that $\neg P \Rightarrow \neg P$ so, by conjunction introduction [52], (3.5) implies that

$$\neg P \Rightarrow (P \wedge \neg P)\ . \tag{3.6}$$

Since the contradiction $(P \wedge \neg P)$ is false, as previously applied to (3.3), *modus tollens* and double negation lead to

$$\neg(P \wedge \neg P) \Rightarrow \neg(\neg P) \Rightarrow P\ , \tag{3.7}$$

and the desired result, i.e. the truth of $P$, is finally concluded.

The sequence (3.5)-(3.7) replicates the logical steps used with (3.2)-(3.4). The only discrepancy is that the new contradiction $(P \wedge \neg P)$ incorporates the initial assumption $\neg P$ in a way that is not required by external (conventional) contradictions. But this might not constitute an impediment to the validity of the proof, provided that the truth of the intermediate statements $Q_1, Q_2, \cdots, Q_n$ in (3.5) is not compromised.

Nevertheless, additional complications can arise when the logical operators connecting the different statements in the proof are not only single conditional for every step, as depicted in (3.5), but also biconditional. The presence of biconditional implications is without effect on the validity of direct proofs (3.1), an observation that can be extended to proofs by external contradiction (3.2). The same cannot be said of proofs by internal (self-referential) contradiction. To evaluate such a scenario, consider a general alternative to (3.5) in which all the connectors from $\neg P$ to the intermediate statement $Q_i$ are biconditional instead of single conditional:

$$\neg P \Leftrightarrow Q_1 \Leftrightarrow Q_2 \Leftrightarrow \cdots \Leftrightarrow Q_{i-1} \Leftrightarrow Q_i \Rightarrow Q_{i+1} \Rightarrow \cdots \Rightarrow Q_n \Rightarrow P\ . \tag{3.8}$$

As an immediate consequence of (3.8), it will be the case for $Q_i$ (as well as for all its preceding intermediate statements) that

$$Q_i \Rightarrow \neg P \ \wedge \ Q_i \Rightarrow P\ , \tag{3.9}$$

which, as similarly deduced for (3.6) and (3.7), leads to

$$Q_i \Rightarrow (P \wedge \neg P)\ . \tag{3.10}$$

The contradiction (3.10) implies the falsehood of $Q_i$. But, as already explained, the validity of (3.8) requires the independent truth of all intermediate statements before hypothetical syllogism can be applied. The inevitable conclusion is that the presence of biconditional connectors immediately after the statement $\neg P$ in (3.8) invalidates such a proof. Furthermore, and because of the Law of Excluded Middle [22,52], it is clear that the falsehood of $Q_i$ implied by (3.10) is unaffected by which one of the two statements $P$ or $\neg P$ is actually true. Hence, the original purpose of the proof (i.e. to determine the truth of $P$) cannot be derived from the logical chain (3.8).

The intrinsic difficulties that constructions like (3.8) bring to the formulation of proofs by internal (self-referential) contradiction have consequences beyond the mere invalidity of (3.8). The assumption that the intermediate statement $Q_i$ is true (when the opposite is the case), and the acceptance of the truth of $P$ (when it might not be the case), could lead to mathematical inconsistency (see Section 1).

3.2. **Cantor's diagonalization argument.** When Georg Cantor wrote his famous diagonalization method of proof in 1891 [8] (an English translation of the original article can be found in [28]), his objective was greater than just corroborating the uncountability of $\mathbb{R}$, which he thought he had proved in 1874 [7]. As his own text indicates, Cantor's purpose was, not only to produce a much simpler argument, but also to provide a more general rationale that could "*be extended immediately to the general theorem that the powers of well-defined manifolds have no maximum, or, what is the same thing, that for any given manifold L we can produce a manifold M whose power is greater than that of L*" [8] [1]. This helps to explain, to a large extent, why the diagonalization argument was so warmly welcomed by many mathematicians, logicians as much as set-theorists [20,59].

[1] For manifold, read set. For power, read cardinality.

For this analysis, the presentation of the proof chosen here [35] is based on the decimal representations of real numbers [2]:

**Theorem (...)** $\mathbb{R}$ *is uncountable.*

***Proof.*** Suppose that $\mathbb{R}$ is countable. Then we can list the reals in the interval $[0,1)$ as

$$a_1, a_2, a_3, \cdots, a_n, \cdots$$

with each real in $[0,1)$ appearing as $a_n$ for exactly one $n \in \mathbb{N}$, $n \geq 1$. We shall represent each such real $r$ by its decimal expansion,

$$r = 0.r_1 r_2 r_3 \cdots r_n \cdots$$

avoiding the use of recurring 9s (so that e.g. we represent 0.2 by 0.2000…, rather than 0.1999…). We can then picture the reals in $[0,1)$ written out in an array:

$$\begin{aligned} a_1 &= 0.a_{1,1} a_{1,2} a_{1,3} \cdots \\ a_2 &= 0.a_{2,1} a_{2,2} a_{2,3} \cdots \\ a_3 &= 0.a_{3,1} a_{3,2} a_{3,3} \cdots \\ &\vdots \\ a_n &= 0.a_{n,1} a_{n,2} a_{n,3} \cdots \\ &\vdots \end{aligned}$$

Now define a real number $r = 0.r_1 r_2 r_3 \cdots r_n \cdots$ by

$$r_n = \begin{cases} 4, & \text{if } a_{n,n} \geq 6, \\ 7, & \text{if } a_{n,n} < 6. \end{cases}$$

Then $r$ belongs to $[0,1)$. However, $r$ has been constructed to disagree with each $a_n$ at the $n$th decimal place, so it cannot equal $a_n$ for any $n$. Thus $r$ does not appear in the list, contradicting that the list contains all reals in $[0,1)$.

Thus (...) , as $[0,1) \subseteq \mathbb{R}$, we have that $\mathbb{R}$ is uncountable. Q.E.D.

The choice made to avoid the use of recurring 9s in the decimal expansions is not more than a technicality to ensure that each real number in the interval $[0,1)$ has a unique representation in the array [2]. It also needs to be pointed out that the selection of values for $r_n$ is just one of many to guarantee that $r_n \neq a_{n,n}$ for any $n$, the key requirement of the proof. Furthermore, the enumeration of real numbers $a_n$ in the array can be formulated in any order, with no change in the validity of the argument.

The diagonalization argument can be refuted by analyzing the type of logical errors and inconsistencies described in Section 3.1.

3.2.1. *Logical objection to the proof.* Returning to the text of the proof used in Section 3.2, the chain of statements could be listed as follows:

- $P =$ 'The unit interval $I := [0, 1)$ is un uncountable set'
- $\neg P =$ 'The unit interval $I := [0, 1)$ is a countable set'
- $Q_1 =$ 'The set $I$ can be written out as an array of decimal expansions
$$a_n = 0.a_{n,1} a_{n,2} a_{n,3} \cdots a_{n,n} \cdots$$
where $n \in \mathbb{N}$, $n \geq 1$'

- $Q_2 =$ 'It is possible to define a real number $r = 0.r_1r_2r_3\cdots r_n\cdots$ by
$$r_n = \begin{cases} 4, & \text{if } a_{nn} \geq 6, \\ 7, & \text{if } a_{nn} < 6. \end{cases}$$
such that $r$ belongs to $[0, 1)$ but is $r \neq a_n$ for all $n \in \mathbb{N}$, so cannot be part of the array'
- $Q_3 =$ 'The array is not a complete listing of the elements of $I$'

This list forms the logical sequence

$$\neg P \Leftrightarrow Q_1 \Leftrightarrow Q_2 \Rightarrow Q_3 \Leftrightarrow P\ , \tag{3.11}$$

where all connectives are biconditional (iffs), apart from the single conditional between statements $Q_2$ and $Q_3$; in principle, there could be other reasons (not addressed by the proof) why the array is not a complete listing of $I$. But it is important to understand that the connective between statements $Q_1$ and $Q_2$ is biconditional: the real number $r$ can only be defined once the array of decimal expansions $a_n$ has been constructed and, in reverse, the definition of $r$ implies the existence of such an array [2].

Given the self-referential nature of Cantor's argument, the construction (3.11) mimics the construction (3.8) described in Section 3.1, so (3.11) is not a valid proof. Statement $Q_1$ is just a different wording for the assumption $\neg P$. Hence, the error must lie with the falsehood of statement $Q_2$, since $Q_2 \Rightarrow (P \wedge \neg P)$.

Statement $Q_2$ deals with the construction of the antidiagonal number $r$. By their own definition, the digits $r_n$ always differ from the digits $a_{nn}$ to be found along the diagonal in the array of decimal expansions $a_n$. So it is undoubtedly the case that $r$ cannot be found among the expansions $a_n$ used to construct $r$. However, statement $Q_2$ also assumes that, in order to construct $r$, it will always be possible to utilise *all* the decimal expansions in the array, provided the array is countable. Since the preceding analysis implies that $Q_2$ is false, the assumption just described appears to be the only possible cause for error, as it is unproven.

3.2.2. *Identification of the erroneous statement* ($Q_2$). As already indicated, the crux of the diagonalization argument lies in the claim that the real number $r = 0.r_1r_2r_3\cdots r_n\cdots$, defined by statement $Q_2$, cannot be found in the array of decimal expansions $a_n$, because *all* the expansions in the array are used for the construction of $r$. This statement, as formulated by Cantor, was assumed to be true without further verification [8]. However, a quantitative analysis can easily be undertaken to reveal the fundamental error of the assumption. In order to facilitate calculations, and rather than using the decimal expansions of Sections 3.2 and 3.2.1, it is preferable to work with binary strings, as originally envisaged by Cantor himself [8].

---

[2] This is not to say that the array described in the proof is the only one that can lead to the "antidiagonal" number $r$. Nevertheless, the construction of $r$ requires the existence of a *countable* array, regardless of the actual sequence of decimal expansions $a_n$, and/or the rules of substitution that generate the digits $r_n$. In this sense, "countable array $\Leftrightarrow$ antidiagonal number", with both statements *simultaneously* true or false, as required by the biconditional connective [2,37].

*Binary strings*. Consider the set of all infinite binary strings, i.e. the set of unending strings like ‘$10100101101101\cdots$’. This set can also be considered to represent all the real numbers in the interval $[0, 1)$, despite the duplications found with, for example, ‘$0.0011000\cdots$’ and ‘$0.0010111\cdots$’, where both binary representations account for the same real number [2]. One advantage of working with an array of binary strings is that, when applying the diagonalization argument to it, there is just one “antidiagonal” string (or real number) that can be obtained, since there is only one rule of substitution that can be applied to the construction of the string ‘$r_1 r_2 r_3 \cdots r_n \cdots$’ (or the real number $r = 0.r_1 r_2 r_3 \cdots r_n \cdots$), i.e. that every ‘0’ is replaced by ‘1’, and every ‘1’ by ‘0’.

With this in mind, the question to ask is not whether the string ‘$r_1 r_2 r_3 \cdots r_n \cdots$’ can be found (or not) in the original array but, accepting that the above string is indeed excluded, how many other binary strings are also excluded by the diagonal enumeration. First consider those strings that differ from the antidiagonal string ‘$r_1 r_2 r_3 \cdots r_n \cdots$’ in only one position, i.e. ‘$\boldsymbol{a_{1,1}} r_2 r_3 \cdots r_n \cdots$’, ‘$r_1 \boldsymbol{a_{2,2}} r_3 \cdots r_n \cdots$’, ‘$r_1 r_2 \boldsymbol{a_{3,3}} \cdots r_n \cdots$’, and so on. Clearly, if any of these strings were to appear at some point in the original array, they could only do it in the position that corresponds to the digit $a_{n,n}$. Hence, the subset consisting of such strings would then account for the whole array (it should also be obvious that this subset is countable) [3]. Next consider those strings that differ from ‘$r_1 r_2 r_3 \cdots r_n \cdots$’ in two digits – with all positions in the array already covered by the previous subset, they would have to be excluded from it, alongside ‘$r_1 r_2 r_3 \cdots r_n \cdots$’ itself. The same conclusion would be reached for the sets of strings that differ from the antidiagonal string ‘$r_1 r_2 r_3 \cdots r_n \cdots$’ in only three digits, four digits, five digits and so on ... all the way to the string that differs from ‘$r_1 r_2 r_3 \cdots r_n \cdots$’ in every digit, i.e. the diagonal string ‘$a_{1,1} a_{2,2} a_{3,3} \cdots a_{n,n} \cdots$’. An attempt to account for all the strings that would be excluded from the original array could be made by considering the cardinalities of the subsets described, which in each case can be equated (by means of appropriate bijections) with the unordered selections without repetition [3] of all positions along the strings taken in groups of none, (not one), two, three, four positions and so on, where such positions are those occupied by the diagonal digits $a_{n,n}$. Accordingly, we could write that the cardinality of the set of all excluded strings $S_E$ is given by

$$|S_E| = |C_{\mathbb{N},0} \cup C_{\mathbb{N},2} \cup C_{\mathbb{N},3} \cup \cdots \cup C_{\mathbb{N},k} \cup \cdots \cup C_{\mathbb{N},\mathbb{N}}|\,, \tag{3.12}$$

where each $C_{\mathbb{N},k}$ symbolizes the set of unordered selections without repetition of all the positions along the infinite strings (thus the suffix $\mathbb{N}$) taken in groups of $k$ elements. The total number of positions occupied by digits $a_{n,n}$ in each case can be finite ($k \in \mathbb{N}$) or infinite, therefore we allow for the “transition” from $C_{\mathbb{N},k}$ all the way to $C_{\mathbb{N},\mathbb{N}}$. Also, the expression (3.12) can be written without having to assert whether the (disjoint) sets $C_{\mathbb{N},k}$ are denumerable or not. To be noted is that $C_{\mathbb{N},1}$, i.e. the set corresponding to strings with only one diagonal digit $a_{n,n}$ present, is missing from (3.12), since these strings are not excluded from the original array. Obviously, $C_{\mathbb{N},0} \leftrightarrow \{`r_1 r_2 r_3 \cdots r_n \cdots'\}$, while $C_{\mathbb{N},\mathbb{N}} \leftrightarrow \{`a_{1,1} a_{2,2} a_{3,3} \cdots a_{n,n} \cdots'\}$.

---

[3] The point of this reasoning is not that the above subset of binary strings is the only one that can constitute the original array, since the antidiagonal string could derive from many other listings. The purpose of this argument is to evaluate how many strings are excluded from the diagonal enumeration.

Inspection of (3.12) shows that, as well as excluding the antidiagonal string '$r_1 r_2 r_3 \cdots r_n \cdots$', the diagonal enumeration of the original array readily excludes many more binary strings than had been anticipated. However, all sets $C_{\mathbb{N},k}$ where $k\,(\in \mathbb{N})$ is a finite number are certainly denumerable [56], despite being excluded. Finally, if the set $S_I$ of all binary strings that are included in the original array is considered, it becomes apparent that $S_I \leftrightarrow C_{\mathbb{N},1}$ and, consequently,

$$|S_I| = |C_{\mathbb{N},1}| \,. \tag{3.13}$$

To further this analysis, and in parallel to determining the size of $S_E$ by means of (3.12), an additional construction could be undertaken in order to generate the complete set of binary strings $S_T$. Consider an initial string formed with '0' as the only digit to appear in all positions, i.e. '$0000000 \cdots$'. Such a string can then be used as a *template* for the construction of every other string in $S_T$. This will be achieved by introducing rules of substitution where the digit '0' in '$0000000 \cdots$' is replaced by the digit '1' in a certain number of positions. When only one position per string is involved, the result will be strings with one '1' among a infinite sequence of 0's, like '$1000000 \cdots$' or '$0100000 \cdots$', for example. When two positions per string are involved, the resulting strings will have two 1's among the 0's, for example '$1010000 \cdots$' or '$0101000 \cdots$'. Similarly for three positions, four, five and so on, all the way to the substitution of all 0's by 1's, i.e. '$1111111 \cdots$'. In parallel to the computation of $S_E$, the sets of binary strings produced in this way can be equated (using bijections) with the unordered selections without repetition of all positions along the strings taken in groups of none, one, two, three, four positions and so on. Accordingly, the cardinality of $S_T$ will be determined by

$$|S_T| = |C_{\mathbb{N},0} \cup C_{\mathbb{N},1} \cup C_{\mathbb{N},2} \cup C_{\mathbb{N},3} \cup \cdots \cup C_{\mathbb{N},k} \cup \cdots \cup C_{\mathbb{N},\mathbb{N}}| \,, \tag{3.14}$$

where each $C_{\mathbb{N},k}$ is interpreted as in (3.12) and (3.13), albeit '1' always replacing $a_{n,n}$ and '0' doing the same for $r_n$.

When comparing the computation undertaken in (3.12) and (3.13) with that used to obtain (3.14), it should be clear that the only real difference lies in the starting *template* used, either the antidiagonal string '$r_1 r_2 r_3 \cdots r_n \cdots$' in the first case, or '$0000000 \cdots$' in the second. Consequently, it is legitimate to state that $S_T$ can be equated with the union of the disjoint sets $S_E$ and $S_I$, i.e. $S_T \leftrightarrow (S_E \cup S_I)$. Hence,

$$|S_T| = |S_E \cup S_I| \,. \tag{3.15}$$

The inevitable conclusion is that, whatever listing of binary strings might be considered as representation of $S_T$, Cantor's diagonalization argument always excludes a myriad of *countable* constituents of the array, where it should not. In fact, as shown here, the subset assimilated by the diagonal enumeration ($S_I$) amounts to little more than a minute proportion of all those elements of $S_T$ certain to be part of a *denumerable* infinity, in no way dependent on the countable or uncountable nature of $S_T$ as a whole.

The above analysis exposes the error embedded into the Cantorian assumption that an antidiagonal string (number) can always be constructed using *all* the elements present in a *countable* array or enumeration, thereby invalidating Cantor's proof.

*Quantification of the diagonal cover.* The argument described in the preceding section can be put on a more rigorous footing by a suitable definition of the scope that a diagonal enumeration has to account for (i.e. cover) the whole of a given array of binary strings [4], either finite or infinite.

**Definition 3.1 (Diagonal Cover).** *Consider an array (either finite or infinite in size) of binary strings of a given length (also finite or infinite). The proportion of elements in the array that a diagonal enumeration can account for is called the diagonal cover (Dc) of the array. A diagonal enumeration is understood as any counting of strings in the array that can be performed, along the whole length of the strings, by moving sequentially from one digital position in one string to another digital position in a different string (with neither position nor string already counted), until all of such positions or all of such strings have been accounted for.*

The above definition conceives diagonal enumerations in a general sense, but it is clear that all the possible counting protocols that abide by the definition will render the same value for the diagonal cover $Dc$. The most usual way to envisage these enumerations is along the diagonal line '$a_{1,1}a_{2,2}a_{3,3}\cdots a_{n,n}\cdots$'.

The value of $Dc$ is very simple to determine when dealing with finite arrays of binary strings of finite length: if the (finite) cardinality of the set is given by the integer $a$, and the total number of positions in the strings is determined by the integer $b$, the diagonal cover $Dc$ will be the fraction $Dc = b/a$ if $b < a$, or $Dc = 1$ whenever $b \geq a$. On the other hand, the value of $Dc$ for a finite array of strings of infinite length will always be $Dc = 1$, since the diagonal enumeration is only to cover a finite number of steps, as for the (finite) cardinality of the set, while the strings provide an unlimited number of positions.

The case of infinite arrays of binary strings of infinite length is obviously different, since the value of $Dc$ cannot be determined directly as the ratio of two integers. Nevertheless, as shown in what follows, $Dc$ can be extracted from a suitable construction of the array and its strings, and the subsequent calculation of the correponding limit, in conjunction with a prerequisite to the implementation of Cantor's diagonalization argument.

**Lemma 3.2.** *Consider an infinite array of binary strings of infinite length. In order to claim that an antidiagonal string (constructed in base to a diagonal*

[4] Similar definitions can be formulated for strings of digits in other bases (e.g. base 10). Neither the application of the definition nor its implications will be affected.

*enumeration) cannot be part of the array, it is essential to demonstrate that the value of the diagonal cover (Dc) for that array is* $Dc = 1$.

*Proof.* To prove the above claim for the antidiagonal string, it is necessary to show that the corresponding diagonal enumeration can account for all the binary strings in the array. This is equivalent to establishing that the ratio between the steps in the diagonal enumeration and the total number of strings in the array takes exactly the value 1. According to Definition 3.1, such a ratio is the diagonal cover $Dc$. Therefore, the required proof is that $Dc = 1$. □

The significance of Lemma 3.2 is better understood by its application to suitable examples that directly illustrate its impact on the diagonalization argument.

**Example 3.3.** Consider the set of all infinite binary strings, $S_T$, as described in the preceding section. According to the diagonalization argument, the assumed denumerability of $S_T$ (depicted by the listing of an array with all the binary strings that constitute $S_T$) can be shown to be false by constructing an antidiagonal string that, while it differs from all the strings in the array in at least one position, still is an element of $S_T$. In order to make this claim, Lemma 3.2 establishes the need to demonstrate that $Dc = 1$. The binary strings that constitute $S_T$ are all the possible arrangements of (infinite) sequential positions, taken by the digits '0' and/or '1' with any degree of order or repetition. The positions in the strings constitute a denumerable sequence.

In order to determine the value of $Dc$ for $S_T$, consider first a set/array $S_k$ of binary strings of finite (and fixed) length,

$$a_n = \text{`}a_{n,1}a_{n,2}a_{n,3}\cdots a_{n,k}\text{'}\,, \tag{3.16}$$

where $n, k \in \mathbb{N}$. Thus, all the strings (3.16) in the array $S_k$ have a fixed number ($k$) of positions occupied by '0' or '1'. The total number of strings that make up $S_k$ is entirely determined by $k$, since $|S_k| = 2^k$. But the length of the strings is $k$ itself. Therefore, the diagonal cover is $Dc = k/2^k$. Consequently, for any set $S_k$ it will always be true that $Dc < 1$. The first few values of $k$ illustrate this point:

(3.17)

| $k=1$ | $k=2$ | $k=3$ | $k=4$ |
|---|---|---|---|
| **0** | **0** 0 | **0** 0 0 | **0** 0 0 0 |
| 1 | 0 **1** | 0 **0** 1 | 0 **0** 0 1 |
| | 1 0 | 0 1 **0** | 0 0 **1** 0 |
| | 1 1 | 0 1 1 | 0 0 1 **1** |
| | | 1 0 0 | 0 1 0 0 |
| | | 1 0 1 | 0 1 0 1 |
| | | 1 1 0 | 0 1 1 0 |
| | | 1 1 1 | 0 1 1 1 |
| | | | 1 0 0 0 |
| | | | 1 0 0 1 |
| | | | 1 0 1 0 |
| | | | 1 0 1 1 |
| | | | 1 1 0 0 |
| | | | 1 1 0 1 |
| | | | 1 1 1 0 |
| | | | 1 1 1 1 |

Examination of (3.17) shows that the diagonal enumerations cannot account for the complete arrays. It can also be noticed that the antidiagonal strings ('1', '10', '111' and '1100', respectively) appear further down in the corresponding lists. As $k$ increases in value, the disparity between the length of the strings and the length of the arrays grows in exponential fashion.

To determine the value of $Dc$ for $S_T$, the only calculation required is to take the limit of $Dc(S_k)$ as $k$ grows without bound, i.e.

$$Dc(S_T) = \lim_{k\to\infty} Dc(S_k) = \lim_{k\to\infty} \frac{k}{2^k} = 0 \ . \tag{3.18}$$

The resulting value of $Dc(S_T)$ confirms the analysis of (3.12)–(3.15), i.e. the impossibility of the diagonal enumeration to account for even the smallest proportion of strings in the array.

Since Cantor's diagonalization argument attempts to prove that the set $S_T$ is uncountable, the fact that the diagonal cover is negligible in value could be simply seen as corroboration of the uncountable nature of the set, and not at all incompatible with the claim of the proof. However, it is still necessary for the diagonal cover $Dc$ of any denumerable set to be $Dc = 1$, as stated by Lemma 3.2. The simplest way to refute this scenario is with suitable counter-examples; two are described here, but others can be constructed.

**Example 3.4.** Consider a set of infinite binary strings, $S_{Q1}$, formed by strings that, after an initial finite number of positions occupied by either '0' or '1', revert to an infinite sequence of 1's. This is,

$$a_n = \text{`}a_{n,1}a_{n,2}a_{n,3}\cdots a_{n,k}\,1\,1\,1\,1\,1\,1\,1\cdots\text{'},\ k\in\mathbb{N} \ . \tag{3.19}$$

Given the construction (3.19) of the strings $a_n \in S_{Q1}$, it is obvious that $S_{Q1} \subseteq S_T$. More importantly, each string $a_n \in S_{Q1}$ can be seen as representing a rational number in the real interval $[0,1)$. Since $S_{Q1}$ can then be taken as representation of a subset of rational numbers, the conclusion is that $S_{Q1}$ is denumerable.

Given that $k$ only takes finite values, the diagonal cover for $S_{Q1}$ cannot be derived from the calculation of a limit, as in (3.18). Nevertheless, for every value of $k$, the total number of possible strings '$a_{n,1}a_{n,2}a_{n,3}\cdots a_{n,k}\,1\,1\,1\cdots$' is equal to $2^k$. This makes it possible to state that

$$Dc(S_{Q1}) < k/2^k\ ,\ \forall k\in\mathbb{N} \ . \tag{3.20}$$

The implication is that $Dc(S_{Q1}) < 1$, although $S_{Q1}$ is a denumerable set.

**Example 3.5.** Consider a set of infinite binary strings, $S_2$, where all the positions of every sequence are occupied by the same digit, '0' or '1', with the exception of a single position per string, taken by the opposite digit — that is, a single '0' per string while all the other digits are 1's, or a single '1' while all the others are 0's:

$$
\begin{array}{ccccccccccccc}
\mathbf{1} & 0 & 0 & 0 & 0 & 0 & 0 & 0 & 0 & 0 & \cdot & \cdot & \cdot \\
0 & \mathbf{1} & 1 & 1 & 1 & 1 & 1 & 1 & 1 & 1 & \cdot & \cdot & \cdot \\
0 & 1 & \mathbf{0} & 0 & 0 & 0 & 0 & 0 & 0 & 0 & \cdot & \cdot & \cdot \\
1 & 0 & 1 & \mathbf{1} & 1 & 1 & 1 & 1 & 1 & 1 & \cdot & \cdot & \cdot \\
0 & 0 & 1 & 0 & \mathbf{0} & 0 & 0 & 0 & 0 & 0 & \cdot & \cdot & \cdot \\
1 & 1 & 0 & 1 & 1 & \mathbf{1} & 1 & 1 & 1 & 1 & \cdot & \cdot & \cdot \\
0 & 0 & 0 & 1 & 0 & 0 & \mathbf{0} & 0 & 0 & 0 & \cdot & \cdot & \cdot \\
1 & 1 & 1 & 0 & 1 & 1 & 1 & \mathbf{1} & 1 & 1 & \cdot & \cdot & \cdot \\
0 & 0 & 0 & 0 & 1 & 0 & 0 & 0 & \mathbf{0} & 0 & \cdot & \cdot & \cdot \\
1 & 1 & 1 & 1 & 0 & 1 & 1 & 1 & 1 & \mathbf{1} & \cdot & \cdot & \cdot \\
0 & 0 & 0 & 0 & 0 & 1 & 0 & 0 & 0 & 0 & \cdot & \cdot & \cdot \\
1 & 1 & 1 & 1 & 1 & 0 & 1 & 1 & 1 & 1 & \cdot & \cdot & \cdot \\
0 & 0 & 0 & 0 & 0 & 0 & 1 & 0 & 0 & 0 & \cdot & \cdot & \cdot \\
1 & 1 & 1 & 1 & 1 & 1 & 0 & 1 & 1 & 1 & \cdot & \cdot & \cdot \\
0 & 0 & 0 & 0 & 0 & 0 & 0 & 1 & 0 & 0 & \cdot & \cdot & \cdot \\
1 & 1 & 1 & 1 & 1 & 1 & 1 & 0 & 1 & 1 & \cdot & \cdot & \cdot \\
0 & 0 & 0 & 0 & 0 & 0 & 0 & 0 & 1 & 0 & \cdot & \cdot & \cdot \\
1 & 1 & 1 & 1 & 1 & 1 & 1 & 1 & 0 & 1 & \cdot & \cdot & \cdot \\
0 & 0 & 0 & 0 & 0 & 0 & 0 & 0 & 0 & 1 & \cdot & \cdot & \cdot \\
1 & 1 & 1 & 1 & 1 & 1 & 1 & 1 & 1 & 0 & \cdot & \cdot & \cdot \\
\cdot & \cdot & \cdot & \cdot & \cdot & \cdot & \cdot & \cdot & \cdot & \cdot & \cdot & \cdot & \cdot \\
\cdot & \cdot & \cdot & \cdot & \cdot & \cdot & \cdot & \cdot & \cdot & \cdot & \cdot & \cdot & \cdot
\end{array}
\tag{3.21}
$$

It is clear that $S_2$ is a denumerable set given that each position in the strings accounts for only two possible sequences, one for either the digit '0' or '1' being the exception, as highlighted in (3.21). Regarding the diagonal cover, it is trivial to establish that

$$Dc(S_2) = \lim_{k \to \infty} \frac{k}{2k} = 1/2 \ , \tag{3.22}$$

hence, once more, $Dc(S_2) < 1$ while the set is denumerable.

Example 3.5 shows that the objections raised against the diagonalization argument are not restricted to sets constructed following exponential patterns, as in Examples 3.3 and 3.4, but can be applied in equal measure to all kind of denumerable sets. In summary, introducing the concept of diagonal cover $Dc$ (Definition 3.1) provides the means to articulate in a more quantitative fashion the arguments presented in the preceding sections against Cantor's method of proof.

*Interpreting bijections.* The problems identified in this article when dealing with Cantor's diagonalization argument raise the immediate prospect of a direct contradiction between ($i$) the criterion normally used to compare the cardinality of infinite sets, i.e. the existence (or not) of a bijection (one-to-one correspondence) between the sets [2,5,35], and ($ii$) the introduction of the concept of diagonal cover ($Dc$) to assess the ability (or not) of a diagonal enumeration to account for all the elements of an infinite set [5]. The standard interpretation given, for example, to a bijection between an infinite set $A$ and the set of natural numbers $\mathbb{N}$ is that, since every element of $A$ can be paired with a different element of $\mathbb{N}$ and, vice versa, every natural number can be linked to an element of $A$, then both sets $A$ and $\mathbb{N}$ have the same cardinality, i.e. $A$ is denumerable. On this basis, the Cantorian assumption that a diagonal enumeration simply mimics the bijection just described appears to be a reasonable idea that precludes the rationale behind the concept of diagonal cover.

[5] When applied to finite sets, it should be pointed out that both criteria reach the same conclusions, so no conflict or contradiction is encountered.

However, it is important to realise that the bijection between $A$ and $\mathbb{N}$, although it fulfills the purpose of demonstrating that $A$ is a denumerable set, does not reflect the actual process of construction needed to generate $A$ in the first place — a process that is always portrayed as the natural extension of a finite process. This much becomes evident when considering, for example, the set of binary strings $S_T$ described previously. The strings that constitute $S_T$ are conceived as an infinite (countable) sequence of digits '$a_{n,1}a_{n,2}a_{n,3}\cdots a_{n,n}\cdots$'; they are constructed as finite strings '$a_{n,1}a_{n,2}a_{n,3}\cdots a_{n,k}$', where the total number of positions in the strings (determined by $k$) is allowed to increase without bound. This is not only a convenient visualization of the strings, it is also the way in which the Axiom of Infinity constructs them and justifies their existence (see Section 3.4). Such a (countable) process of construction, by using $k$, also legitimizes the concept of diagonal cover ($Dc$; Definition 3.1).

These two approaches can be reconciled by understanding that, while the bijection between $A$ and $\mathbb{N}$ characterizes the former as a denumerable set, it does it in a qualitative fashion that does not acknowledge the relative size of $A$ (compared with another denumerable set). In this sense, the concept of denumerability becomes an absolute property assigned to the set, deprived of precise numerical meaning. On the other hand, and introducing a rationale that becomes complementary instead of contradictory, denumerable sets might also be compared in term of their relative sizes (as exemplified by the concept of diagonal cover). This second criterion adds a new perspective to the study of denumerable sets and provides a secondary degree of quantification, albeit in relative terms. Section 6.3 expands on this point and introduces the concept of relative cardinality for countable sets (both finite and infinite), a parameter that is readily subject to calculation.

3.3. **Cantor's Theorem: the power set.** The fundamental pillar that provides support to Cantor's theory of transfinite numbers is the proof(s) for the cardinality of the power set: $|A| < |\mathcal{P}(A)|$, true for any set $A$, either finite or infinite. The axiomatic importance of this statement cannot be underestimated. As well as underpinning the uncountable nature of $\mathbb{R}$, clearly stated by Cantor [8], it also contributes the only viable justification for the existence of higher cardinalities beyond $\aleph_0$, the cardinality of $\mathbb{N}$. When considering power sets, $\mathcal{P}(\mathbb{N})$ as well as subsequent ones, the critical difficulty is that no piecemeal process can be envisaged for the construction of any infinite set supposed to be nondenumerable [19,39]. It is this impossibility that makes the formulation of the Power Set Axiom so important for set theory [39,56].

Historically, two proofs have been put forward to demonstrate the countable nature of $\mathcal{P}(\mathbb{N})$ and subsequent power sets. The first proof, the diagonalization argument, was originally presented by Cantor in 1891 [8]. Cantor's original paper used this argument to determine the cardinality of the power set when the initial set is countable, and when it is assumed to be nondenumerable. The second proof, more set-theoretical but equally simple, was published by Zermelo in 1908 [62].

3.3.1. *First proof for the power set* $\mathcal{P}(\mathbb{N})$. As previously mentioned in Section 3.2, Cantor proposed his diagonalization argument as a general proof of the existence of nondenumerable sets [8]. Compared with his original argument of 1874 [6], he also saw it as a simpler proof of the uncountable nature of $\mathbb{R}$ that did not require any consideration of irrational numbers. It is this format that is often described in the literature, since the connection between $\mathbb{R}$ and $\mathcal{P}(\mathbb{N})$ is clearly apparent [59].

It suffices to consider the set of infinite binary strings, i.e. the set of unending strings like '$10100101101101\cdots$', described previously. The rationale of the proof is the same as described in Section 3.2, where the set of reals in the interval $[0,1)$ are represented by their decimal expansion. However, the added advantage of using binary expansions is that each unending binary string can also be thought to characterize a specific subset of natural numbers, where '0' indicates the absence of a given $n \in \mathbb{N}$ in such subset, and '1' indicates its presence. Thus, for example, the string written above represents the subset $\{0,2,5,7,8,10,11,13,\cdots\}$. The collection of unending binary strings, where all possible arrangements of 0s and 1s are envisaged, can be interpreted as the complete set of all possible subsets of $\mathbb{N}$, i.e. the power set $\mathcal{P}(\mathbb{N})$. Consequently, the same refutations and arguments included in Sections 3.2.1 and 3.2.2 can be applied here.

3.3.2. *First proof for higher powers*. As discussed in Section 3.2, Cantor's objective when introducing the diagonalization argument was to extend his proof of nondenumerability to sets of any cardinality, particularly power sets constructed from infinite sets already supposed uncountable, starting with the set of the real numbers, $\mathbb{R}$. This would allow him to claim the existence of a hierarchy of ever larger transfinite numbers [8]. When the initial set (used to construct the corresponding power set) is already assumed to be nondenumerable, the format of the proof has to differ from the versions described so far, because no enumeration can be constructed before defining the diagonal number (or function) critical to the argument. This version of the proof, as published by Cantor (English translation found in [28]) is reproduced below:

> This proof is remarkable not only because of its great simplicity, but more importantly because the principle followed therein can be extended immediately to the general theorem that the powers of well-defined manifolds have no maximum, or, what is the same thing, that for any given manifold $L$ we can produce a manifold $M$ whose power is greater than that of $L$.
>
> Let, for instance, $L$ be a linear continuum, say the totality of all real numbers which are $\geq 0$ and $\leq 1$.

> Let $M$ be the totality of all single-valued functions $f(x)$ which take only the values 0 or 1, while $x$ runs through all real values which are $\geq 0$ and $\leq 1$.
>
> That $M$ does *not* have a *smaller* power than $L$ follows from the fact that subsets of $M$ can be given which have the same power as $L$ – for instance, the subset which consists of all functions of $x$ which have the value 1 for a single value $x_0$ of $x$, and for all other values of $x$ have the value 0.
>
> But $M$ does *not* have *the same* power as $L$: for otherwise the manifold $M$ could be brought into a reciprocal one-to-one correspondence with the variable $z$, and $M$ could be thought of in the form of a single-valued function of the two variables $x$ and $z$
>
> $$\phi(x, z)$$
>
> such that to every value of $z$ there corresponds an element $f(x) = \phi(x, z)$ of $M$, and, conversely, to every element $f(x)$ of $M$ there corresponds a single determinate value of $z$ such that $f(x) = \phi(x, z)$. But this leads to a contradiction. For if one understands by $g(x)$ the single-valued function of $x$ which takes on only the values 0 and 1 and is different from $\phi(x, x)$ for every value of $x$, then on the one hand $g(x)$ is an element of $M$, and on the other hand $g(x)$ cannot arise from any value $z = z_0$ of $\phi(x, z)$, because is different from $g(z_0)$.
>
> But if the power of $M$ is neither smaller than nor equal to that of $L$, it follows that it is greater than the power of $L$. Q.E.D.

Although it is not explicitly stated in Cantor's presentation, the set ('manifold') $M$ can be shown to have the same cardinality as the power set of $L$: the single-valued functions $f(x)$ take only the values 0 or 1. As $x$ runs through all the real values included in $L$, each specific function $f(x)$ effectively maps a unique subset of $L$ and, *vice versa*, every subset of $L$ can be uniquely assigned to a single function $f(x)$ of set $M$. This implies that the cardinalities of the two sets, $m = |M|$ and $l = |L|$, are related to each other by the already familiar exponential $m = 2^l$, equivalent to writing $|M| = |\mathcal{P}(L)|$.

Again, the criticisms of this proof are similar to those detailed in Sections 3.2.1 and 3.2.2. First, the dissection of the proof according to logical connectives and statements shows how acceptance of its validity amounts to acceptance of inconsistency. The key chain of statements could be listed as follows:

- $P$ = '$M$ does not have the same power as $L$'
- $\neg P$ = '$M$ has the same power as $L$'
- $Q_1$ = '$M$ can be brought into a reciprocal one-to-one correspondence with $L$ by means of the single-valued function of two variables $\phi(x, z)$'
- $Q_2$ = 'It is possible to define the single-valued function $g(x)$, different from $\phi(x, x)$ for every value of $x$'
- $Q_3$ = 'The function $g(x)$ leads to contradiction because, although it is an element of $M$, equally it cannot be'

Such list forms the logical sequence

$$\neg P \Leftrightarrow Q_1 \Leftrightarrow Q_2 \Rightarrow Q_3 \Leftrightarrow P \tag{3.23}$$

where all connectives are biconditional (iffs), apart from the single conditional that links $Q_2$ and $Q_3$ (in principle, there might be other arguments that also lead to this contradiction). As with the analysis of Sections 3.2.1 and 3.2.2, the objections to the proof stem from the inconsistency that results from accepting the truth of $Q_2$, since the sequence (3.23) readily implies that $Q_2 \Rightarrow (P \wedge \neg P)$, which makes $Q_2$ false. It is also possible to identify the flaw in (3.22) as statement $Q_2$, here the definition of the function $g(x)$.

Because of the impossibility of using enumerations, Cantor made use of the single-valued function $\phi(x, z)$, where both variables $x$ and $z$ take the values of all the real numbers in the interval $[0, 1]$ defined by $L$. While this is a precise way of constructing the bijection between $M$ and $L$, the problem with the proof arises when it is claimed that a diagonal function $g(x) \neq \phi(x, x)$ can be constructed that differs from all the functions $f(x)$ included in $M$.

The claim can be refuted by considering the construction of the functions $f(x)$. The functions $f(x)$ can only take the values 0 or 1. If an initial value of the variable $x$ were considered, $x_1 \in L$, the set $M$ would be divided into two disjoint subsets, one subset containing those functions $f(x)$ that make $f(x_1) = 0$, the second subset containing those functions for which $f(x_1) = 1$. If a second value $x_2$ of $x$ were considered, each of the above subsets could now in turn be subdivided into two new subsets containing those functions for which $f(x_2) = 0$ or $f(x_2) = 1$. Hence, so far $M$ will have been divided into $2 \times 2 = 4$ disjoint subsets. The same process could be carried out a number of times, for a finite set of values of $x$, $\{x_1, x_2, \cdots, x_n\}$, resulting in a total of $2^n$ disjoint subsets of $M$, each defined by its unique combination of outputs (0s and 1s) for the range $f(x)_{(1,n)}$.

If an attempt is made to construct a diagonal function $g_n(x) \neq \phi(x, x)$, with the second variable $z$ taking values only from the same finite set $\{x_1, x_2, \cdots, x_n\}$, it is clear that, since the function $g_n(x)$ can only take a total of $n$ different values, $g_n(x)$ will be made to differ from the functions $f(x)$ described by only $n$ subsets of the total $2^n$ disjoint subsets constructed. In other words, $g_n(x)$ can only be claimed to disagree with a very small fraction $(n/2^n)$ of the total number of functions $f(x)$ included in $M$.

The construction of the functions $f(x)$ requires the progressive and complete subdivision of the subsets in the manner described for the whole range of real values of $x$ in the closed interval $[0, 1]$. Each function $f(x)$ is defined by a unique (infinite) combination of outputs (0s, 1s), as $x$ takes all the real values included in $L$. The construction of the final function $g(x)$ also requires that $z$ takes all the real values in $L$. Taking $l = |L|$, as with the analysis described in Section 3.2.2, $g(x)$ can only be claimed to differ from a final fraction $(l/2^l)$ of the totality of functions $f(x)$ included in $M$. Such a fraction stands as the limiting value of $n/2^n$ when $n \to l$ and, consequently, is infinitely small. Therefore, the contradiction claimed by statement $Q_3$ is untrue and the proof is shown to be invalid.

The refutation outlined here focuses on the factual errors of Cantor's argument. However, another observation merits attention and establishes an important connection with the more set-theoretical rationale used by Zermelo [62], discussed below (Section 3.3.3): the construction of the diagonal function $g(x) \neq \phi(x, x)$ is largely self-referential. The proof relies on a mathematical statement that is intrinsically contradictory.

3.3.3. *Second proof*. The argument for the nondenumerability of $\mathcal{P}(\mathbb{N})$, which Zermelo integrated into his original axiomatization of set theory in 1908 [62], is a particular case of a general statement regarding the relative cardinality of any set (finite or infinite) and its power set. By making a direct appeal to set-theoretical concepts, the proof diverges marginally from the diagonalization argument used by Cantor, but it still has a remarkably simple and concise formulation, as seen below [2]:

> (...) **Cantor's Theorem** *If $A$ is any set, then there is no surjection of $A$ onto the set $\mathcal{P}(A)$ of all subsets of $A$.*
>
> ***Proof.*** Suppose that $\varphi: A \to \mathcal{P}(A)$ is a surjection. Since $\varphi(a)$ is a subset of $A$, either $a$ belongs to $\varphi(a)$ or it does not belong to this set. We let
> $$D := \{a \in A : a \notin \varphi(a)\}.$$
> Since $D$ is a subset of $A$, if $\varphi$ is a surjection, then $D = \varphi(a_0)$ for some $a_0 \in A$.
> We must have either $a_0 \in D$ or $a_0 \notin D$. If $a_0 \in D$, then since $D = \varphi(a_0)$, we must have $a_0 \in \varphi(a_0)$, contrary to the definition of $D$. Similarly, if $a_0 \notin D$, then $a_0 \notin \varphi(a_0)$ so $a_0 \in D$, which is also a contradiction. Q.E.D.

Since a bijection is both injective and surjective [2,49], it is sufficient to prove that a surjection can never be constructed from $A$ onto $\mathcal{P}(A)$, in order to conclude that a one-to-one correspondence is impossible between a set and its power set. Consequently, the power set always has a greater cardinality, i.e. $|A| < |\mathcal{P}(A)|$. Cantor's Theorem is obviously true when the set considered is finite: if the cardinality of $A$ is $|A| = n$, the cardinality of the power set $\mathcal{P}(A)$ will be $|\mathcal{P}(A)| = 2^n$ and, trivially, $n < 2^n$ for all $n \in \mathbb{N}$, so no surjection is possible from $A$ onto $\mathcal{P}(A)$. Given that no additional evidence is required for finite sets, the proof reproduced above only becomes relevant to the case of infinite sets.

The second proof of Cantor's Theorem initially appears perfectly sound; its simplicity makes it difficult to identify potential pitfalls. However, the dissection of the logical structure of the proof, as undertaken in Sections 3.1, 3.2.1 and 3.3.2, raises doubts about the rigorous implementation of the *reductio* method. On this occasion, the short chain of statements can be described as follows:

- $P$ = 'For any set $A$, there is no surjection of $A$ onto its power set $\mathcal{P}(A)$'
- $\neg P$ = 'For any set $A$, a surjection $\varphi$ exists onto its power set $\mathcal{P}(A)$'
- $Q_1$ = 'There is a subset $D$ of $A$, defined as $D := \{a \in A : a \notin \varphi(a)\}$'
- $C$ = 'The existence of $D$ leads to contradiction: since $D = \varphi(a_0)$ for some $a_0 \in A$, if $a_0 \in D$ then $a_0 \notin D$ and, *vice versa*, if $a_0 \notin D$ then $a_0 \in D$' [6]

This list forms the logical sequence

$$\neg P \Leftrightarrow Q_1 \Leftrightarrow C \tag{3.24}$$

where both connectives are biconditional (iffs). Consequently, (3.24) equates to writing $\neg P \Leftrightarrow C$, and this is taken as a sufficient argument supporting the truth of the theorem. Therefore, the proof depends entirely on the valid definition of subset $D$, despite its self-referential nature [56].

[6] $C$ is a composite statement that could be dissected further into its more elementary parts. However, the main objective of this analysis is to highlight the logical rationale of the proof, hence using $C$ to encapsulate the contradiction serves the purpose.

The formulation of subset $D$ appears entirely (and exclusively) linked to the original assumption, i.e. the existence of surjection $\varphi$. Since the proof claims that subset $D$ (which leads directly to contradiction) is firm evidence of the incorrect formulation of surjection $\varphi$, it is assumed that the truth/falsehood of subset $D$ and surjection $\varphi$ are intrinsically linked. While this might be seen as a natural conclusion, it overlooks the distinct alternative of subset $D$ being a logical impossibility (because of its self-referential nature), despite surjection $\varphi$ being true. Due to this possible scenario, the proof should not be considered sufficiently rigorous and, up to this point, the true/false nature of surjection $\varphi$ should remain undetermined.

The incorporation of a self-referential statement into a mathematical proof is undesirable (see Section 5 for further discussion). The original attempts to establish Cantorian set theory as the foundational bedrock of mathematics were severely inconvenienced by the emergence of numerous paradoxes and antinomies, many arising from self-referential statements [49, 56]. The formulation of self-referential statements (leading to contradiction), within the boundaries of standard rules and definitions, was judged sufficiently troublesome to require the introduction of axiomatic principles specifically designed for their eradication [56]. In view of these past efforts, it is rather puzzling that Zermelo's proof has survived the test of time and is so widely recorded in modern textbooks [2,35,49,56]. It can only be assumed that the resilience of the proof has been due, more than to its own merits, to the perceived strength of the diagonalization argument as support of the same conclusion.

Although the analysis described in this section does not amount to a solid refutation, it raises serious (and more general) concerns about the inclusion of self-referential statements in the construction of mathematical proofs. The case of Cantor's Theorem and its second proof provides an opportunity to address such a general issue. An appropriate way to corroborate the invalidity of the formulation of subset $D$ (and properly refute the proof's argument) would be to provide a satisfactory proof for the existence of surjection $\varphi$. This is presented in Section 4.

3.4. **The Axiom of Infinity.** As mentioned in Section 1, the objective of this article is to undertake a systematic review of the arguments supporting transfinite number theory. The preceding sections have addressed Cantor's diagonalization argument and other closely linked results. In order to perform a similar exercise with Cantor's other important results, it is necessary to examine the way in which axiomatic set theory deals with the concept of infinity [35,56]. Particular attention has to be given to issues of consistency (or lack of it) in the interpretation of the axiom of infinity, which can formally stated as [35]

$$\exists y \Big( \varnothing \in y \wedge \forall x \big( x \in y \Rightarrow x \cup \{x\} \in y \big) \Big). \tag{3.25}$$

The axiom postulates that there is a set $y$ (known as an *inductive* set) such that the empty set $\varnothing$ is an element of $y$ and, whenever any $x$ is a member of $y$, the same applies to the union of $x$ with its singleton $\{x\}$.

The axiom of infinity cannot be derived from the other axioms of ZF, which cannot imply that there is an infinite set [35]. Therefore, (3.25) plays a crucial role in the handling of mathematical infinities, because it allows the definition of the set of natural numbers $\mathbb{N}$ as "the intersection of all inductive subsets of any inductive set $y$" [35]. From this, a natural number $\boldsymbol{n}$ is a member of $\mathbb{N}$, and $\mathbb{N}$ itself is an inductive set. The natural numbers can then be constructed using the *successor* function $x \to x^+$, where $x^+ = x \cup \{x\}$, identifying the empty set $\varnothing$ as the first natural number $\mathbf{0}\,(\in \mathbb{N})$, so that

$$\begin{aligned} \varnothing &\leftrightarrow \mathbf{0} \\ \varnothing \cup \{\varnothing\} = \{\varnothing\} = \varnothing^+ &\leftrightarrow \mathbf{1} = \mathbf{0}^+ = \mathbf{0} \cup \{\mathbf{0}\} = \{\mathbf{0}\} \\ \{\varnothing\} \cup \{\{\varnothing\}\} = \{\varnothing, \{\varnothing\}\} = \varnothing^{++} &\leftrightarrow \mathbf{2} = \mathbf{1}^+ = \mathbf{1} \cup \{\mathbf{1}\} = \{\mathbf{0}\} \cup \{\mathbf{1}\} = \{\mathbf{0},\mathbf{1}\} \\ \{\varnothing, \{\varnothing\}\} \cup \{\{\varnothing, \{\varnothing\}\}\} = \{\varnothing, \{\varnothing\}, \{\varnothing, \{\varnothing\}\}\} = \varnothing^{+++} &\leftrightarrow \mathbf{3} = \mathbf{2}^+ = \mathbf{2} \cup \{\mathbf{2}\} = \{\mathbf{0},\mathbf{1}\} \cup \{\mathbf{2}\} = \{\mathbf{0},\mathbf{1},\mathbf{2}\} \\ \cdots &\;\cdots\; \cdots \end{aligned} \tag{3.26}$$

From (3.26), it is interesting to observe that every natural number $\boldsymbol{n} \in \mathbb{N}$ is not just a member of $\mathbb{N}$, but also the set of all its predecesors thus containing, intuitively, $n$ elements [35]. $\mathbb{N}$ becomes the smallest set that contains $\mathbf{0}$ $(= \varnothing)$ and is closed under the successor function,

$$\mathbb{N} = \{\mathbf{0}, \mathbf{1}, \mathbf{2}, \mathbf{3}, \ldots, \boldsymbol{n}, \ldots\}. \tag{3.27}$$

Finally, the above construction provides a criterion that enables determination of whether a set is finite or infinite [35]:

> **Definitions**
>
> A set $X$ is *finite* if there is a bijection $f\colon n \to X$ for some $\boldsymbol{n} \in \mathbb{N}$. If there is no such bijection for any $\boldsymbol{n} \in \mathbb{N}$, $X$ is *infinite*.
>
> If there is a bijection $f\colon \mathbb{N} \to X$, then $X$ is *countably infinite* [7]. A *countable* set is one which is either finite or countably infinite.

A conventional view endorsed by the above definitions is that "a set is infinite if it is not finite" [2,35]. As discussed in the following section, this consideration translates into a viable and simple tool to determine when a given construction generates infinite set(s).

The formulation of the axiom of infinity presented in (3.25), and exemplified by (3.26) and (3.27), provides the backdrop from which to examine the concept of actual (completed) infinity [5,35].

3.4.1. *Constructing countable infinities.* Since other axioms of set theory cannot generate infinite sets, the correct interpretation of (3.25) becomes critical to the understanding of actual infinities. The starting point is the fact that the axiom of infinity is constructive — the inductive set $y$ exists because of the certainty that every member of the set exists, i.e. every element of $y$ can be constructed by applying the inductive process described by (3.25). This extends to the set of the natural numbers $\mathbb{N}$, i.e. $\mathbb{N}$ exists because the the successor function (starting from the empty set $\varnothing$) can be used to construct every single natural number $\boldsymbol{n}$ belonging to $\mathbb{N}$. To accept $\mathbb{N}$ as an actual (completed) infinity, it is necessary to acknowledge that every natural number $\boldsymbol{n}$ exists because it can be constructed. Therefore, completion means existence, unaffected by the fact that an infinite process can never be completed.

[7] For countably infinite, read denumerable.

As soon as the word "never" is introduced in the discussion of infinite sets (and how they are constructed), the dimension of time is forced upon the human perception of mathematical infinity, despite the fact that the algebraic and formal language used to write expressions like (3.25) does not operate on a temporal basis. The existence of every element of the inductive set $y$ (and, consequently, of every member of $\mathbb{N}$) does not depend on how long it would take to construct such an element. Instead, the inductive algorithm simply shows how the construction could be done — and this is sufficient to imply that elements (and set) exist.

In historic terms [5,35], the concept of actual infinity was seen as the only way to explain the existence of irrational numbers: from a set-theoretical perspective, real numbers are conceived as sets/sequences of rational numbers. If the real numbers under consideration are irrational, such sets/sequences are unavoidably infinite. In fact, nobody has ever managed to account for irrational numbers without making use of infinities [28,35].

Therefore, it is this need to accept the existence of actual infinities that acts as the key factor behind the interpretation of the axiom of infinity described in the preceding paragraphs. It is important that such an interpretation is applied consistently. The following theorem acts a suitable platform to enforce consistency in the construction of denumerable sets.

**Theorem 3.6 (of Actual Countable Infinity).** *The construction of a set A via the application of an inductive algorithm Φ is sufficient to establish the existence of all the members of A and characterise such a set A as countably infinite. This is,*

$$\exists a_0 \exists \Phi\Big(\Phi(a_0) = a_1 \wedge \forall n \in \mathbb{N}\big(\Phi(a_n) = a_{n+1}\big)\Big) \Rightarrow \exists A\left(A = \bigcup_{n=0}^{\infty}\{a_n\} \wedge |A| = |\mathbb{N}|\right). \tag{3.28}$$

*An inductive algorithm Φ is understood as any function or routine that assigns each element $a_n$ of A (in a one-to-one correspondence) to a different element of $\mathbb{N}$, based on an order determined by the construction of such elements. While set A has a least-element ($a_0$) that acts a the initial input of Φ but can never be the output of this algorithm, every other element $a_n \in A$ can be either input or output: as input, $a_n$ leads to the next element $a_{n+1}$ in the ordering; as output, $a_n$ is constructed by Φ using as input the preceding element $a_{n-1}$ in that ordering.*

*Proof.* By contradiction. First, consider a specific example of (3.28) where $a_0$ is the empty set $a_0 = \varnothing$, and the inductive algorithm $\Phi$ is the successor function, i.e. $\Phi(a_n) = a_{n+1} = a_n^{+} = a_n \cup \{a_n\}$. Accordingly, (3.28) becomes the construction of the natural numbers, as illustrated by (3.26) and (3.27). The definitions given before for finite and infinite sets can now be used to determine the actual nature of set $A$: if this set were finite, there would be a natural number $\boldsymbol{n} \in \mathbb{N}$ such that a bijection exists between $A$ and $\boldsymbol{n}$, $A \leftrightarrow \boldsymbol{n}$. And the implication would be that $A$ has a greatest element, $n_g$ (according to the ordering of the bijection) such that $n \leq n_g$, $\forall n \in A$. However, the definition of the inductive algorithm $\Phi$ means that every element $n \in A$ can be the input of $\Phi$, with output $\Phi(n) = n \cup \{n\} = n+1 > n$. When this is applied to $n_g$, the result is $\Phi(n_g) = n_g \cup \{n_g\} = n_g+1 > n_g$, contradicting the fact that $A$ cannot have any element greater than $n_g$. Hence, $A$ is not finite, i.e. $A$ is a countably infinite set, as stated by the theorem.

Furthermore, if $A$ is infinite, instead of a bijection $A \leftrightarrow \boldsymbol{n}$, the construction (3.28) leads to a one-to-one correspondence between $A$ and $\mathbb{N}$ itself. Given that, by virtue of (3.28) in the selected example, every element of $A$ is a natural number (i.e. an element of $\mathbb{N}$) and, *vice versa*, every element of $\mathbb{N}$ is constructed by (3.28) (i.e. also an element of $A$), we conclude that $A = \mathbb{N}$. According to the axiom of infinity (3.25), $\mathbb{N}$ is an infinite set. And, since all the elements of $\mathbb{N}$ (i.e. the natural numbers) are finite, the construction (3.26) guarantees their existence. The same conclusion applies then to the elements of $A$, not only for the selected example (when $A = \mathbb{N}$), but also for any other inductive algorithm $\Phi$ that satisfies the definition given by the theorem, since it is obvious that the bijection $A \leftrightarrow \mathbb{N}$ will exist in all these cases. □

Theorem 3.6 (of Actual Countable Infinity) addresses the need for consistency in the construction of infinite sets. Theorem 3.6 is of particular relevance to subroutines that, while satisfying the conditions imposed by the theorem on the inductive algorithm $\Phi$, form part of a much larger construction (see Section 4).

3.4.2. *The timeless existence and reordering of completed infinities.* As already mentioned, the construction of actual (completed) infinities is not a process that can be judged in terms of temporal execution. While it is clearly true that the mechanical completion of an infinite process cannot be achieved (for example, any inductive algorithm carried out by a computerised system), this time dependence is of no relevance to the mathematical description of such a process. The existence of a mathematical entity generated by an inductive algorithm (as in Theorem 3.6) is a reality unaffected by how long a human brain or any mechanical device would need to perform the calculations or constructions required. Instead, it suffices to have a valid description of how that algorithm would construct every element of a denumerable set. Such a description acts as a proof of existence. In other words, what the concept of infinite completion means is, not the temporal completion of the construction, but the certain existence of every element of the countably infinite set.

A significant consequence of Theorem 3.6 also results from the application of the theorem to the task of reordering. The infinite alteration of the order in which an enumeration presents the elements of a denumerable set is a process that can be conceived as a denumerable set of sequential transformations, making it again a completed infinity in which every transformation can be executed.

**Theorem 3.7.** *Let A be a denumerable set, and B an infinite subset of A, $B \subseteq A$. It will always be possible to construct a listing or enumeration of all the elements of A such that every element of B appears in that listing before any element of A that is not a member of B.*

*Proof.* Consider the denumerable set $A = \{a_0, a_1, a_2, a_3, \ldots, a_n, \ldots\}$, and an infinite string $s_0 =$ '$a_0\, a_1\, a_2 \ldots\, a_n \ldots$' that depictes all the elements of $A$ according to a given order $<_0$. Also consider an infinite subset of $A$, $B = \{b_0, b_1, b_2, \ldots, b_k, \ldots\}$, $B \subseteq A$. Thus, every element of $B$ will also be a member of $A$, i.e. $\forall b_k\, \exists a_n (a_n = b_k)$, and all the elements of $B$ will make one appearance in the infinite string $s_0$. The elements of $B$ have been written in a way that is compatible with the order used to construct $s_0$, i.e. $b_k <_0 b_{k+1}\ \forall b_k, b_{k+1} \in B$.

In line with Theorem 3.6, now conceive an inductive algorithm $\Phi$ that, using the string $s_0$ as initial input, undertakes a sequential process of reordering. $\Phi$ performs two operations: *(i)* first, $\Phi$ examines the first position of $s_0$ (i.e. $a_0$), in order to determine whether $a_0$ is also a member of $B$. *(ii)* Secondly, if $a_0$ is not an element of $B$, no operation takes place, since $a_0$ already occupies the first position in the string, and the output of this first step is $\Phi(a_0, s_0) = s_1$, where $s_1 = s_0$; however, if $a_0$ also belongs to $B$ (i.e. $a_1 = b_0$), then $\Phi$ relabels $a_0$ as $b_0$ and the resulting output is $s_1 =$ '$b_0 a_1 a_2 a_3 \ldots a_n \ldots$'. This first transformation can be named $t_0 = \Phi(a_0, s_0)$. In the next step, $\Phi$ uses the resulting string $s_1$ as input and performs the corresponding operations: *(i)* first, $\Phi$ examines the second position of $s_0$ (i.e. $a_1$), to establish whether $a_1$ is a member of $B$. *(ii)* Secondly, if $a_1$ is not part of $B$, no reordering is carried out and the output of the second step is $\Phi(a_1, s_1) = s_2$, where $s_2 = s_1$; but, if $a_1$ is a member of $B$ (i.e. $a_1$ will be $b_0$ or $b_1$, depending on the outcome of the first step), then $\Phi$ alters the ordering in $s_1$ by placing $b_0$ or $b_1$ $(= a_1)$ in front of $a_0$ or after $b_0$, respectively, i.e. $\Phi(a_1, s_1) = s_2$, where $s_2 =$ '$b_0 a_0 a_2 a_3 \ldots a_n \ldots$' or $s_2 =$ '$b_0\, b_1\, a_1\, a_2\, a_3 \ldots a_n \ldots$'. This second transformation is named $t_1 = \Phi(a_1, s_1)$.

The inductive algorithm $\Phi$ continues to use the output of each transformation as input for the next step of inspection and reordering (if required, and always immediately after the latest element of $B$ to have been relocated). In the way described, every time an element $a_n$ of $A$ is reordered, $a_n$ is also relabelled as the corresponding element $b_k$ of $B$. The sequential application of $\Phi$ generates a set $T$ of transformations $T = \{t_0, t_1, t_2, t_3, \ldots, t_n, \ldots\}$, where $t_n = \Phi(a_n, s_n) = s_{n+1}$. According to Theorem 3.6, $T$ is a denumerable set, with a bijection $T \leftrightarrow \mathbb{N}$ in place.

Consider $B$'s complementary set $C = \{c_0, c_1, c_2, \ldots, c_k, \ldots\}$, i.e. $C = A \setminus B$ and $A = B \cup C$ with $B \cap C = \varnothing$, such that $\forall c_k\ (c_k \in A \wedge c_k \notin B)$. A second consequence of Theorem 3.6, implied by the bijection $T \leftrightarrow \mathbb{N}$, is that the inductive algorithm $\Phi$ is capable of performing the operations of inspection and reordering for every element of $A$. Accordingly, a transformation $t_f \in T$ *exists* that has for output an infinite string $s_f$ where *all* the elements of $B$ are enumerated before any element of $C$, i.e. $s_f =$ '$b_0\, b_1\, b_2 \ldots b_n \ldots; c_0\, c_1\, c_2 \ldots c_n \ldots$'. Hence, string $s_f$ is the listing of all the elements of $A$ that completes the proof. □

Theorem 3.7 can be illustrated with a simple example based on the set of all the natural numbers, $\mathbb{N}$.

**Example 3.8.** Consider the set $\mathbb{N}$ of the natural numbers, listed according to the usual order $\leq$, $\mathbb{N} = \{0, 1, 2, 3, 4, \ldots, n, \ldots\}$. As a representative application of Theorem 3.7, take an algorithm $\Phi$ that inspects every natural number $n$ to determine whether $n$ is either even, i.e. $n \equiv 0 \pmod 2$, or odd, i.e. $n \equiv 1 \pmod 2$. If the number is odd, $\Phi$ leaves it where it is in the sequence. However, if $n$ is even, $\Phi$ places it immediately in front of '1'. In this way, the set of transformations $T$ will systematically reorder all even numbers in front of the odd ones. The sequence of strings $s_n$ will look as follows:

$$
\begin{array}{llllllllllllllll}
s_0 & \mathbf{0} & 1 & 2 & 3 & 4 & 5 & 6 & 7 & 8 & \cdot & \cdot & \cdot & n & \cdot & \cdot \; \cdot \\
s_1 & 0 & \mathbf{1} & 2 & 3 & 4 & 5 & 6 & 7 & 8 & \cdot & \cdot & \cdot & n & \cdot & \cdot \; \cdot \\
s_2 & 0 & \mathbf{2} & 1 & 3 & 4 & 5 & 6 & 7 & 8 & \cdot & \cdot & \cdot & n & \cdot & \cdot \; \cdot \\
s_3 & 0 & 2 & 1 & \mathbf{3} & 4 & 5 & 6 & 7 & 8 & \cdot & \cdot & \cdot & n & \cdot & \cdot \; \cdot \\
s_4 & 0 & 2 & \mathbf{4} & 1 & 3 & 5 & 6 & 7 & 8 & \cdot & \cdot & \cdot & n & \cdot & \cdot \; \cdot \\
s_5 & 0 & 2 & 4 & 1 & 3 & \mathbf{5} & 6 & 7 & 8 & \cdot & \cdot & \cdot & n & \cdot & \cdot \; \cdot \\
s_6 & 0 & 2 & 4 & \mathbf{6} & 1 & 3 & 5 & 7 & 8 & \cdot & \cdot & \cdot & n & \cdot & \cdot \; \cdot \\
s_7 & 0 & 2 & 4 & 6 & 1 & 3 & 5 & \mathbf{7} & 8 & \cdot & \cdot & \cdot & n & \cdot & \cdot \; \cdot \\
s_8 & 0 & 2 & 4 & 6 & \mathbf{8} & 1 & 3 & 5 & 7 & \cdot & \cdot & \cdot & n & \cdot & \cdot \; \cdot \\
\cdot & \cdot & \cdot & \cdot & \cdot & \cdot & \cdot & \cdot & \cdot & \cdot & \cdot & \cdot & \cdot & \cdot & \cdot & \cdot \; \cdot \\
s_{\mathrm{f}} & 0 & 2 & 4 & 6 & 8 & \cdot & \cdot & \cdot & ; & 1 & 3 & 5 & 7 & \cdot & \cdot \; \cdot
\end{array}
\tag{3.29}
$$

Within each string, the number highlighted in bold is the one inspected and reordered (when required) by the inductive algorithm $\Phi$.

It is interesting to observe that the existence of the final string $s_{\mathrm{f}}$ can be recognised independently of the set of transformations $T$, as a well-ordered series of all the elements of $\mathbb{N}$, formed by two fundamental segments, one for each of the subsets of all even or odd numbers [44]. Furthermore, since all strings $s_n$ are representations of the same set ($\mathbb{N}$), it should be legitimate to state that, as a further consequence of Theorem 3.7, bijections can be claimed to exists between any two strings, treated as sets. Selecting $s_0$ and $s_{\mathrm{f}}$ on that basis:

$$
\begin{array}{lllllllllllllll}
s_0 & 0 & 1 & 2 & 3 & 4 & \cdot & \cdot & \cdot & \cdot & \cdot & \cdot & \cdot & \cdot & \cdot \; \cdot \; \cdot \\
\updownarrow & \updownarrow & \updownarrow & \updownarrow & \updownarrow & \updownarrow & & & & & \updownarrow & \updownarrow & \updownarrow & \updownarrow & \updownarrow \\
s_{\mathrm{f}} & 0 & 2 & 4 & 6 & 8 & \cdot & \cdot & \cdot & ; & 1 & 3 & 5 & 7 & 9 \; \cdot \; \cdot \; \cdot
\end{array}
\tag{3.30}
$$

Although, according to (3.30), it is certain that the odd numbers in string $s_{\mathrm{f}}$ have counterparts in string $s_0$, it is obvious that the counterpart numbers cannot be specified, since the infinite enumeration of all the even numbers that precede them in $s_{\mathrm{f}}$ cannot be manually completed. However, those counterparts exist.

The interpretation of the axiom of infinity proposed here challenges the current orthodoxy. However, the formulation of Theorems 3.6 and 3.7 indicates that consistency in the treatment of actual (completed) infinities can be achieved by simply acknowledging that the existence of a denumerable set, like $\mathbb{N}$, implies the potential construction of every member of that set and, *vice versa*, the existence of all the elements of the set is essential for the existence of the set itself. This is something that current textbooks readily state regarding $\mathbb{N}$, for example, since every natural number $n$ is finite [2,35]. The conclusion derived from Theorem 3.6 is that the same can be said of any countably infinite (i.e. denumerable) set which forms part of a larger aggregate (see Section 4).

3.5. **Cantor's first proof of the nondenumerability of $\mathbb{R}$.** Cantor published his first proof for the uncountability of $\mathbb{R}$ (the set of real numbers) in 1874 [6], after repeated correspondence with Dedekind [5,28]. An English translation of the original article can be found in [28]. The relevant section is reproduced here:

> Suppose we have an infinite sequence of real numbers,
>
> $$\omega_1,\ \omega_2,\ \ldots\ \omega_\nu,\ \ldots \qquad (4)$$
>
> where the sequence is given according to any law and where the numbers are distinct from each other. Then in any given interval $(\alpha \ldots \beta)$ a number $\eta$ (and consequently infinitely many such numbers) can be determined such that it does not occur in the series (4); this shall now be proved.
>
> We go to the end of the interval $(\alpha \ldots \beta)$, which has been given to us arbitrarily and in which $\alpha < \beta$; the first two numbers of our sequence (4) which lie in the interior of this interval (with the exception of the boundaries), can be designated by $\alpha', \beta'$, letting $\alpha' < \beta'$; similarly let us designage the first two numbers of our sequence which lie in the interior of $(\alpha' \ldots \beta')$ by $\alpha'', \beta''$, and let $\alpha'' < \beta''$; and in the same way one constructs the next interval $(\alpha''' \ldots \beta''')$, and so on. Here therefore $\alpha', \alpha'' \ldots$ are by definition determinate numbers of our sequence (4), whose indices are continually increasing; the same goes for the sequence $\beta', \beta'' \ldots$; furthermore, the numbers $\alpha', \alpha'', \ldots$ are always increasing in size, while the numbers $\beta', \beta'', \ldots$ are always decreasing in size. Of the intervals $(\alpha \ldots \beta)$, $(\alpha' \ldots \beta')$, $(\alpha'' \ldots \beta'')$, ... each encloses all of those that follow.— Now here only two cases are conceivable.
>
> *Either* the number of intervals so formed is finite; in which case, let the last of them be $(\alpha^{(\nu)} \ldots \beta^{(\nu)})$. Since in its interior there can be at most one number of the sequence (4), a number $\eta$ can be chosen from this interval which is not contained in (4), thereby proving the theorem in this case.—
>
> *Or* the number of constructed intervals is infinite. Then the numbers $\alpha, \alpha', \alpha'', \ldots,$ because they are always increasing in size without growing into the infinite, have a determinate boundary value $\alpha^\infty$; the same holds for the numbers $\beta, \beta', \beta'', \ldots$ because they are always decreasing in size. Let their boundary value be $\beta^\infty$. If $\alpha^\infty = \beta^\infty$ (a case that constantly occurs with the set $(\omega)$ of all real algebraic numbers), then one easily persuades oneself, if one only looks back to the definition of the intervals, that the number $\eta = \alpha^\infty = \beta^\infty$ *cannot* be contained in our sequence;[1] but if $\alpha^\infty < \beta^\infty$ then every number $\eta$ in the interior of the interval $(\alpha^\infty \ldots \beta^\infty)$ or also on its boundaries satisfies the requirement that it not be contained in the sequence (4).—
>
> ---
>
> [1] If the number $\eta$ were contained in our sequence, then one would have $\eta = \omega_p$, where $p$ is a definite index. But this is not possible, for $\omega_p$ does not lie in the interior of the interval $(\alpha^{(p)} \ldots \beta^{(p)})$, while by definition the number $\eta$ does lie in the interior of the interval.

Cantor's presentation of the above proof has been substantially simplified and polished since 1874 [2]. The description of a particularly interesting version of this proof using Cantor's ternary sets can be found in [36]; this proof, which defines a so called *diagonal* function, is seen as an instance of a *diagonal argument* [56] and, as such, has conceptual affinity with the proofs of nondenumerability described (and refuted) in Section 3. As seen in what follows, these similarities can also be appreciated in the arguments used for refutation.

In all cases, the rationale of the proof is underpinned by the critically important Completeness Property of $\mathbb{R}$, i.e. "every nonempty set of real numbers that has an upper bound also has a *supremum* in $\mathbb{R}$" [2]. It is thanks to the completeness property that $\mathbb{R}$ can be considered a complete ordered field, with significant analytical implications. The Cantorian argument makes use of the well-known Nested Intervals Property. This states that, for every nested sequence of closed bounded intervals, there is at least a real number that belongs to every interval. Furthermore, if the lengths of the intervals constitute a set of real values that is bounded below by an *infimum* of value zero, $\inf = 0$, then the common real number belonging to every interval is also unique [2]. This could not be claimed without the completeness property.

As the text reproduced here shows, Cantor's own proof attempts to cover three possible scenarios: *(i)* a finite sequence of closed bounded intervals such that the common real values constitute an infinite set; *(ii)* an infinite sequence of intervals but the common real value is not unique; *(iii)* an infinite sequence of intervals with a common real value that is also unique. Given that the initial enumeration of the real numbers described by the proof is assumed to be complete and does not have any restrictions, it should be clear that only the third case merits consideration. Indeed, modern presentations consider only this third case [2].

If the real number line is visualized, the Cantorian construction of a sequence of closed intervals (where each pair of values $\alpha^{(\nu)}$ and $\beta^{(\nu)}$ are the endpoints that define each interval) can be represented as:

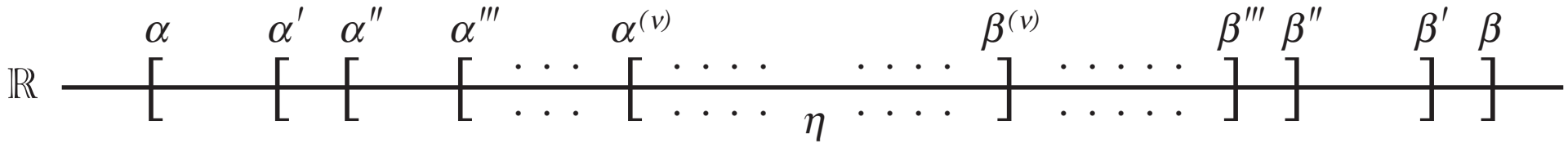

The unique common real value ($\eta$) becomes the supremum (infimum) of the two sequences, such that

$$\lim_{\nu\to\infty} \alpha^{(\nu)} = \lim_{\nu\to\infty} \beta^{(\nu)} = \alpha^{\infty} = \beta^{\infty} = \eta \,. \tag{3.31}$$

The endpoints of the nested intervals can be taken as part of two distinct sets, $A = \{\alpha, \alpha', \alpha'', \alpha''', \ldots, \alpha^{(\nu)}, \ldots\}$ and $B = \{\beta, \beta', \beta'', \beta''', \ldots, \beta^{(\nu)}, \ldots\}$, so $A, B \subseteq \mathbb{R}$ and $A \cap B = \varnothing$, such that $\forall \alpha_i \in A, \forall \beta_j \in B \; (\alpha_i < \beta_j)$. The crucial argument used by the proof is that, because of the completeness property, $\eta$ cannot belong to either set $A$ or $B$ and, accordingly, cannot be extracted from the initial enumeration of $\mathbb{R}$, despite being a real number. Nevertheless, a detailed examination of this claim quickly shows that the argument is flawed.

Based on the rationale described in Section 3.1 and applied to the previous proofs reviewed in this article, the text of Cantor's proof can be reduced to the following chain of statements:

- $P =$ '$\mathbb{R}$, the set of real numbers, is uncountable'
- $\neg P =$ '$\mathbb{R}$, the set of real numbers, is countable'
- $Q_1 =$ '$\mathbb{R}$ can be written as an infinite sequence
$$\omega_1, \omega_2, \omega_3, \ldots, \omega_\nu, \ldots,$$
where the sequence is given according to any law and where the numbers are distinct from each other'
- $Q_2 =$ 'By starting from any given interval $[\alpha, \beta]$, and orderly extracting from the above sequence upper and lower endpoints, it is possible to construct a nested sequence of closed bounded intervals, such that the endpoints of those intervals constitute two infinite sequences $A = \{\alpha, \alpha', \alpha'', \ldots, \alpha^{(\nu)}, \ldots\}$ and $B = \{\beta, \beta', \beta'', \ldots, \beta^{(\nu)}, \ldots\}$, where $\forall \alpha_i \in A, \forall \beta_j \in B\ (\alpha_i < \beta_j)$. The sequences $A$ and $B$ converge on an unique real number $\eta$ that, by virtue of the completeness property, can neither be part of the sequences $A$ or $B$, nor of the infinite sequence of real numbers described by $Q_1$'
- $Q_3 =$ 'The infinite sequence of real numbers described by $Q_1$ cannot be a complete listing of all the members of $\mathbb{R}$'

This list forms the logical sequence

$$\neg P \Leftrightarrow Q_1 \Leftrightarrow Q_2 \Rightarrow Q_3 \Leftrightarrow P\ , \tag{3.32}$$

where all connectives are biconditional (iffs), apart from the single conditional between statements $Q_2$ and $Q_3$ — *a priori*, it cannot be claimed that the argument used by the proof is the only one that might imply nondenumerability. Although the wording of statement $Q_2$ is lengthy, it only describes a single construction (the nested sequence of closed bounded intervals) and its intrinsic property (the unique common real number, $\eta$); both are a direct consequence of the initial assumption stated by $Q_1$.

It is important to understand why the connective between statements $Q_1$ and $Q_2$ is biconditional, $Q_1 \Leftrightarrow Q_2$. The construction of the nested sequence of intervals could not be performed without extracting their endpoints, $\alpha^{(\nu)}$ and $\beta^{(\nu)}$, from the assumed complete enumeration of $\mathbb{R}$. The specific listing of all the real numbers is not relevant to the argument, what matters is that such a listing exists and can be used. This rationale is similar to that offered in Footnote 2 for the diagonalization argument and reflects the formal similarities between both proofs.

As in previous sections, (3.32) is an example of an internal (self-referential) proof by contradiction. As before, (3.32) means that $Q_2 \Rightarrow (P \wedge \neg P)$. Hence, statement $Q_2$ cannot be true. However, this does not mean that the Completeness Property of $\mathbb{R}$ [2] needs to be questioned. Instead, everything points to the fact that the final assertion of statement $Q_2$, i.e. that $\eta$ cannot be found in the enumeration '$\omega_1, \omega_2, \omega_3, \ldots, \omega_\nu, \ldots$', is incorrect. The reason why this is the case is provided by an obvious corollary of Theorem 3.7: it is possible for *all* the elements of the subsets $A$ and $B$ or, if preferred, their union $(A \cup B) \subseteq \mathbb{R}$, to appear in the listing

'$\omega_1, \omega_2, \omega_3, \ldots, \omega_\nu, \ldots$' before other members of $\mathbb{R}$ that are not part of those subsets are encountered ... including the endpoint $\eta$.

It is of interest to examine how Cantor handles the exclusion of $\eta$ from the original listing in Footnote 1 of his proof. The original text clearly shows the circular character of Cantorian thinking. According to it, if $\eta$ were to appear in the initial listing of all elements of $\mathbb{R}$, then $\eta$ would also have to be a member of $A$ or $B$. However, by its own definition, $\eta$ cannot be a member of either set, which is claimed to be a contradiction. But, if $\eta$ would have been in either $A$ or $B$, then the converging endpoint of the two sequences would have been, not $\eta$, but a different number, say $\eta'$, so that $\eta' \neq \eta$. Once more, the same rationale would follow and a new and different number would emerge $\eta'' \neq \eta' \neq \eta$, *ad infinitum*. In this proof, as in his other arguments supporting nondenumerability, Cantor used self-referential contradictions which failed to avoid the pitfull of $Q_2 \Rightarrow (P \wedge \neg P)$.

It is of historic interest to observe that Cantor was not alone in following this potentially flawed rationale. Even Henry Poincaré, otherwise a very strong critic of Cantorian set theory [28]), described informally his own proof of the uncountable nature of $\mathbb{R}$ in a lecture dedicated to transfinite numbers which he delivered at Göttingen in 1909 (an English translation of the transcript can be found in [28]). Poincaré's argument, like Cantor's proof, fails when confronted by Theorem 3.7 (Section 3.4.2), and a similar refutation applies.

## 4. Proofs of Denumerability: $\mathcal{P}(\mathbb{N})$ and $\mathbb{R}$

The refutations included in Section 3 cast substantial doubts on the basic concept of nondenumerability, as originally presented by Cantor [6,8]; the construction of a convincing proof (or proofs) of denumerability will be sufficient to resolve this issue completely. That such proofs are impossible has been the universally shared view for many years, an obvious consequence of accepting Cantor's results. Certainly, Cohen's often quoted opinion was that "*... This point of view regards* [the continuum] *as an incredibly rich set given to us by one bold new axiom* [the Power Set Axion], *which can never be approached by any piecemeal process of construction*" [19,39]. While Cohen's belief went further by expecting the cardinality of the power set, i.e. $|\mathcal{P}(\mathbb{N})| = 2^{\aleph_0}$, to be larger than any transfinite aleph [19], common ground is kept by the conviction that no constructive proof of denumerability for $\mathcal{P}(\mathbb{N})$ (and, hence, $\mathbb{R}$) can ever be found.

But convictions can sometimes be changed. As shown below, once Cantor's proofs are discounted, the construction of direct proofs of denumerability for the power set $\mathcal{P}(\mathbb{N})$ is both possible and relatively simple. Theorem 3.6 (of Actual Countable Infinity) is the pivotal result that facilitates the process.

### 4.1. The denumerability of $\mathcal{P}(\mathbb{N})$.

As previously discussed in Section 3.3.3, the power set of any finite set always has a greater cardinality than the set itself. However, the situation is very different when the set considered is denumerable, i.e. countable and infinite.

Since what defines a set $A$ as denumerable is the construction of a bijection

(one-to-one correspondence) between $A$ and $\mathbb{N}$, it is obvious that a bijection will also exist between the power sets of $A$ and $\mathbb{N}$, i.e. between $\mathcal{P}(A)$ and $\mathcal{P}(\mathbb{N})$. Hence, determining the cardinality of $\mathcal{P}(\mathbb{N})$ will be sufficient to establish the cardinality of the power set of any denumerable set.

Before proceeding with the theorem and its proofs, a number of preliminary considerations are required. Firstly, it is important to realise that one way of constructing $\mathcal{P}(\mathbb{N})$ is to proceed through progressive stages of fragmentation into subsets of the final power set, according to parameters which can be shown to be countable. On this basis, the initial stage should be the classification of the subsets of $\mathbb{N}$ (the elements which constitute the power set), according to their cardinality. This can be illustrated by expressing the exponential term $2^n$ on the basis of the binomial expansion [1]

$$2^n = (1+1)^n = \binom{n}{0} + \binom{n}{1} + \cdots + \binom{n}{n-1} + \binom{n}{n} = \sum_{i=0}^{n} \binom{n}{i}, \tag{4.1}$$

where the binomial coefficients take special significance: if $2^n$ is understood as the cardinality of the power set for a finite set with $n$ elements and equates to the total possible subsets of the initial set, then every binomial coefficient in (4.1) accounts for the specific number of subsets with $i$ elements. Therefore, the summation of coefficients will replicate the enumeration of all possible subsets grouped according to their cardinality.

Similar fragmentation can be considered with $\mathcal{P}(\mathbb{N})$. The key consideration of such a classification is that $\mathcal{P}(\mathbb{N})$ can be subdivided into disjoint subsets, each comprised of subsets of $\mathbb{N}$ with a specific cardinality. Such a cardinality takes all possible values from $i = 0$ (the empty set, $\varnothing$) to the cardinality of $\mathbb{N}$ itself ($\aleph_0$). This listing is denumerable, in the same way as the listing of all the elements of $\mathbb{N}$ is denumerable [8].

A second observation, equally relevant to the construction of proofs for the denumerability of $\mathcal{P}(\mathbb{N})$, relates to the following theorem [2]:

> **Theorem** *If $A_m$ is a countable set for each $m \in \mathbb{N}$, then the union $A := \bigcup_{m=1}^{\infty} A_m$ is countable.*

This property of the union of countable sets is important for the construction of $\mathcal{P}(\mathbb{N})$ because it indicates that, once the enumeration based on cardinalities has been established, it is sufficient to show that each subset of $\mathcal{P}(\mathbb{N})$ (comprising of subsets of $\mathbb{N}$ with a given cardinality) can be constructed according to an equally countable enumeration. Such a procedure is reminiscent of the diagonal construction originally undertaken by Cantor to prove the denumerability of the set $\mathbb{Q}$ of all rational numbers [2,6,35]; in Cantor's proof, the mapping of the Cartesian product $\mathbb{N}\times\mathbb{N}$ onto $\mathbb{Q}^+$, the set of the positive rationals, is combined with the surjection of $\mathbb{N}$ onto $\mathbb{N}\times\mathbb{N}$ to compose a new surjection of $\mathbb{N}$ onto $\mathbb{Q}^+$, hence concluding that $\mathbb{Q}^+$ is a countable set (and so, by extension, is $\mathbb{Q}$) [2].

---

[8] It is important to understand that this enumeration covers the cardinalities of all possible subsets of $\mathbb{N}$. Since $\aleph_0$ is not an element of $\mathbb{N}$, it could be argued that such a listing cannot reach any of the infinite subsets of $\mathbb{N}$, nor $\mathbb{N}$ itself. But this objection fails to acknowledge that the construction of $\mathbb{N}$ is achieved by the same process of enumeration. By definition, the construction of infinite sets cannot be manually completed; it is precisely this property that makes them infinite and not finite. But, as extensively argued in Section 3.4.1 and enshrined in Theorem 3.6, accepting first that the listing (countable and without bound) of subsets of $\mathbb{N}$, grouped according to their cardinality, replicates the construction of $\mathbb{N}$ itself, the implication is that such an enumeration accounts for all subsets of $\mathbb{N}$, infinite as well as finite.

With these two observations in mind, the key result can now be presented.

**Theorem 4.1 (Denumerability of the Power Set).** *Let $\mathbb{N}$ be the set of all natural numbers, $\mathbb{N}=\{0, 1, 2, 3, \cdots, n, \cdots\}$. Its power set $\mathcal{P}(\mathbb{N})$, the set of all subsets of $\mathbb{N}$, is denumerable.*

Three independent proofs of the denumerability of $\mathcal{P}(\mathbb{N})$ are presented that are based on different ways of constructing the subsets of N: *(i)* Proof 1, based on a recursive algorithm; *(ii)* Proof 2, based on Cartesian products; *(iii)* Proof 3, based on exponential subsets. The direct proofs of Theorem 4.1 reported here are an immediate consequence of a rigorous interpretation of the axiom of infinity with regard to the construction of completed infinities (Section 3.4). They are also totally compatible with the refutations of Cantor's proofs of nondenumerability described in Section 3. It is legitimate to claim that both sets of results are, therefore, internally consistent.

4.1.1. *Proof 1 (Recursion).* The first step is to subdivide $\mathcal{P}(\mathbb{N})$ into a collection of disjoint subsets, $S_i$, each comprised of subsets of $\mathbb{N}$ with a given cardinality. Thus, it can be written that

$$\mathcal{P}(\mathbb{N})=\bigcup_{i=0}^{\infty} S_i \,, \tag{4.2}$$

where the subsets $S_i$ are defined as

$$S_i :=\left\{a_i : a_i \subseteq \mathbb{N}\,,\ |a_i| = i\ (i\in\mathbb{N}) \vee |a_i| = |\mathbb{N}|\ \right\}. \tag{4.3}$$

In order to satisfy (4.2), $i$ runs first through all the natural (finite) numbers $n\in\mathbb{N}$, and then extends to the subsets $a_i$ of infinite cardinality. The validity of this assertion depends on the ability that the algorithm used to construct the subsets of $\mathbb{N}$ has to deliver infinite subsets within a countable enumeration. As described in what follows, this requirement is fulfilled by the application of Theorem 3.6 where necessary. The next step is, therefore, to show that the subsets $S_i$ are also the result of a countable composition. This is achieved by subdividing each subset $S_i$ further into a collection of new subsets of $\mathcal{P}(\mathbb{N})$, initially denoted $A_i^g$. Each of these subsets is comprised of subsets of $\mathbb{N}$ ($a_i^g$), with the cardinality already assigned by (4.3), that are also characterised by the inclusion of a specific greatest element ($g\in\mathbb{N}$) in the subset. The smallest natural number $g$ that can be present in $a_i^g$ is $g=i-1$, as is the case for the subset $a_i^g=\{0, 1, 2, \cdots, i-1\}$. For each $S_i$ made up of finite subsets of $\mathbb{N}$, the enumeration of disjoint subsets $A_i^g$ progresses through all the natural numbers starting with $g=i-1$. A single exception to this rule is the first $S_i$ in the sequence, when $i=0$ and the only possible subset of $\mathbb{N}$ is $\varnothing$, for which no enumeration of $g$ and $A_i^g$ applies. Consequently, it can be written that

$$S_i = \bigcup_{g=i-1}^{\infty} A_i^g \,, \tag{4.4}$$

where the subsets $A_i^g$ (when $i \in \mathbb{N}$) are defined as

$$A_i^g := \left\{ a_i^g : a_i^g \subseteq \mathbb{N},\ n,i,g \in \mathbb{N},\ i \neq 0,\ |a_i^g| = i,\ n,g \in a_i^g,\ \forall n\,(n \leq g) \right\}. \tag{4.5}$$

As outlined by (4.4), the starting point of the enumerations ($g = i-1$) for the different subsets $S_i$ varies systematically with $i$. To construct the bijection of $\mathbb{N}$ onto $\mathcal{P}(\mathbb{N})$, a greater consistency can then be achieved by a convenient change of variable, which normalises the enumerations: since $i$ has a fixed value for each $S_i$, it is sufficient to introduce $j = g-i+1$ in (4.4), which then becomes

$$S_i = \bigcup_{j=0}^{\infty} A_i^{j+i-1} \,, \tag{4.6}$$

given that $g = j+i-1$. Therefore, the combination of (4.2) and (4.6) will eventually allow the construction of a surjection of $\mathbb{N} \times \mathbb{N}$ onto $\mathcal{P}(\mathbb{N})$, where each subset $A_i^g$ of $\mathcal{P}(\mathbb{N})$ is the unique representation of an ordered pair $(i, j)$ in the Cartesian product $\mathbb{N} \times \mathbb{N}$ [9]. As previously discussed, given that $\mathbb{N} \times \mathbb{N}$ is a denumerable set, this is sufficient to prove that $\mathcal{P}(\mathbb{N})$ is also denumerable, provided it can be shown that all the subsets $A_i^g$ are countable. This is first shown to be the case when $i \in \mathbb{N}$, i.e. when the subsets of $\mathbb{N}$ constructed by the algorithm are finite. Once this has been done, the proof needs to extend this construction to all infinite subsets of $\mathbb{N}$ (by applying Theorem 3.6).

The construction of each subset $A_i^g$ (when $i \in \mathbb{N}$) is a piecemeal process. In order to prove this fact, it is of value to consider an auxiliary type of subset of $\mathcal{P}(\mathbb{N})$, $B_i^g$, defined as

$$B_i^g := \left\{ b_i^g : b_i^g \subseteq \mathbb{N},\ n,i,g \in \mathbb{N},\ i \neq 0,\ |b_i^g| = i,\ n \in b_i^g,\ \forall n\,(n \leq g) \right\}. \tag{4.7}$$

Therefore, each subset $B_i^g$ is comprised of all those subsets $b_i^g$ of $\mathbb{N}$ with a given cardinality ($|b_i^g| = i$) that contain elements of $\mathbb{N}$ not greater than a certain value $g$. The absence of the condition $g \in b_i^g$ in (4.7) implies that $A_i^g \subseteq B_i^g$.

The introduction of the subsets $B_i^g$ provides the mechanism to construct each $A_i^g$ ($i \in \mathbb{N}$). Based on the definition described by (4.5), it can be seen that, since the subsets of $\mathbb{N}$ which make up $A_i^g$ always contain $g$, the total of such subsets ($a_i^g$) amounts to all the possible combinations of the members of $\mathbb{N}$ smaller than $g$ taken in sets of $i-1$ elements. This is,

$$|A_i^g| = \binom{g}{i-1} , \tag{4.8}$$

where $g \geq i-1$. Hence, each member $a_i^g$ of $A_i^g$ can be constructed by the union of the singleton $\{g\}$ with one of the ($i$-1)-element subsets of $\mathbb{N}$ comprised of numbers smaller than $g$. Furthermore, the definition (4.7) shows that such subsets constitute the members of the auxiliary subset $B_{i-1}^{g-1}$. Therefore,

---

[9] A bijection between $\mathbb{N} \times \mathbb{N}$ and $\mathcal{P}(\mathbb{N})$ is made impossible by the fact that no subsets $A_i^g$ are available for the ordered pairs $(0, j)$. As discussed in the main text, the construction of a surjection of $\mathbb{N} \times \mathbb{N}$ onto $\mathcal{P}(\mathbb{N})$ is still sufficient to support the overall proof [2].

$$a_i^g = \{g\} \cup b_{i-1}^{g-1} \ , \ \forall b_{i-1}^{g-1} \in B_{i-1}^{g-1} \ , \tag{4.9}$$

with $B_{i-1}^{g-1}$ defined as

$$B_{i-1}^{g-1} := \{ b_{i-1}^{g-1} : b_{i-1}^{g-1} \subseteq \mathbb{N},\ n,i,g \in \mathbb{N},\ i>1,\ |b_{i-1}^{g-1}| = i-1,\ n \in b_{i-1}^{g-1},\ \forall n \ (n \leq g-1) \}. \tag{4.10}$$

The condition $i>1$ is now necessary because, when $i=1$, the construction of $a_i^g$ is reduced to $\{g\}$ and no auxiliary subset $B_{i-1}^{g-1}$ is required. The analysis of (4.9) and (4.10) also shows that (4.8) can be extended to write

$$|A_i^g| = |B_{i-1}^{g-1}| = \binom{g}{i-1} . \tag{4.11}$$

The definition of $B_i^g$ given by (4.7) can now be compared with the definition of $S_i$ given by (4.3) for a finite cardinality. The similarity between the two definitions (all subsets $a_i^g$ and $b_i^g$ have the same cardinality) underlines the fact that the subsets $B_i^g$ are simply a special case of the subsets $S_i$ where all the subsets of $\mathbb{N}$ members of $B_i^g$ are restricted by not including any natural numbers greater than $g$. Consequently, taking into account (4.4), it is possible to write

$$B_i^g = \bigcup_{k=i-1}^{g} A_i^k , \tag{4.12}$$

where $i \neq 0$. Equation (4.12) can then be applied to the subset $B_{i-1}^{g-1}$ used in (4.9), and the result is

$$B_{i-1}^{g-1} = \bigcup_{k=i-2}^{g-1} A_{i-1}^k , \tag{4.13}$$

with the condition $i>1$. The relationship between subsets outlined by (4.13) is reinforced further by a well-known property of binomial coefficients [1] which, in combination with (4.11), leads to

$$|A_i^g| = |B_{i-1}^{g-1}| = \binom{g}{i-1} = \sum_{k=i-2}^{g-1} \binom{k}{i-2} = \sum_{k=i-2}^{g-1} |A_{i-1}^k| . \tag{4.14}$$

According to (4.9), the construction of $A_i^g$ requires the auxiliary subset $B_{i-1}^{g-1}$. And the construction of $B_{i-1}^{g-1}$ requires, according to (4.13), the construction of all the subsets $A_{i-1}^k$, where $k$ takes all the values from $i-2$ to $g-1$. Adapting (4.5), these subsets $A_{i-1}^k$ are defined as

$$A_{i-1}^k := \{ a_{i-1}^k : a_{i-1}^k \subseteq \mathbb{N},\ n,i,g,k \in \mathbb{N},\ i>1,\ i-2 \leq k \leq g-1, \\ |a_{i-1}^k| = i-1,\ n,k \in a_{i-1}^k,\ \forall n \, (n \leq k) \}. \tag{4.15}$$

Therefore, for each auxiliary subset $B_{i-1}^{g-1}$, a total of $g-i+2$ different subsets $A_{i-1}^k$ will have to be constructed. The construction of each $A_{i-1}^k$ can be achieved by adapting (4.9) such that

$$a_{i-1}^k = \{k\} \cup b_{i-2}^{k-1} \ , \ \forall b_{i-2}^{k-1} \in B_{i-2}^{k-1} \ , \tag{4.16}$$

with $B_{i-2}^{k-1}$ defined as

$$B_{i-2}^{k-1} := \{ b_{i-2}^{k-1} : b_{i-2}^{k-1} \subseteq \mathbb{N},\ n,i,k \in \mathbb{N},\ i>1,\ |b_{i-2}^{k-1}| = i-2,\ n \in b_{i-2}^{k-1},\ \forall n \ (n \leq k-1) \}. \tag{4.17}$$

According to (4.17), the cardinality of $b_{i-2}^{k-1}$, the subsets of $\mathbb{N}$ which constitute the auxiliary subset $B_{i-2}^{k-1}$, has been reduced by one unit ($i-2$).

The construction of each subset $B^{k-1}_{i-2}$ can be achieved by adapting (4.13), such that

$$B^{k-1}_{i-2} = \bigcup_{l=i-3}^{k-1} A^{l}_{i-2}\ , \tag{4.18}$$

which introduces a new layer of subsets of reduced cardinality $A^{l}_{i-2}$.

The construction of all the subsets $A^{l}_{i-2}$ requires a new round of definitions like (4.15)-(4.18), repeating the process and reducing the cardinality again by one unit. Such a process can be repeated as many times as necessary, until a final aggregate of singletons $A_1$ is reached.

The important comment is that the overall scheme outlined for the construction (up to this point in the proof) of the finite subsets $a^g_i$, members of $A^g_i$ ($i \in \mathbb{N}$), is essentially recursive [1,59]: for each subset $A^g_i$, the procedure consists of a series of iterations, with precisely defined operations, to be carried out until all relevant (finite) subsets $a^g_i$ of $\mathbb{N}$ are generated, a total determined by (4.8).

So far, and revisiting (4.3), the described construction can readily account for all the main subsets $S_i$ of $\mathcal{P}(\mathbb{N})$ that contain finite subsets of $\mathbb{N}$. If $\mathcal{F}(\mathbb{N})$ is the set of all finite subsets of $\mathbb{N}$, the conclusion is that $\mathcal{F}(\mathbb{N}) \subseteq \mathcal{P}(\mathbb{N})$ is a denumerable set, a result already known [56]. However, examination of the construction described by (4.2)-(4.18) in the light of Theorem 3.6 clearly shows that such a construction extends to subsets of infinite cardinality. Consider as an example – within the mapping of $\mathbb{N}\times\mathbb{N}$ onto $\mathcal{P}(\mathbb{N})$ that the construction leads to – the ordered pairs $(i,0)$ in $\mathbb{N}\times\mathbb{N}$ that map the subsets of $\mathbb{N}$ constituted by the first $i$ natural numbers, $a^g_i = \{0,\ 1,\ 2, \cdots,\ i-1\}$. As $i$ takes ever greater values, the subsets $a^g_i$ become increasingly similar to $\mathbb{N}$ itself. The increase in the value of the variable $i$ can be associated with an inductive algorithm $\Phi(i)$, as defined by Theorem 3.6, such that the set of ordered pairs $(i,0)$ in $\mathbb{N}\times\mathbb{N}$ forms a one-to-one correspondence with $\mathbb{N}$. Consequently, this bijection implies the construction by $\Phi(i)$ of an infinite set containing all the natural numbers: $\mathbb{N}$.

Similar sets of ordered pairs in $\mathbb{N}\times\mathbb{N}$ could be selected to apply Theorem 3.6 in the same way, always implying the construction of infinite subsets of $\mathbb{N}$. The inevitable conclusion is that the construction (4.2)-(4.18) can account for every possible infinite (as well as finite) subset of $\mathbb{N}$, i.e. the whole of $\mathcal{P}(\mathbb{N})$. □

Given the above proof of Theorem 4.1, it is appropriate to elaborate further on some of the determining aspects of Proof 1.

The construction of the subsets $A^g_i$ is shown to be a piecemeal process, the cardinality of each subset $A^g_i$ determined by (4.8). The complete enumeration of the members of $A^g_i$ (i.e. $a^g_i$, subsets of $\mathbb{N}$) can be integrated into the surjection of $\mathbb{N}\times\mathbb{N}$ onto $\mathcal{P}(\mathbb{N})$, such that each subset $a^g_i$ is the unique representation of an ordered pair in $\mathbb{N}\times\mathbb{N}$. This is achieved by expanding the variable $j = g - i + 1$, introduced in (4.6), into a new variable $j'$, designed to incorporate the enumeration of all the subsets $a^g_i$ which share the same greatest element $g$ defining each $A^g_i$. Hence, after rewriting (4.8) as

$$|A^g_i| = \binom{g}{i-1} = \binom{i+j-1}{i-1} = \binom{i+j-1}{j}\ , \tag{4.19}$$

$j'$ can be defined according to

$$j' := \begin{cases} \sum_{t=0}^{j-1} \binom{i+t-1}{i-1} + m, \ \forall m \left(0 \le m \le \binom{i+j-1}{j} - 1\right), & \text{if } i, j \ne 0, \\ 0 & , \text{ if } i \ne 0 \text{ and } j = 0. \end{cases} \tag{4.20}$$

The definition of $j'$ given by (4.20) facilitates the mapping of each subset $a_i^g$ by a unique ordered pair $(i, j')$ in $\mathbb{N}\times\mathbb{N}$. This mapping explains why proving that the subsets $A_i^g$ are countable is sufficient to complete the proof of Theorem 4.1 with regard to subsets $a_i$ of finite cardinality.

The introduction of the variable $j'$ is of help when describing the piecemeal construction of the subsets $A_i^g$. A key point is that the construction of the subsets $a_i^g$, for each $i$, is readily achieved by application of (4.9), needing only the list of subsets corresponding to the preceding value of $i$, that is, those subsets with cardinality just one unit lower. The starting enumeration is made of all the singletons $\{j'\}$, $j' \in \mathbb{N}$, which correspond to $i = 1$, i.e. the ordered pairs $(1, j')$ in the mapping. The next enumeration (for $i = 2$) is readily obtained by following the rule (4.9), so that the subsets with two elements are constructed. Similar construction follows for other subsets with higher cardinalities. Therefore, for the generation of all the subsets that are members of any $A_i^g$, the listing of subsets with the preceding cardinality is always available, since such an enumeration is required only up to the greatest element $g-1$. An algorithm designed to carry out this protocol will proceed uninterrupted and without bound, and be able to construct the surjection of $\mathbb{N}\times\mathbb{N}$ onto $\mathcal{P}(\mathbb{N})$ (infinite subsets included thanks to Theorem 3.6) without any need of external input.

A further observation on the countable nature of the subsets $A_i^g$, as the values of $i$ and $j$ are allowed to increase without bound, can be derived from (4.19). The binomial coefficient which accounts for the cardinality of these subsets can be expanded according to

$$|A_i^g| = \binom{i+j-1}{j} = \frac{i}{1} \times \frac{(i+1)}{2} \times \frac{(i+2)}{3} \times \cdots \times \frac{(i+j-1)}{j} . \tag{4.21}$$

(4.21) implies that the cardinality $|A_i^g|$ will still be finite while $i$ and $j$ remain confined to the domain of the natural (finite) numbers. However, if either (or both) variables are allowed to extend beyond their finite limits and into the domain of transfinite numbers, the calculation of the cardinality $|A_i^g|$ will be determined by the basic properties of transfinite arithmetic concerning the first transfinite cardinal ($\aleph_0$). Of such properties, it is sufficient to resort to those described by the three identities $\aleph_0 + n = \aleph_0$ $(n \in \mathbb{N})$, $\aleph_0 + \aleph_0 = \aleph_0$, and $\aleph_0 \times \aleph_0 = \aleph_0$ [35,49]. On this basis, it is apparent that the multiplication outlined by (4.21) can never take $|A_i^g|$ beyond the boundaries of denumerability. Since (4.2), (4.4) and (4.6) illustrate how $\mathcal{P}(\mathbb{N})$ can be broken down into a collection of disjoint subsets $A_i^g$, the above analysis is another way of arguing the case for Theorem 4.1, given that Theorem 3.6 allows the construction of $A_i^g$ to be extended to the domain of subsets $a_i$ with infinite cardinalities. Furthermore, the dissection of $\mathcal{P}(\mathbb{N})$ used in Proof 1 replicates the transformation of the exponential function $2^n$ into the summation of binomial coefficients described by (4.1). This confirms that the premise on transfinite exponentiation defended by Cantor [9,10,11] is incorrect.

4.1.2. *Proof 2 (Cartesian products).* One of the virtues of the algorithm used in the preceding proof of Theorem 4.1 is that every subset of $\mathbb{N}$ is constructed only once, without repetition. An alternative protocol can be used that, although generating repetitions in certain points, facilitates the required elimination of superfluous repeats by an appropriate mapping at each step. First, the subsets $S_i$ of $\mathcal{P}(\mathbb{N})$ introduced in (4.2) can be used to define a new aggregate when $i \in \mathbb{N}$:

$$T_k = \bigcup_{i=0}^{k} S_i \ , \tag{4.22}$$

such that

$$T_k := \left\{ a_i : a_i \subseteq \mathbb{N} \, , \ |a_i| = i \leq k, \ i,k \in \mathbb{N} \right\} . \tag{4.23}$$

It should be obvious that $T_k \subseteq \mathcal{P}(\mathbb{N})$, $\forall k \in \mathbb{N}$. Also, if the construction of $T_k$ was allowed to extend to the domain of infinite subsets of $\mathbb{N}$, then the result would be the power set $\mathcal{P}(\mathbb{N})$ itself.

As the title of the proof indicates, the construction of $T_k$ is based on Cartesian products. We start with the bijection between the set of natural numbers $\mathbb{N}$ and the Cartesian product of $\mathbb{N}$ with itself, i.e. $\mathbb{N} \times \mathbb{N} \leftrightarrow \mathbb{N}$ [2,35]. Each ordered pair $(i,j) \in \mathbb{N} \times \mathbb{N}$ will be assigned to a unique natural number $n \in \mathbb{N}$. The next step is a mapping $F_2$ such that:

$$\begin{aligned} F_2 : \mathbb{N} \times \mathbb{N} \ &\rightarrow \ T_2 \backslash \varnothing \\ (i,j) \ &\rightarrow \ \{i\} \cup \{j\} \ . \end{aligned} \tag{4.24}$$

It is clear that the surjection $F_2$ renders $T_2 \backslash \varnothing$ : when $i \neq k$, the union of the singletons $\{i\}$ and $\{j\}$ results in $\{i,j\}$; when $i = k$, the mapping will account for all subsets $a_i$ of cardinality $|a_i| = 1$. To complete $T_2$ we only need $T_2 = (T_2 \backslash \varnothing) \cup \varnothing$. Given the bijection $\mathbb{N} \times \mathbb{N} \leftrightarrow \mathbb{N}$ and the surjection $\mathbb{N} \times \mathbb{N} \rightarrow T_2 \backslash \varnothing$, we can conclude that $T_2$ is denumerable and a bijection exists $b_2 : T_2 \leftrightarrow \mathbb{N}$ [2]. The mapping $F_2$ eliminates all repetitions because ordered pairs like $(i, j)$ and $(j, i)$ render the same subset, i.e. $\{i,j\}$.

The next step is the Cartesian product $T_2 \times \mathbb{N}$ and the corresponding bijection $T_2 \times \mathbb{N} \leftrightarrow \mathbb{N}$. A new mapping (surjection) $F_3$ is required such that:

$$\begin{aligned} F_3 : \ T_2 \times \mathbb{N} \ &\rightarrow \ T_3 \backslash \varnothing \\ (a_i, j) \ &\rightarrow \ a_i \cup \{j\} \ , \end{aligned} \tag{4.25}$$

where $|a_i| \leq 2$. It is again easy to see that mapping $F_3$ renders $T_3 \backslash \varnothing$ since, when $j \in a_i$, the result is $a_i$ itself and, when $j \notin a_i$, the result is an increase in cardinality by one unit. The surjective mapping eliminates all repetitions. Since, as before, $T_3 = (T_3 \backslash \varnothing) \cup \varnothing$, we can conclude that there is a bijection $b_3 : T_3 \leftrightarrow \mathbb{N}$.

The construction proceeds as described, using the Cartesian products $T_{k\text{-}1} \times \mathbb{N}$ and the corresponding surjections $T_{k\text{-}1} \times \mathbb{N} \rightarrow T_k \backslash \varnothing$, thus rendering the bijections $b_k : T_k \leftrightarrow \mathbb{N}$. This protocol generates a set of bijections $B = \{b_2, b_3, b_4, \ldots, b_k, \ldots\}$, where $k$ takes the value of every natural number $k \geq 2$. Consequently, there is also a bijection $B \leftrightarrow \mathbb{N}$.

The algorithm described is certain to construct all subsets $a_i \subseteq \mathbb{N}$ with finite cardinality. To prove that the construction extends to infinite subsets $|a_i| = |\mathbb{N}|$, Theorem 3.6 can be applied. Consider, for example, an inductive algorithm $\Phi(k)$, as defined by the theorem, that extracts from every bijection $b_k$ the largest subset $a_k$ available ($|a_k| = k$) with the lowest possible value for the arithmetic sum of all its elements, i.e. $a_k = \{0, 1, 2, \ldots, k\text{-}1\}$. As in the proof of Theorem 3.6, the assignment of $a_k$ to each bijection $b_k$ implies the construction by $\Phi(k)$ of an infinite set containing all the natural numbers ($\mathbb{N}$).

Other selections of subsets $a_k$ assigned to the bijections $b_k$ could have been made, always implying the construction of infinite subsets of $\mathbb{N}$. Consequently, the final conclusion is that the construction rendering the set $B$ of bijections $b_k : T_k \leftrightarrow \mathbb{N}$ can account for every possible infinite (as well as finite) subset of $\mathbb{N}$, i.e. the construction of $T_k$, $\forall k \in \mathbb{N}$, extends eventually to $\mathcal{P}(\mathbb{N})$. □

The algorithms that Proof 1 (Recursion) and Proof 2 (Cartersian products) apply to construct the subsets of $\mathbb{N}$ illustrate the fact that there are multiple ways of "counting" all the elements of a denumerable set, i.e. multiple ways of constructing a bijection $A \leftrightarrow \mathbb{N}$ ($A$ countably infinite). It is also significant that Proofs 1 and 2 effectively use Theorem 3.6 an infinite number of times. Given that the proof of Theorem 3.6 (Section 3.4.1) is based on the distinction to be made between finite and infinite sets according to their own definition, it is interesting that neither Proof 1 nor Proof 2 take explicit advantage of the key characteristic of all infinite sets, i.e. that it is always possible to construct a bijection between the set itself and at least one of its proper subsets [5,26,27] (Section 2). By contrast, Proof 3 (Exponential subsets) uses a construction of $\mathcal{P}(\mathbb{N})$ that makes specific use of this defining property of infinite sets and, in doing so, simplifies both the overall algorithm and the application of Theorem 3.6.

4.1.3. *Proof 3 (Exponential subsets).* Consider $\mathbb{N}$, the set of natural numbers. Consider also, within the enumeration of the elements of $\mathbb{N}$, all the powers of 2:

$$\mathbb{N} = \{0,\ \mathbf{1},\ \mathbf{2},\ 3, \mathbf{4},\ 5,\ 6,\ 7,\ \mathbf{8},\ 9,\ 10,\ 11,\ 12,\ 13,\ 14,\ 15,\ \mathbf{16},\ 17,\ \cdots\}\ . \tag{4.26}$$

The powers of 2 can be selected to form the corresponding subset $N_2 \subseteq \mathbb{N}$:

$$N_2 = \{2^0,\ 2^1,\ 2^2,\ 2^3,\ 2^4,\ \cdots\} = \{2^k,\ \forall k \in \mathbb{N}\}\ . \tag{4.27}$$

It is evident that subset $N_2$ is an infinite set (same rationale as that used in the proof of Theorem 3.6). Hence, there is a bijection $N_2 \leftrightarrow \mathbb{N}$. The exponents of all the elements of $N_2$ can now be taken to construct another set:

$$N_k = \{0_k,\ 1_k,\ 2_k,\ 3_k,\ \cdots,\ n_k,\ \cdots\} = \{n_k,\ \forall n \in \mathbb{N}\}\ . \tag{4.28}$$

The only purpose of the subscript $k$ is to distinguish between the elements of $N_k$ and the elements of $\mathbb{N}$. It is obvious that $N_k \leftrightarrow N_2 \leftrightarrow \mathbb{N}$.

A partial listing of the elements of $N_k$ will render a finite subset $N_i \subseteq N_k$:

$$N_i = \{0_k,\ 1_k,\ 2_k,\ 3_k,\ \cdots,\ (i\text{-}1)_k\} = \{n_k \in N_i,\ n, i \in \mathbb{N},\ n < i\}\ . \tag{4.29}$$

The finite subset $N_i$ can then be used to construct its power set $\mathcal{P}(N_i)$:

(4.30) $$\mathcal{P}(N_i) = \{a_n : a_n \subseteq N_i\}\,,$$

such that $|\mathcal{P}(N_i)| = 2^i$.

Since, according to (4.28), $i$ can be take the value of every natural number, an inductive algorithm $\Phi(i)$ can be conceived in line with Theorem 3.6 that, taking as input the power set $\mathcal{P}(N_i)$, produces as output the power set $\mathcal{P}(N_{i+1})$ of the set $N_{i+1} = N_i \cup \{i_k\}$ :

(4.31) $$\mathcal{P}(N_{i+1}) = \{b_n : b_n \subseteq N_{i+1},\ (b_n = a_n \subseteq N_i) \vee (b_n = a_n \cup \{i_k\})\}\,,$$

such that $|\mathcal{P}(N_{i+1})| = 2 \times 2^i = 2^{i+1}$. At the same time, $\Phi(i)$ assigns each output, i.e. each power set $\mathcal{P}(N_i)$, to the corresponding integer $i \in \mathbb{N}$, thus constructing a bijection between the set $P$ of all the power sets $\mathcal{P}(N_i)$ and $\mathbb{N}$, $P \leftrightarrow \mathbb{N}$.

Examination of (4.26), (4.27) and (4.28) indicates that, when constructing each power set $\mathcal{P}(N_i)$, $\Phi(i)$ also establishes a one-to-one correspondence between the the members of $\mathcal{P}(N_i)$ and elements of $\mathbb{N}$. Any specific case can illustrate this fact, for example $\mathcal{P}(N_3)$:

(4.32)
$$\begin{array}{rcl@{\qquad\qquad}rcl}
\boldsymbol{a_n} & \leftrightarrow & \boldsymbol{n} & \boldsymbol{a_n} & \leftrightarrow & \boldsymbol{n} \\
\varnothing & \leftrightarrow & 0 & \{0_k, 1_k\} & \leftrightarrow & 4 \\
\{0_k\} & \leftrightarrow & 1 & \{0_k, 2_k\} & \leftrightarrow & 5 \\
\{1_k\} & \leftrightarrow & 1 & \{1_k, 2_k\} & \leftrightarrow & 6 \\
\{2_k\} & \leftrightarrow & 3 & \{0_k, 1_k, 2_k\} & \leftrightarrow & 7
\end{array}$$

Furthermore, the output of $\Phi(3)$ will be $\mathcal{P}(N_4)$ and, according to (4.31), the subsets of $\mathbb{N}$ to add to the list (4.32) will be:

(4.33)
$$\begin{array}{rcl@{\qquad\qquad}rcl}
\boldsymbol{b_n} & \leftrightarrow & \boldsymbol{n} & \boldsymbol{b_n} & \leftrightarrow & \boldsymbol{n} \\
\{3_k\} & \leftrightarrow & 8 & \{0_k, 1_k, 3_k\} & \leftrightarrow & 12 \\
\{0_k, 3_k\} & \leftrightarrow & 9 & \{0_k, 2_k, 3_k\} & \leftrightarrow & 13 \\
\{1_k, 3_k\} & \leftrightarrow & 10 & \{1_k, 2_k, 3_k\} & \leftrightarrow & 14 \\
\{2_k, 3_k\} & \leftrightarrow & 11 & \{0_k, 1_k, 2_k, 3_k\} & \leftrightarrow & 15
\end{array}$$

this is, the power set $\mathcal{P}(N_4)$ results from amalgamating (4.32) and (4.33). It should be observed that the natural numbers $n \in \mathbb{N}$ required to construct each bijection $\mathcal{P}(N_i) \leftrightarrow \boldsymbol{i}$ are precisely all those such that $n < 2^i$.

For each step, the inductive algorithm $\Phi(i)$ contructs a set with all possible subsets of $\mathbb{N}$ containing elements of $\mathbb{N}$ up to $i$-1. Theorem 3.6 implies that $\Phi(i)$ is capable to perform the same task to account for every element of $\mathbb{N}$, and that translates into the complete power set $\mathcal{P}(\mathbb{N})$. This conclusion is accompanied by the fact that, at every point, $\Phi(i)$ achieves its goal of constructing the bijection $\mathcal{P}(N_i) \leftrightarrow \boldsymbol{i}$ using only elements of $\mathbb{N}$ bounded by $n < 2^i$. When $\Phi(i)$ extends its output to account for infinite subsets of $\mathbb{N}$, this is therefore achieved without having to abandon the numerical domain set by $\mathbb{N}$ itself. Consequently, the final bijection is obtained, $\mathcal{P}(\mathbb{N}) \leftrightarrow \mathbb{N}$. □

4.2. **Corollaries**. Theorem 4.1 has far-reaching consequences. Some of the most immediate ones can be formalised in the following corollaries.

**Corollary 4.2.** *Let $A$ be a denumerable set. Its power set $\mathcal{P}(A)$, the set of all subsets of $A$, is also denumerable.*

*Proof.* The above statement is a trivial consequence of Theorem 4.1. A bijection exists between $A$ and $\mathbb{N}$, by definition. A bijection also exists between $\mathbb{N}$ and $\mathcal{P}(\mathbb{N})$, as a result of Theorem 4.1. On this basis, the first bijection can be used to replace every element of $\mathbb{N}$ by the corresponding element of $A$, and every subset of $\mathbb{N}$ by the corresponding subset of $A$. These changes transform the second bijection into a one-to-one correspondence between $A$ and $\mathcal{P}(A)$, as required. □

The above corollary can be applied to the power set of any denumerable set, including the power set of any power set already shown to be denumerable, e.g. $\mathcal{P}(\mathcal{P}(\mathbb{N}))$ or any subsequent power set. The unavoidable implication is that the Cantorian claim of existence for an infinite sequence of greater and greater transfinite cardinal numbers, alephs [8,23], which lies at the heart of both the Continuum Hypothesis (CH) and the Generalized Continuum Hypothesis (GCH), is left unsubstantiated.

**Corollary 4.3.** *The reals in the interval $[0,1)$ constitute a denumerable set.*

*Proof.* The connection between $I := [0,1) \subseteq \mathbb{R}$ and $\mathcal{P}(\mathbb{N})$, discussed in Section 2.4.1, is widely documented in the literature, and confirms the surjection of $\mathcal{P}(\mathbb{N})$ onto $I$ [35,49]. Such mapping combines with the denumerability of $\mathcal{P}(\mathbb{N})$ and, therefore, corroborates the countable nature of $I$. □

The denumerability of $I = [0,1)$ extends naturally to the totality of $\mathbb{R}$ [35,49], so it is a result of great importance, analytical as much as set-theoretical.

4.3. **The Continuum.** Sections 4.1 and 4.2 provide evidence on the countable nature of $\mathbb{R}$. Such results translate into a formal statement on the mathematical nature of the Continuum which is at odds with Cantor's claims. The constructions of the real numbers, independently developed by Dedekind and Cantor [23,35], led to the formulation of the Cantor-Dedekind Axiom, establishing that the points on a line can be aligned in a one-to-one correspondence with $\mathbb{R}$ [5,23]. Acceptance of the embodiment in the real number line of this intimate connection between arithmetics and geometry is a basic prerequisite of Theorem 4.4.

**Theorem 4.4 (of the Continuum).** *Seen as a numerical reality, the Continuum is a countable infinity. If $\aleph_0$ symbolizes the cardinality of $\mathbb{N}$, the set of the natural numbers, and $\boldsymbol{c}$ (the power of the Continuum) signifies the cardinality of the set of the real numbers $\mathbb{R}$, then*

$$\boldsymbol{c} = 2^{\aleph_0} = \aleph_0 . \tag{4.34}$$

*Proof.* See the proofs of Theorem 4.1, and Corollaries 4.2 and 4.3. □

The Theorem of the Continuum supersedes the Continuum Hypothesis.

## 5. Implications for Proofs by *Reductio (ad Absurdum)*

The refutations of Cantor's proofs of nondenumerability (Section 3) do not question the reliability of all mathematical proofs by contradiction, but only what constitutes a logically acceptable implementation of such a method of proof. Given that *reductio* methods account for a large share of all published mathematical proofs, it could prove unwise to assume that the errors made by Cantor in his proofs are an isolated occurrence.

The formal simplicity of the logical evaluations of Cantor's proofs presented in Section 3 suggests that the implementation of a suitable quality control for this type of proof does not need to be an onerous task. The first requirement, a method to identify incorrect mathematical statements, is addressed by the following definition.

**Definition 5.1.** *A mathematical statement Q is said to be inconceivable when there is another statement P such that*

***i)*** $(Q \Rightarrow P) \wedge (Q \Rightarrow \neg P)$, *or* ***ii)*** $Q \Rightarrow ((P \Rightarrow \neg P) \vee (\neg P \Rightarrow P))$.

*Otherwise, the statement Q is considered conceivable.*

The mathematical statement $Q$ can itself define a given mathematical object or entity. In such a case, this object is also considered inconceivable or conceivable, in line with the statement which defines it.

When examining a mathematical proof, the above definition can provide the means to recognise any statement that undermines consistency. The first part of Definition 5.1 might be considered unnecessary, since it is obvious that any statement $Q$ implying both the truth and the falsehood of another statement $P$, i.e. $Q \Rightarrow (P \wedge \neg P)$, is contradictory and, consequently, false. However, Section 3.1 points to the irreparable damage that the inclusion of an inconceivable statement causes to the validity of any proof by internal (self-referential) contradiction (i.e. a proof where the assumption $\neg P$ leads to its own negation $P$ through a direct chain of inference). The first part of Definiton 5.1 provides, therefore, a reminder of the need to avoid such a situation (see Principle 5.2 below). It is apparent that the statements found at fault in Cantor's various proofs of nondenumerability (Section 3) should be classified as inconceivable, according to Definition 5.1.

The second part of Definition 5.1 is aimed at a different type of mathematical statement, i.e. sentences of a self-referential nature. As discussed in Section 3.3.3, a final refutation of Cantor's Theorem could only be claimed once the denumerability of $\mathcal{P}(\mathbb{N})$ had been established (by Theorem 4.1). An important implication of such a refutation is that the falsehood of a self-referential statement (which invariably leads to contradiction) is not necessarily linked to the falsehood of the original assumption. Although the self-referential statement appears to be a logical consequence of the original assumption, this is not the case. This means that the initial statement could be true (as in the proof of Cantor's Theorem), despite the contradiction generated by the self-referential statement (which can never be true). Thus, the logical connection between assumption and self-referential statement (which, unfortunately, can seem most reasonable) is not a sound one.

As proposed here, a suitable way of dealing with self-referential statements is to classify them as inconceivable, since they undermine consistency. This is achieved by the second part of Definition 5.1. Again, Theorem 4.1 is the pivotal justification of the definition.

The main objective of formulating the two parts of Definition 5.1 is, therefore, to avoid mistakes. The refutations presented in Section 3 provide sufficient justification to search for, and exclude, both types of inconceivable statements from mathematical proofs, with one exception presented below.

**Principle 5.2 (of Conceivable Proof).** *No mathematical proof can be judged valid if its construction includes an inconceivable statement; the exception is if the purpose of the proof is to demonstrate the falsehood of an inconceivable statement, provided that the resulting contradiction is not conceptually linked to the initial assumption of the proof.*

As already indicated, the fact that (by definition) an inconceivable statement leads to contradiction implies that such a statement is always false. Such an observation justifies the exception included in Principle 5.2, required to account for the possibility of using an inconceivable statement as the starting point of a *reductio* proof. Situations could be envisaged whereby, in order to prove the truth of a given statement $P$, the usual assumption is made (i.e. that the negation $\neg P$ is true), only to discover from the workings of the proof that $\neg P$ is inconceivable and, therefore, false. It should be no surprise that proofs of this nature can be found in the mathematical literature [16,45] [10], an observation that reinforces the logical strength of Principle 5.2.

The classification of mathematical statements as conceivable/inconceivable facilitated by Definition 5.1 is not necessarily equivalent to the more elementary one of true/false [11] although, as stated above, every inconceivable statement is false. The opposite notion, i.e. "every conceivable statement is true", would imply that "every false statement is inconceivable" (i.e. always leads to contradiction). While this might well be the case, the classification of self-referential statements as inconceivable (and the consequent blockade to their incorporation into mathematical proofs) is a key justification of Definition 5.1, regardless of any potential duplicity of classifications.

In any mathematical system, the truth or falsehood of a statement is best asserted by proof. As Cantor's arguments illustrate very well, it may be difficult to be absolutely certain that a given statement is conceivable and, consequently, suitable for inclusion in valid proofs. A solution could be to propose that, if the examination of the logical form of the constructed proof does not indicate any of the inconsistencies covered by Definition 5.1, then no inconceivable statements play a part. It should also be noted that such a test cannot rule out other (unrelated) causes of invalidity.

[10] A particularly illustrative example was provided by James Clarkson in his proof (published in 1966) of the divergence of the series of prime reciprocals [16]. By making the initial assumption that the series of prime reciprocals is convergent, Clarkson was able to conclude that another suitably constructed series is both convergent and divergent [16,45]. Clearly, this contradiction could only be reconciled with the falsehood of the starting assumption (i.e. the convergence of the series of prime reciprocals), thus completing the proof. Furthermore, and according to Definition 5.1, the false assumption has to be classified as an inconceivable statement.

[11] To affirm that any given statement is either true or false implies the tacit acceptance of the Law of Excluded Middle [22], also known as Principle of Bivalence [52].

While the inconceivable nature of statements based on the first part of Definition 5.1 appears to be caused by recognisable logical errors, a similar conclusion for self-referential statements is more problematic. As discussed in Section 3, the conventional way of dealing with these statements has been to reformulate the axiomatic and hierarchical structure of set theory, in an attempt to eliminate any associated paradoxes/antinomies, with a predictable and inevitable increase in complexity [35,56]. The introduction of Definition 5.1 and Principle 5.2 offers an alternative strategy and preserves an appealing level of simplicity. Furthermore, as embodied by the following conjecture, any attempt to eradicate from mathematical theory all contradictory elements could well be condemmed to eventual failure.

**Conjecture 5.3 (of Logical Imperfection).** *Any sound and/or consistent system of mathematics is capable of generating inconceivable statements.*

The purpose of the above conjecture is to acknowledge that any mathematical language, whatever degree of sophistication and coherence it might possess, can be used to construct inconceivable statements. The faithful adherence to the semantic rules of such a language does not provide a guarantee that logical impossibilities are totally avoided. Self-referential statements illustrate this point very well. By analogy, there are countless examples of spoken language where, despite a complete respect for the grammatical rules, sentences are constructed which are clearly nonsensical (e.g. the red ball is blue) [12].

The limitations that Conjecture 5.3 impose on the formal scope of mathematical language are not an impediment to the generation of mathematical truths. Given that these depend entirely on the construction of valid proofs, the Principle 5.2 (of Conceivable Proof) provides a precise framework in which to operate safely. Although this is particularly relevant to the use of *reductio* proofs, it is important to point out that Principle 5.2 does not preclude the construction of satisfactory arguments by contradiction.

Perhaps an appropriate way to illustrate this last point is to describe a classical example of a *reductio* proof that complies comfortably with Principle 5.2: the Euclidean proof of the infinitude of primes (text taken from [2]).

> **Theorem** (Euclid's Elements, Book IX, Proposition 20.) *There are infinitely many prime numbers.*

---

[12] It is of interest that the view expressed here on the limitations of mathematical languages are not that dissimilar from the intuitionist stand on this matter [28].

***Proof.*** If we suppose by way of contradiction that there are finitely many prime numbers, then we may assume that $S = \{p_1, \cdots, p_n\}$ is the set of *all* prime numbers. We let $m = p_1 \cdots p_n$, the product of all the primes, and we let $q = m + 1$. Since $q > p_i$ for all $i$, we see that $q$ is not in $S$, and therefore $q$ is not prime. Then there exists a prime $p$ that is a divisor of $q$. Since $p$ is prime, then $p = p_j$ for some $j$, so that $p$ is a divisor of $m$. But if $p$ divides both $m$ and $q = m + 1$, then $p$ divides the difference $q - m = 1$. However, this is impossible, so we have obtained a contradiction. Q.E.D.

As in any conventional proof by contradiction, the above proof starts by assuming that the negation of the statement to be proved is the true one. The proof can be then broken down into its individual statements:

- $P$ = 'There are infinitely many prime numbers'
- $\neg P$ = 'There are finitely many prime numbers'
- $Q_1$ = '$S = \{p_1, \cdots, p_n\}$ is the set of *all* prime numbers'
- $Q_2$ = '$m = p_1 \cdots p_n$ is the product of all the primes, and $q = m + 1$'
- $Q_3$ = '$q > p_i$ for all $i$, so $q \notin S$'
- $Q_4$ = '$q$ is not prime'
- $Q_5$ = 'There exists a prime $p$ that is a divisor of $q$'
- $Q_6$ = '$p = p_j$ for some $j$, so that $p$ is a divisor of $m$'
- $C$ = '$p$ divides the difference $q - m = 1$, which is impossible'

Accordingly, the above list of statements constitute the logical sequence

$$\neg P \Leftrightarrow Q_1 \Leftrightarrow Q_2 \Leftrightarrow Q_3 \Rightarrow Q_4 \Leftrightarrow Q_5 \Leftrightarrow Q_6 \Rightarrow C\,. \tag{5.1}$$

The connectives between $\neg P$, $Q_1$, $Q_2$ and $Q_3$ are all biconditional (iff), since the four statements say the same thing (the number $q$ carries implicitly its own definition). The connective between $Q_3$ and $Q_4$ is not biconditional because the property of not being prime does not necessarily make a number greater than all the primes $p_i$ in $S$. However, it can be argued that the connectives between $Q_4$, $Q_5$ and $Q_6$ are all iffs, because of the definitions of $q$ and $m$. The connective between $Q_6$ and $C$ is not biconditional; it is possible for a number to be a divisor of a difference between two other numbers, without dividing either of them. Finally, the falsehood of the last statement $C$ can be independently assessed (proved) and does not depend on $q$ being prime or not.

What this analysis illustrates is that the proof sequence (5.1) does not contain any statement which could be classified as inconceivable. This is made easier by the fact that the final contradiction ($C$) has no connection with the original assumption ($\neg P$), i.e. the argument is an instance of proof by external contradiction (Section 3.1). In summary, this proof complies entirely with Principle 5.2.

## 6. Set Theoretical Implications of Denumerability

Although Cantor's initial attempts to construct a theory of transfinite numbers did not instantly receive the general acclaim that they later enjoyed [5,23], the natural beauty and simplicity of his original theory of sets gave a younger generation of mathematicians the tools needed to build the formal foundations of arithmetics and, by extension, the rest of mathematics [23,35]. However, the initial simplicity did not last for long, as numerous and varied paradoxes soon troubled set theory [35,49,56].

A systematic and detailed examination of all these paradoxes is beyond the scope of this paper. Nevetheless, two main reasons for the occurrence of such contradictions can be mentioned: firstly, the construction of self-referential objects (e.g. Russell's Antinomy [22, 35, 56]); secondly, the construction of transfinite cardinal/ordinal entities (e.g. Cantor's Paradox and Burali-Forti Paradox [35, 49, 56]). Other paradoxes are both set theoretical and logical in nature (e.g. Skolem's Paradox [37,56]).

From Russell onwards, logicians and mathematicians have faced up to the challenge presented by the set theoretical paradoxes with relative success [35,56]. But this success has not been achieved without significantly increasing the complexity of the axiomatic formulation of set theory, through a hierarchy of levels and histories [35,56].

The suggestion made here is that the concepts introduced in the preceding section, encapsulated in Principle 5.2 (of Conceivable Proof), go a considerable way to alleviate the problems generated by the self-referential paradoxes, while maintaining the simplicity of a more naïve formulation of set theory [2]. Also, Principle 5.2 stems from the most fundamental need to eliminate erroneous proofs, something that more elaborate set theories have failed to do.

The second source of paradoxes in set theory is the postulation of transfinite numbers. However, the refutation of Cantor's proofs of nondenumerability (Section 3) and construction of proofs of denumerability (Section 4) leave little scope for the existence of transfinite numbers, prompting the proposal of the following conjecture.

**Conjecture 6.1 (of Countable Infinity).** *All sets are countable, so all infinite sets are denumerable.*

The elimination of transfinite numbers implies the disappearance of most (if not all) of the paradoxes not caused by self-referential statements, and reinforces the subsequent call for a simplification of the axiomatic principles of set theory. Focusing on the standard ZF axioms [35, 56], it could be suggested that any axioms only required to facilitate the formulation of transfinite number theory might become redundant or peripheral.

6.1. **The Axiom of Choice (AC)**. AC is the source of some puzzling results, not only in set theory, but also in other areas of mathematics, such as analysis, geometry and game theory [39]. However, mathematics without AC is equally complicated by paradoxical outcomes [39]. Consequently, the use of AC is now regarded as an acceptable necessity, though to be avoided if at all possible [56].

The refutations of Cantor's proofs change the circumstances regarding the use of AC. It can be argued that the problems traditionally associated with AC have little to do with AC itself, but are the consequence of nondenumerable infinities, particularly $\mathbb{R}$ [39]. Therefore, it is suggested that the above problems will disappear, once nondenumerability is abandoned. If this is the case, AC will have a simpler role to play, free of controversy or paradox.

6.2. **Actual (completed) infinity**. Until the second half of the 19$^{th}$ Century, mathematicians remained reluctant to accept the inclusion of actual completed infinities into the workings of mathematics [5]. That attitute eventually changed when the work of Cantor and Dedekind with infinite sets gained general acceptance. However, the refutations of Cantor's arguments for nondenumerability (Section 3) are in part the result of a rigorous interpretation of the Axiom of Infinity [35,56] with regard to the construction and reordering of countable infinities (Theorems 3.6 and 3.7).

Theorem 3.6 (of Actual Countable Infinity) demonstrates that the concept of denumerability and, with it, the process of construction of (countably) infinite sets stem directly from the well-ordering property of $\mathbb{N}$ and the resulting principles of mathematical induction [2]. The significance of Theorem 3.6 resides in accepting that all the elements of an infinite set exist, despite the fact that the construction of the set can never be manually completed. Acceptance of the existence of actual (completed) infinities, as encapsulated in Theorem 3.6 appears to be a necessary prerequisite for the manipulation of infinite sets.

6.3. **Cardinality and infinite sets**. It is important to remember that the apparent benefit of Cantor's transfinite theory was to provide the means to compare the size of infinite sets [23,35]. The concept of nondenumerability was sufficient to establish, for example, that $|\mathbb{R}| > |\mathbb{Q}|$ [2]. The results presented in Section 4 imply that such means for comparison are no longer available, particularly if Conjecture 6.1 (of Countable Infinity) is accepted.

Dedekind defined a set as infinite when it is possible to put all its elements into one-to-one correspondence with the elements of at least one of its proper subsets (Section 2). However, this definition does not provide a method able to determine the actual cardinality of the set. Cantor made use of Dedekind's definition to compare infinite sets, by the implicit assumption that some sets could be nondenumerable. But, if all infinite sets are denumerable, such a way of comparing them can only confirm that they are indeed countable. In summary, using Cantor's criterion for comparing the cardinality of infinite sets, transfinite cardinals (alephs) reduce to just one ($\aleph_0$), so all infinite sets have the same cardinality.

In the above scenario, $\aleph_0$ is the universal cardinality of all infinite sets. That makes it an absolute value that conveys a defining property of the sets, hence a qualitative property deprived of numerical meaning. If the requirement is to compare infinite sets in a more meaningful way, a different criterion is needed.

The concept of diagonal cover (*Dc*), defined in Section 3.2.2 to analyse Cantor's diagonalization argument, suggests the viability of an alternative rationale to compare the cardinality of infinite sets. Although the cardinality of denumerable sets ($\aleph_0$) cannot be used as a comparator, their relative constructions may fulfil this role. Consider two infinite sets. Provided there are formulae $\phi(n)$, available for both sets and relative to each other, that translate their cardinalities from the finite form into the final infinite value, it will always be possible to calculate the limiting value of the ratio, as $n$ grows without bound ($n \to \infty$). Although the limiting value of both formulae $\phi(n)$ is always infinite ($\aleph_0$), the calculated ratio can take various values, so a comparison arises.

Two definitions are needed.

**Definition 6.2 (Relative Cardinality of Finite Sets).** *Consider two finite sets, A and B, with cardinalities* $|A| = a$ and $|B| = b$. *Their relative cardinality is defined as the ratio* $\rho = |A|/|B| = a/b$.

**Definition 6.3 (Relative Cardinality of Infinite Sets).** *Consider two sets, A and B, both denumerable. Assume their constructions generate formulae,* $\phi_A(n)$ *and* $\phi_B(n)$, $\forall n \in \mathbb{N}$, *which render the cardinalities of the respective interim finite sets, in relation to each other. The relative cardinality of A and B is defined as the limiting ratio*

$$\rho = \lim_{n \to \infty} \frac{\phi_A(n)}{\phi_B(n)} \quad . \tag{6.1}$$

The implementation of Definition 6.3 will not always be easy. The critical difficulty lies, not only in obtaining the formulae $\phi_A(n)$ and $\phi_B(n)$, but also in being certain that the finite cardinalities they deliver can be genuinely compared. Since all denumerable sets can be put on a one-to-one correspondence with $\mathbb{N}$, the constructions used to establish the bijections $\mathbb{N} \to A$ and $\mathbb{N} \to B$ could also be used to extract the formulae $\phi_A(n)$ and $\phi_B(n)$, without a guarantee that such formulae are truly comparable. This essential point can be illustrated with a very simple example.

**Example 6.4.** Consider two subsets of $\mathbb{N}$, the set of all the natural numbers divisible by 2, $A = \{2, 4, 6, 8, 10, \cdots\}$, and the set of all the natural numbers divisible by 3, $B = \{3, 6, 9, 12, 15, \cdots\}$. According to the conventional way of determining the cardinality of infinite sets, the bijections $\mathbb{N} \to A$ and $\mathbb{N} \to B$ are constructed by $n \to n \times 2$ and $n \to n \times 3$ ($n \in \mathbb{N}$, $n \neq 0$), respectively, meaning that both sets $A$ and $B$ have the same cardinality as $\mathbb{N}$ ($\aleph_0$). However, to determine the relative cardinality $\rho$ between $A$ and $B$, a different rationale needs to be followed. Consider a finite set $N$, consisting of only the natural numbers up to a certain value $n$, $N = \{0, 1, 2, 3, \cdots, n\}$. For $N$, the total number of elements divisible by 2 is given by $\phi_A(n) = \lfloor n/2 \rfloor$, the greatest integer in the fraction $n/2$. In similar fashion, the total number of elements divisible by 3 results from $\phi_B(n) = \lfloor n/3 \rfloor$. Applying (6.1) results in $\rho = 3/2$, which states that $\mathbb{N}$ has, in relative terms, 1.5-fold more numbers divisible by 2 than divisible by 3.

This example shows that Definition 6.3 provides a method to calculate the cardinality of infinite sets which is a credible alternative to Cantor's criterion (of establishing or not a bijection between all the elements of the sets). The new method is precisely determined by Definition 6.3, implying that no conflict of interpretation is possible, so no contradiction arises. Since equation (6.1) determines a limiting value, it transforms the output of the calculation from one that applies to two finite sets, to one associated with the infinite sets $A$ and $B$. Therefore, $\rho$ is a genuine property of the two sets, taken in conjunction as the ordered pair $(A, B)$.

The above analysis makes $\rho$ a key comparator between denumerable sets. This observation can be encapsulated by a new definition, derived from Definitions 6.2 and 6.3.

**Definition 6.5.** *Two sets, A and B, are said to be equicardinal when their relative cardinality is* $\rho = 1$.

This definition, which applies to either finite or infinite sets, can be particularly useful in the evaluation of infinite sets. In this context, the refutation of Cantor's diagonalization argument acquires a special significance by highlighting the importance of $\rho$ in certain constructions. For example, as noted in Section 3.2.2 for the set of all infinite binary strings $S_T$, the diagonal cover $Dc$ (i.e. the relative cardinality $\rho$ between the set of diagonal positions and $S_T$ itself) takes the limiting value of the equation $\phi_A(n)/\phi_B(n) = n/2^n$, as $n \to \infty$. Since $Dc$ (i.e. $\rho$) reaches an infinitely small value (not only is $\rho < 1$, as a limit $\rho$'s value is 0), it is impossible for the enumeration along the diagonal (represented in its finite form by $\phi_A(n) = n$) to ever match the listing of all the binary strings in $S_T$ (represented by $\phi_B(n) = 2^n$), despite $S_T$ being a countable set, as proved by Theorem 4.1 (Denumerability of the Power Set). This impossibility is a direct consequence of the fact that the sets under consideration, $\mathbb{N}$ and $S_T$, are not equicardinal.

Infinite sets with relative cardinality $\rho = 1$ (i.e. equicardinal) are not common, though not impossible (e.g. the two sets of all even/odd natural numbers). The conclusion is that a range of relative cardinalities, naturally expected of finite sets, can also be extended to infinite sets. Furthermore, they can be calculated.

One case that stands out in particular is the power set, $\mathcal{P}(\mathbb{N})$. The proofs of Theorem 4.1 imply that the formulae to use in equation (6.1), $\phi_A(n)/\phi_B(n) = n/2^n$, can be appropriately compared. Hence, the value of $\rho$ is infinitely small when comparing $\mathbb{N}$ against $\mathcal{P}(\mathbb{N})$. The following example, which replicates the third proof of Theorem 4.1 illustrates how the denumerability of $\mathcal{P}(\mathbb{N})$ is compatible with this observation.

**Example 6.6.** Consider $\mathbb{N}$, the set of the natural numbers. Consider also, within the enumeration of the elements of $\mathbb{N}$, all the powers of 2:

$$\mathbb{N} = \{0, \mathbf{1}, \mathbf{2}, 3, \mathbf{4}, 5, 6, 7, \mathbf{8}, 9, 10, 11, 12, 13, 14, 15, \mathbf{16}, 17, 18, 19, 20, \cdots\} .$$

The powers of 2 can be selected to form the corresponding subset $N_2 \subseteq \mathbb{N}$:

$$N_2 = \{2^0, 2^1, 2^2, 2^3, 2^4, \cdots\} = \{2^k, \forall k \in \mathbb{N}\} . \tag{6.2}$$

It is trivial that the subset $N_2$ is an infinite set. Hence, the enumeration of the powers of 2 reproduces the enumeration of all elements $k$ of $\mathbb{N}$ (the usual bijection). Calculate now the relative cardinality of $N_2$ and $\mathbb{N}$, by applying equation (6.1); $A$ becomes $N_2$ and $B$ becomes $\mathbb{N}$. The result is

(6.3) $$\phi_{N_2}(k)/\phi_{\mathbb{N}}(k) = (k+1)/(2^k+1) \ .$$

The extra unit of both terms of the fraction needs to be included to account for the number 0 in both sets. The limiting value of equation (6.3), as $k$ grows without bound, becomes the corresponding value of $\rho$, which is infinitely small (effectively, 0).

Now consider $\mathbb{N}$ and its power set, $\mathcal{P}(\mathbb{N})$. To determine, as before, the relative cardinality of these two sets, equation (6.1) is applied, taking $A$ to be $\mathbb{N}$, and $B$ to be $\mathcal{P}(\mathbb{N})$. This produces

(6.4) $$\phi_{\mathbb{N}}(n)/\phi_{\mathcal{P}(\mathbb{N})}(n) = n/2^n,$$

where $n = 2^k+1$. The relative cardinality $\rho$ takes the same value as derived from equation (6.3).

This example illustrates that the relationship between $N_2$ and $\mathbb{N}$ is, in practice, the same as between $\mathbb{N}$ and $\mathcal{P}(\mathbb{N})$, as used by the third proof of Theorem 4.1 (Section 4.1.3). $\mathbb{N}$ is not the power set of $N_2$, but the above proof shows that the power set of $N_2$, $\mathcal{P}(N_2)$, has a one-to-one correspondence with $\mathbb{N}$. This bijection, together with the acknowledged bijection between $N_2$ and $\mathbb{N}$, implies that the power set of $\mathbb{N}$, $\mathcal{P}(\mathbb{N})$, has a one-to-one correspondence with the power set of $N_2$, $\mathcal{P}(N_2)$, i.e. with $\mathbb{N}$ itself. Consequently, $\mathcal{P}(\mathbb{N})$ is a denumerable set [13].

Example 6.6 highlights the usefulness of the relative cardinality defined by equation (6.1). When applied to a denumerable set and its power set, it is apparent that the sequence of larger and larger cardinalities that Cantor aimed to substantiate with his theory of transfinite numbers can be replicated by the relative cardinalities introduced in Definition 6.3. This is achieved without violating the countable nature of all infinite sets.

## 7. Logical Implications of Denumerability

Of the various mathematical disciplines affected by the concept of nondenumerability, logic and metamathematics are particularly significant examples [12,20,37]. The well-established concepts of incompleteness and undecidability are a direct consequence of the assumed reality of uncountable sets [12,59]. In turn, arithmetic incompleteness has had its own impact, not only on mathematics, but also in areas further afield such as computational science [12] and artificial intelligence [53,54]. This implies that the refutations of the proofs of nondenumerability presented in Section 3 are likely to have a very profound effect on the above areas.

---

[13] The construction of the subset $N_2$, as described in Example 6.7, can be repeated for those values of $k$ that are themselves powers of 2. Thus, the new infinite subset of $N_2$ constructed in such way will be equinumerous with $\mathbb{N}$. $\mathbb{N}$, in turn, will also have a one-to-one correspondence with the power set of the power set of the new set. The process can be repeated again and again. Consequently, the rationale of Example 6.6 reinforces the result already encapsulated by Corollary 4.2, i.e. the power set of any denumerable set is itself denumerable.

Using nondenumerability as a working concept [20,37], logicians were quick to adopt the diagonal arguments used by Cantor (in various formats) for the construction of their proofs of incompleteness and/or undecidability [12,31,59,60]. Consequently, a re-examination of the above proofs is required, which should be based on the refutations described in Section 3, combined with the analysis presented in Section 5 on the use of conceivable statements in the construction of valid proofs by contradiction (Principle 5.2 of Conceivable Proof). In the case of all diagonal arguments, the key requirement for validity is that the infinite sets involved in the construction of such arguments have to be shown to be equicardinal (according to Definitions 6.3 and 6.5). But this cannot be achieved, since the relation between the two sets is invariably one of linear *vs* exponential magnitudes [20, 59]; this observation explains why all the examples examined in Section 3 fail to comply with Principle 5.2.

The diagonal arguments used in mathematical logic can be classified into two broad but distinctive categories [20,59].

- Firstly, the proofs which closely follow Cantor's method and implement the original rationale of the diagonalization argument. Using a certain enumeration as a starting platform, they attempt to construct a new entity (set, function, etc) which then is claimed not to be part of the enumeration. Informally, it could be said that such an entity has been 'diagonalized out'.
- Secondly, a very important family of proofs (Gödel's proofs of incompleteness [31,38] among them) turn the previous argument on its head and seek to construct a vital mathematical object (e.g. a primitive recursive function), whose existence is claimed and leads to the final conclusion. Informally, it could be said that such an object has been 'diagonalized in'. This argument was formalised as the Diagonalization Lemma [59].

As described in Sections 7.1 and 7.2, a broad analysis of the two types of diagonal arguments is sufficient to refute them. Inevitably, such refutations also imply the refutation of the proofs they support. The historic significance of Gödel's Incompleteness Theorems [31, 59] merits further analysis which, based on Principle 5.2, is described in Section 7.3. The introduction by Gödel's theorems of logical statements which are 'true-but-unprovable in the system' [59] has been followed over the years by the description of supporting examples, which are discussed in Section 7.4. Gödel's work had a negative impact on Hilbert's Programme [41] that prompted others to expand on Gödel's results. Proofs of mathematical undecidability soon followed which also relied on diagonal arguments, as outlined in Section 7.5.

The refutations of the proofs of incompleteness and undecidability based on diagonal arguments are not sufficient to eradicate these concepts completely, since alternative proofs have been proposed that avoid diagonalization [4, 59]. However, self-referential statements are in conflict with Principle 5.2 (of Conceivable Proof). Only a proof of the completeness of first-order arithmetic can totally disprove the opposing concept of incompleteness and undecidability. Section 7.6 describes the construction of such a proof, based on the denumerability of $\mathcal{P}(\mathbb{N})$ (Theorem 4.1), which fully reinstates Hilbert's Programme [41].

7.1. **'Diagonalizing out' (the Diagonalization Argument)**. Section 3 provided a comprehensive refutation of this type of proof. However, there are additional considerations to be included. The system of logic deals with the expression (within a formal language $L$) of numerical properties by open well-formed formulae (wffs), $\varphi(x)$, with a direct influence on the enumeration of such wffs. Adhering to the standard terminology, a general definition states that [59]

> A property $P$ is *expressed* by the open wff $\varphi(x)$ with one free variable in an arithmetical language $L$ iff, for every $n$,
>
> if $n$ has the property $P$, then $\varphi(\mathrm{n})$ is true,
>
> if $n$ does not have the property $P$, then $\neg\varphi(\mathrm{n})$ is true.

n represents the standard numeral expressed in the arithmetical language in use. This definition requires the evaluation of $\varphi(n)$ [14] for all possible values of $n$ ($\in \mathbb{N}$), which have (or not) the property $P$. Accordingly, the total expression of $P$ by $\varphi(n)$ can be represented by an unending binary string of 0s and 1s (as described in Section 3.3.1), where 0 and 1 indicate the opposing outcomes; the choice of symbols is obviously arbitrary. If this evaluation were implemented for a certain set of wffs, the result would be a collection of binary strings, each corresponding to a specific mapping $\mathbb{N} \to \{0,1\}$. Although it is possible for two or more different wffs to be associated with the same binary string, this scheme (allowing for repetitions) produces a workable enumeration of the whole set of wffs, provided such a set is countable [59].

A numerical property $P$ can also be associated with its characteristic function, which is then defined as the one-place function, $f(n)$, such that $f(n)=0$ when the numeral $n$ has the property $P$, and $f(n)=1$ when it has not [59]. This definition allows for a scheme of representation similar to the one just described for wffs, where the characteristic functions are total (i.e. the mapping is defined for all possible values of $n \in \mathbb{N}$). Again, the resulting collection of binary strings acts as the enumeration of the set of characteristic functions.

Although greatly simplified [15], this analysis provides a sufficient platform to identify the flaw in the diagonal arguments used by logicians [20,59]. When the purpose of the argument is to show that it is not possible to obtain a countable enumeration of the given set of binary strings, the construction of the diagonal object replicates Cantor's original rationale in 'diagonalizing out' such object, thus repeating Cantor's error (the '$n$ *vs* $2^n$' disparity).

A rather different scenario is encountered when it is known (by independent means) that the set represented by the collection of binary strings is effectively enumerable [59]. Since this implies that the set under consideration is countable, the purpose of the diagonal argument is no longer to prove nondenumerability, but to construct a diagonal object that, although it is computable, cannot be part of the effective enumeration. Some very significant theorems use this kind of proof on their way to incompleteness, for example the claim that "*there are denumerable sets that are not effectively enumerable*" [59], so its refutation has far reaching consequences.

[14] The definition uses $x$ as the free variable which takes all the possible values of $n$, expressed as the corresponding standard numeral in the formal language used [59]. This distinction is not necessary for this discussion on enumerations, so it can be simplified.

[15] Open wffs can also express two-place (or many-place) relations [59], but the incorporation of more than one argument (free variable) is not required for the analysis that follows.

The flaw in this set of proofs can be illustrated with an important argument involving primitive recursive functions [59]. It is widely accepted that there are effectively computable numerical functions that are not primitive recursive, with verifiable examples readily available, such as the Ackermann-Péter functions [59]. However, this observation has also been claimed with a specifically engineered example which deploys a diagonal argument [59]. The starting point in the argument is to prove that the set of all primitive recursive (p.r.) functions is effectively enumerable, and to use the same construction to build a corresponding table in which these p.r. functions, $f_m(n)$, are enumerated (repetitions allowed) against the free variable, $n$. As discussed previously, such effective enumeration equates each p.r. function to a binary string, with $n$ covering its whole domain ($n \in \mathbb{N}$):

(7.1)

| | 0 | 1 | 2 | 3 | ... | $n$ | ... |
|---|---|---|---|---|---|---|---|
| $f_0$ | $f_0(0)$ | $f_0(1)$ | $f_0(2)$ | $f_0(3)$ | ... | $f_0(n)$ | ... |
| $f_1$ | $f_1(0)$ | $f_1(1)$ | $f_1(2)$ | $f_1(3)$ | ... | $f_1(n)$ | ... |
| $f_2$ | $f_2(0)$ | $f_2(1)$ | $f_2(2)$ | $f_2(3)$ | ... | $f_2(n)$ | ... |
| ... | ... | ... | ... | ... | ... | ... | ... |
| $f_m$ | $f_m(0)$ | $f_m(1)$ | $f_m(2)$ | $f_m(3)$ | ... | $f_m(n)$ | ... |
| ... | ... | ... | ... | ... | ... | ... | ... |

Accordingly, for each specific value (0 or 1) of $f_m(n)$, $n$ signifies the position of the digit in the binary string, and $m$ signals the corresponding p.r. function(s) in the enumeration.

The proof continues with the construction of a diagonal function defined as $\delta(n) = f_n(n) + 1 \pmod 2$. Since $\delta(n)$ differs from each function $f_m$ in at least one digit of the binary string, the proof concludes that $\delta(n)$ is not part of the enumeration, thus cannot be a p.r. function. But, on the other hand, since the computation of $\delta(n)$ follows the same construction used to evaluate the whole set of p.r. functions $f_m(n)$, i.e. a step-by-step mechanical procedure, such function $\delta(n)$ is effectively computable [59].

In line with the analysis of the preceding sections, it can be concluded that the above proof fails when defining the diagonal function $\delta(n)$. The definition assumes that it is possible to compute $f_n(n)$ for all values of $m$ $(= n)$. However, table (7.1) uses two infinite enumerations, the sequence of digits ($n$) in the binary strings, and the collection of p.r. functions ($f_m$), which by virtue of the construction used cannot constitute equicardinal sets (Definition 6.5) [16].

7.2. '**Diagonalizing in' (the Diagonalization Lemma)**. The refutation described above provides a convenient link between proofs that 'diagonalize out' and those that 'diagonalize in' a given mathematical entity. 'Diagonalizing in' represents a complete reversal of Cantor's original rationale, since the starting point is the effective (thus countable) enumeration of a certain collection of functions, wffs or sets used to claim the existence of a suitably defined object which is then extracted from the enumeration. This is the object that takes the rationale of the proof to its final conclusion. It is important to emphasize this point because, without the sound existence of such an object, the proof loses all validity.

The collections under consideration are those of wffs (with one free variable), expressed in the language of a formal theory $T$, as described previously (Section 7.1). If a suitable method were available to effectively enumerate these wffs, $\varphi_m(n)$, a table similar to (7.1) could be constructed, of $\varphi_m$ *vs* $n$, where $n$ is the free variable and the index $m$ signifies the effective enumeration of the wffs. As before, the correct description of each wff is carried out by the corresponding binary string, where each digit (0 or 1) indicates whether the specific value of $n$ satisfies (or not) the wff $\varphi_m$.

The enumeration of the wffs $\varphi_m$ requires an appropriate numbering system, ideally a coding method which translates the syntactic meaning of each wff into a numerical code uniquely assigned to it. This process (The Arithmetization of Syntax [59]) could be carried out in different ways, but the original numbering method was developed by Gödel for his famous Incompleteness Theorems [31,38]. The intricacies of Gödel's numbering scheme (described in [59]) do not require further explanation; it is sufficient to understand that the effective enumeration of the collection of wffs $\varphi_m$ is achieved using as indexes $m$ the Gödel numbers (g.n.) of the wffs in the set, represented by $\ulcorner\varphi\urcorner$.

The diagonalization at the heart of these proofs operates a substitution where the free variable $n$ takes the value of the $\varphi$'s g.n., i.e. $\varphi(n) = \varphi(\ulcorner\varphi\urcorner)$. The critical step is to assimilate this operation on expressions into a p.r. function (in the case of Gödel's proof), $diag(n)$, appropriately defined by the following theorem [59]:

> **Theorem** *There is a p.r. function diag*$(n)$ *which, when applied to a number n which is the g.n. of some wff, yields the g.n. of that wff's diagonalization.*

Writing $D(\varphi) = \varphi(\ulcorner\varphi\urcorner)$ as the diagonalization of $\varphi$, the above theorem claims the existence of a p.r. function $diag(n)$ such that $diag(\ulcorner\varphi\urcorner) = \ulcorner D\urcorner$. Once the validity of this p.r. function, $diag(n)$, is accepted, it becomes a very powerful tool for the construction of additional p.r. relations and wffs, leading to the final conclusion.

---

[16] The discrepancy '$n$ *vs* $2^n$' cannot be eliminated by possible repetitions, i.e. by the fact that the same binary string can be the representation of more than one p.r. function. As emphasized in previous sections, the relative cardinality of the two sets ($\rho$), not only is smaller than 1, it is infinitely small. Only a very small subset of p.r. functions can be equicardinal to the set of digits in the binary strings (e.g. consider the functions '$s$ is the successor of $n$', which correspond to strings with only one 1 and an infinity of 0s). Consequently, the diagonal enumeration leaves out an infinitely larger collection of p.r. functions, which $\delta(n)$ can never account for.

Gödel's proof is one particular example of the general theorem encapsulating this process, termed the Diagonalization Lemma [59]:

> **Theorem** *If $T$ is a nice theory and $\varphi(x)$ is any wff of its language with one free variable, then there is a sentence $\gamma$ of $T$'s language such that* $T \vdash \gamma \leftrightarrow \varphi(\ulcorner\gamma\urcorner)$. [17]

The sentence $\gamma$ is the diagonalization of a preceding wff, i.e. $\gamma(\psi) = \psi(\ulcorner\psi\urcorner)$, while the definition of the intermediate wff, $\psi$, makes use of $diag(n)$ [59]. This construction implies that the validity of the theorem depends directly on the validity of the diagonalization function, $diag(n)$.

Although $diag(n)$ has been conventionally proved to be primitive recursive [59], the existence of this function relies on the possibility of executing the operation of diagonalization, i.e. $D(\varphi) = \varphi(\ulcorner\varphi\urcorner)$, for the complete set of wffs $\varphi_m(n)$. But such a completeness cannot be achieved, as a refutation twinned with the analysis described in the preceding section (the disparity $n$ *vs* $2^n$) readily shows. Hence, $diag(n)$ has no conceivable existence and violates Principle 5.2.

The construction of $diag(n)$, has been termed 'diagonalizing in' to denote the slight differences between this kind of reasoning and Cantor's original rationale (termed 'diagonalizing out'). However, the fact that the definition of $diag(n)$ is effectively an inconceivable statement suggests an even closer link between both approaches. Such a link can be revealed by the implicit consequences of formulating $diag(n)$, which can be expressed as

$$diag(n) \Leftrightarrow \forall m = \ulcorner\varphi_m\urcorner \big(D(\varphi_m) = \varphi_m(m)\big)\,. \tag{7.2}$$

Thus, $diag(n)$ implies that the operation of diagonalization can be carried out for all the wffs $\varphi_m$, enumerated according to their own g.n. ($\ulcorner\varphi_m\urcorner$) [18]. This assumption would allow the construction of a new diagonal wff, $\chi(m)$, such that

$$\chi(m) = \varphi_m(m) + 1 \ (\mathrm{mod}\ 2)\ ,\ \forall m = \ulcorner\varphi_m\urcorner\,. \tag{7.3}$$

If a claim were that the new wff $\chi(m)$ differed from all the wffs $\varphi_m$ in the enumeration, the construction (7.3) could be used to 'diagonalize out' $\chi(m)$, becoming a target of the refutations of Section 7.1. Consequently, the Diagonalization Lemma could be said to 'diagonalize in' a function $diag(n)$, which is then used to construct the final sentence $\gamma$ that 'diagonalizes out' the desired conclusion. Gödel's First Incompleteness Theorem fits this description precisely, within the framework of a formal theory $T$, by constructing a sentence that is "*true because it is unprovable-in-T*" [59].

7.3. **Gödel and incompleteness**. Gödel's famous 1931 paper [31] brought to light the perceived limitations of any formal system of first-order arithmetics and,

---

[17] The term *nice* is applied to a formal theory $T$ when it is consistent, p.r. axiomatized, and extends Q (Robinson Arithmetic) [59].

[18] The universal quantifier $\forall m = \ulcorner\varphi_m\urcorner$ is not intended to introduce a self-referential character to the expression, but to highlight the fact that the indexing of the wffs $\varphi_m$ is based on their own g.n.

by extension, the whole of mathematics. Gödel's First Incompleteness Theorem can be formulated as [59]:

> **Theorem** *If T is a p.r. adequate, p.r. axiomatized theory whose language includes $L_A$, then there is an $L_A$-sentence $\varphi$ of Goldbach type such that, if T is consistent then $T \nvdash \varphi$, and if T is $\omega$-consistent then $T \nvdash \neg\varphi$.*

A detailed and comprehensive explanation of the technical aspects of the above theorem can be found elsewhere [59], and is not required here. It suffices to say that the theorem is a classical implementation of the Diagonalization Lemma, refuted above (Section 7.2). Nevertheless, the significance of Gödel's result is such that some additional considerations are worth discussing. For this purpose, the all-important Gödel sentence, $G$, can be defined in a more common, but informal way [59]:

$$G := \text{“}G \textit{ is true iff is unprovable in } T\text{”}. \tag{7.4}$$

As stated in the first half of the theorem, the definition (7.4) of $G$ holds if the theory $T$ is consistent (the weaker condition of $\omega$-consistency is not relevant to this analysis [59]). Although Gödel's formulation required that consistency, and not soundness, is a sufficient condition for incompleteness [31,59], the more general prerequisite does not imply that the theorem is seeking to deal only with unsound theories. A theory cannot be sound (i.e. it only proves true statements) without being consistent (i.e. it cannot simultaneously prove both the truth and falsehood of the same statement), so Gödel's result appears to be adequate. However, since it is an obvious and fundamental necessity to work with mathematical theories that are (or are believed to be) sound, it should also be expected that the proof of Gödel's theorem is compatible with the soundness of theory $T$.

The Principle 5.2 (of Conceivable Proof), introduced in Section 5, provides a method to evaluate any mathematical proof, particularly those dependent on the implementation of a *reductio* method. Section 3, as well as Sections 7.1 and 7.2, have shown that Principle 5.2 can be useful in the identification of inconceivable statements (Definition 5.1). Self-referential statements are a particular target of Principle 5.2. Therefore, this principle could be used to analyse the Gödel sentence $G$ (7.4). If doubts are raised on the validity of $G$, given its self-referential character, the classical answer has been that the truth of this sentence is determined by its construction, so contradictions do not arise [59]. However, the p.r. function $diag(n)$ is an inconceivable object, as Section 7.2 has previously shown, hence the Diagonalization Lemma is not a valid method of proof; this conclusion is sufficient to invalidate sentence $G$ and contradicts the classical defence of its validity. A striking observation is that, by applying the second part of Principle 5.2 to $G$ itself, this sentence can also be found to be inconceivable: $G \Rightarrow ((P \Rightarrow \neg P) \vee (\neg P \Rightarrow P))$, which corroborates the refutation of the Diagonalization Lemma. The analysis is facilitated by writing $P :=$ “*G is unprovable-in-T*”, since the sentence $G$ provides all the information required.

The first step is to assume that theory $T$ is sound and negation-complete [59]. On this basis, the sentence $P$ can be examined to identify any contradictions.

Any contradictions can then be considered as the initial conditions, i.e. soundness and completeness, are removed. Accordingly:

*i)* $G \Rightarrow (P \Rightarrow \neg P)$. Since $T$ is negation-complete, the fact that $G$ is unprovable-in-$T$ implies that $G$ is false (otherwise, if $G$ were true, it would be provable-in-$T$). Since $G$ is false, it cannot be unprovable-in-$T$ (because that is what $G$ says). Therefore, $G$ is probable-in-$T$. Contradiction.

*ii)* $G \Rightarrow (\neg P \Rightarrow P)$. $\neg P =$ "*G is provable-in-T*". Since $T$ is sound, the fact that $G$ is provable-in-$T$ implies that $G$ is true. But the sentence $G$ says that, if true, $G$ is unprovable-in-$T$. Also a contradiction.

The first contradiction depends on $T$ being negation-complete. However, this is the property of $T$ that $G$ claims to prove false. Consequently, if the condition of completeness is removed, the contradiction does not arise.

The second contradiction depends on $T$ being sound. Therefore, to eliminate the contradiction, the condition of soundness has to be removed. This is the only way to avoid the contradiction. This implies that, if $T$ is sound, then the sentence $G$ is inconceivable and fails to comply with Principle 5.2. The fact that $\neg P$ is the negation of $G$ does not rescue the proof, since the truth (or falsehood) of the assumption ($\neg P$) has no bearing on the consistency of the logical argument deployed (otherwise, no *reductio* proof could ever be constructed). It is important to reiterate that the contradiction identified by $G \Rightarrow (\neg P \Rightarrow P)$ is intrinsically linked to the soundness of theory $T$. The implication is that, if the objective is to have a conceivable proof, soundness and incompleteness (as asserted by Gödel's First Theorem) have to be mutually exclusive.

The logical evaluation of Gödel's sentence $G$ with Principle 5.2 confirms that the refutation of the Diagonalization Lemma also refutes Gödel's First Theorem. The above analysis illustrates how simple and effective the application of Principle 5.2 (of Conceivable Proof) can be as a tool for proof invalidation.

The refutation of Gödel's First Incompleteness Theorem has key implications. An immediate consequence is the refutation of Gödel's Second Incompleteness Theorem [31,59]. This theorem states that "*assuming a theory T is consistent, such a consistency is unprovable-in-T*" [59]. The second theorem is almost a corollary of the first one, so the refutation of one leads to the refutation of the other. A substantial number of proofs and results which extend these theorems are thereby compromised [59].

Gödel's use of consistency (rather than soundness) for his proofs was not an attempt to disguise the problems described, but a considered attempt to provide a comprehensive answer to the formalistic ideal embodied by the set of goals known as Hilbert's Programme [41]. Concentrating on a formal system of first-order arithmetic like PA, the first two goals of the programme were, firstly, to prove that PA is negation-complete and, secondly, to show that PA is capable of proving its own consistency. Gödel's theorems were designed to answer these two aims in the negative, in a way that did not depend on proving that the theory under examination is sound. Hence, the limitations supposedly exposed by Gödel would have been intrinsic to the fabric of arithmetics and, by extension, to the whole of mathematics. The refutations of Gödel's theorems reported in this article challenge current views. Therefore, they should bring the first two objectives of Hilbert's Programme back under consideration, at least regarding Gödel's contribution.

7.4. **True-but-unprovable**. Initially, the incompleteness that Gödel's theorems brought to mathematical theory relied upon true-but-unprovable statements of a primarily logical nature, constructed using coding devices [59]. Therefore, it was important to discover other more natural sentences which, as well as being of more intrinsic mathematical interest, could be shown to be both true and unprovable in first-order arithmetic (commonly PA). The verification of the truth and unprovability of such examples has to depend on mathematical tools that, while unavailable to basic arithmetics, are also non-Gödelian in character. The published examples achieve this goal by using transfinite numbers [59]. Yet the refutations and proofs of Sections 3 and 4 have left the concept of transfinite numbers (both ordinal and cardinal) bereft of mathematical proof. Such a lack of evidence raises fundamental objections to the inclusion (direct or implicit) of transfinite numbers in the construction of any subsequent proofs.

As already mentioned in Section 1, a case of particular interest is Goodstein's Theorem, which is purely number-theoretic in character [36,56,59]. If a first-order arithmetic proof of Goodstein's Theorem could be constructed, contradicting the alleged unprovability of the theorem [46], it would strongly challenge the veracity and validity of transfinite number theory or Gödel's incompleteness. Furthermore, since Gödel's proofs of incompleteness depend effectively on the existence of nondenumerable sets, and these are essential to the construction of transfinite numbers, a proof-in-PA of Goodstein's Theorem would ultimately undermine transfinite number theory as well as incompleteness, corroborating the results reported in this article. Such a proof has already been constructed and was reported in 2009 [55].

7.5. **Church, Turing and undecidability**. As discussed in Section 7.3, Gödel's theorems undermined the first two goals of Hilbert's Programme. A third goal, the so-called *Entscheidungsproblem* [59], still remained unresolved. As reiterated by Hilbert and Ackermann in 1928 [41], this was the problem of determining the existence of an effective method to decide whether an arbitrary sentence of first-order logic is a theorem or not [41,59]. Two independent and parallel avenues answered this question in the negative, completing the apparent destruction of Hilbert's dreams.

One strategy, developed by Church [13,14,59], provided a natural expansion to the concept of recursiveness. This was to introduce the notion of $\mu$-recursive functions which, unlike primitive recursive functions, can be constructed from initial functions by a chain of definitions, not only by composition and primitive recursion, but also by regular minimization [59]. Therefore, p.r. functions are a subset of all $\mu$-recursive functions. Since Church's Thesis presumes that the total numerical functions which are effectively computable by some algorithmic routine are just the $\mu$-recursive functions [59], $\mu$-recursiveness made it possible to build on the rationale developed by Gödel using p.r. functions and derive a key result, namely that "*the property of being a theorem of first-order logic is recursively undecidable*" (Church's Theorem) [13,14,59].

However, like Gödel's incompleteness theorems, the proof of Church's Theorem uses the Diagonalization Lemma [59]. Thus, based on the refutation outlined in Section 7.2, Church's proof has to be considered flawed and, with it, any subsequent results derived from Church's Theorem [59].

The second independent strategy used to unravel the *Entscheidungsproblem* was championed by Turing in his paper on computable numbers [38,60] which, like Church's results, was published in 1936. Turing conceived of mechanical devices (Turing machines) capable of performing arithmetical tasks. Therefore, Turing machines compute numerical functions. The extent of this property is encapsulated in Turing's Thesis, by which the numerical functions that are effectively computable by some algorithm are computable by a suitable Turing machine [59,60]. The similarities between Church's Thesis and Turing's Thesis are apparent. They are further strengthened by the proven facts that every $\mu$-recursive function is Turing-computable and, in reverse, that all Turing-computable functions are $\mu$-recursive [59]. Their equivalence led directly to the formulation of the Church-Turing Thesis which claims that "*the effectively computable total numerical functions are the $\mu$-recursive/Turing computable functions*" [15,59]. To date, the Thesis has withstood the test of time and is widely accepted.

Turing's strategy was to formulate a universal Turing machine which could be programmed to replicate every conceivable Turing machine and used to generate an effective enumeration of all Turing machines (identified by their programs). Turing's first result was to show that there is no effective procedure to decide whether a given Turing program halts when set to work with its own index-number as input, i.e. the so-called self-halting problem [59,60]. Turing's initial objective had been to apply Cantor's diagonalization argument to the enumeration of Turing machines, in order to prove by means of a diagonal construction that no Turing program exists that is capable of deciding the self-halting problem [48,59]. The refutations of Section 7.1 imply that reliance on Cantor's erroneous method of proof invalidates Turing's results.

In Turing's strategy, the claim that the self-halting problem is not solvable provides the platform to prove that other halting problems are also unsolvable, including the *Entscheidungsproblem* [59,60]. The refutation of Turing's proof for the self-halting problem inevitably undermines all his additional results. Such an outcome applies also to all other incompleteness results derived from the self-halting problem, including the reformulation of Gödel's First Theorem [59].

Another significant casualty of these refutations which should be mentioned is the solution to Hilbert's 10th problem, the existence or not of an algorithm capable of determining whether a given diaphontine equation is solvable [5,12,40]. This problem was finally answered in the negative, on the basis of the undecidability results previously mentioned [24,25,51]. In the absence of an alternative proof, its resolution and that of related problems [12] become again open questions.

7.6. **Proof of the completeness of first-order arithmetic**. The adaptations of Cantor's methods of proof used by Gödel and others show that the concept of arithmetical incompleteness has largely been linked to the notion of nondenumerability [59]. Given the extent to which incompleteness has penetrated into the fabric of mathematics, the only definitive answer to the questions raised by the refutations reported in Sections 7.1, 7.2, 7.3 and 7.5 would be to provide a sound proof for the completeness of first-order arithmetic.

The possibility of a theory of first-order arithmetic being negation-complete would appear more likely when considering that a host of reduced but formalised arithmetical theories have been shown to be negation-complete. This can be said of *Baby Arithmetic* (BA), which uses the successor, addition and multiplication functions, but is quantifier-free and lacks induction [59]. It also applies to *Presburger Arithmetic* (P), theory which amounts to PA stripped of multiplication [42,57,59]. Finally, the existence of a complete theory for the truths expressible in a first-order language with multiplication but lacking addition (or the successor function) was shown by Skolem in 1930 [58,59]. Hence, Gödel's work was responsible for the current belief that the combination of multiplication with addition and successor in one single theory produces incompleteness [31,59]. If proven, the completeness of first-order arithmetic (e.g. PA) would be placed alongside the completeness of predicate calculus [20] and the consistency of Euclidean and alternative geometries [5,44], fully reviving the ambitious goals of Hilbert's Programme [41].

A proof of the completeness of first-order arithmetic is presented below, derived from results found in the literature, combined with the denumerability of $\mathcal{P}(\mathbb{N})$ (Theorem 4.1). This proof shows that arithmetic completeness is a corollary of denumerability.

In order to prove that there is a suitable theory of first-order arithmetic that is negation-complete, the first requirement is to show that True Basic Arithmetic, $\mathcal{T}_A$, is effectively enumerable. Using $L_A$ to symbolize the formalised version of the language of basic arithmetic, $\mathcal{T}_A$ is the set of truths of $L_A$, i.e. the set of sentences which are true on the standard interpretation built into $L_A$. A suitable description of effective enumerability can be taken from [59]:

> "A set $\Sigma$ is *effectively enumerable* if an (idealized) computer could be programmed to generate a list of its members such that any member will eventually be mentioned – the list may be empty, or have no end, and may contain repetitions, so long as any item in the set eventually makes an appearance."

This definition determines what is expected of an algorithm with the capability to construct an effectively enumerable set.

Given a syntax $\mathcal{L}$ (i.e. the syntactically defined system of expressions which form part of a formal language $L$), there is a need to identify those strings of symbols that constitute the well-formed formulae of $\mathcal{L}$ (its *wff*s), particularly the closed $\mathcal{L}$-wffs, i.e. the $\mathcal{L}$-sentences. Quoting again from [59]:

> "So, whatever the details, for a properly formalized syntax $\mathcal{L}$, there should be clear and objective procedures, agreed on all sides, for *effectively deciding* whether a putative constant-symbol really is a constant, etc. Likewise we need to be able to effectively decide whether a string of symbols is an $\mathcal{L}$-wff/$\mathcal{L}$-sentence. It goes almost without saying that the formal languages familiar from elementary logic have this feature."

In other words, the ability to identify $\mathcal{L}$-sentences is an inherent capability of any formal language $L$. Therefore, it can be stated with confidence that $L_A$ has the capability to identify all its $\mathcal{L}_A$-sentences, by the operation of a suitable algorithm. Some of these sentences will be the truths of $L_A$, i.e. elements of the set $\mathcal{T}_A$.

The proof of completeness presented below uses a powerful result available in the literature. This is a corollary derived from Craig's Re-axiomatization Theorem [21,59], which can be stated as [59]:

> **Theorem** *If $\Sigma$ is a set of wffs which can be effectively enumerated, then there is a p.r. (primitive recursive) axiomatized theory $T$ whose theorems are exactly the members of $\Sigma$.*

One of the defining characteristic of a p.r. axiomatized theory $T$ is that the numerical properties of being the Gödel number (g.n.) of a $T$-wff/$T$-sentence are primitive recursive, i.e. decidable by p.r. functions [59]. This requirement is perfectly satisfied by PA, standard formulation of first-order arithmetic [59].

Having in mind the definitions and results just described, a theorem of completeness for first-order arithmetic follows.

**Theorem 7.1 (Arithmetic Completeness).** *Let $\mathcal{T}_A$ (True Basic Arithmetic) be the set of truths of $L_A$. It is the case that:*
*i) $\mathcal{T}_A$ is an effectively enumerable set.*
*ii) There is a p.r. axiomatized theory $T_A$ whose theorems are exactly the members of $\mathcal{T}_A$. Theory $T_A$ is sound (hence consistent), negation-complete, decidable and (assuming its p.r. adequacy) able to prove its own consistency.*

*Proof.* The crucial requirement is the construction of effective enumerations.

*i)* Consider $\Sigma_A$, the set of all $\mathcal{L}_A$-sentences (i.e. closed $\mathcal{L}_A$-wffs) of the formal language $L_A$. The first step is to show that the set $\Sigma_A$ is effectively enumerable. This can be achieved by making use of the Gödel numbering scheme since every g.n. is a natural number, i.e. $\forall \varphi \in \Sigma_A, \ulcorner \varphi \urcorner \in \mathbb{N}$. The aim will be to construct $\Sigma_A$.

An algorithm can be designed to: *a)* run through all the natural numbers $n \in \mathbb{N}$; *b)* treat them as potential g.n. and decode them into a string of symbols of $\mathcal{L}_A$ whenever possible; *c)* if a given string accounts for a $\mathcal{L}_A$-sentence (discerning such a fact is a basic capability of any formal language), incorporate it into $\Sigma_A$; *d)* identify each $\varphi \in \Sigma_A$ by its g.n. ($\ulcorner \varphi \urcorner$) and use this labelling as a one-to-one correspondence between $\Sigma_A$ and $N_A$, the set of all those natural numbers that double as the g.n. of a $\mathcal{L}_A$-sentence, i.e. $N_A = \{ n \in \mathbb{N}, \forall n ( n = \ulcorner \varphi \urcorner, \varphi \in \Sigma_A ) \}$. Since it is apparent that $N_A \subseteq \mathbb{N}$, the construction described produces an effective enumeration of the elements of $\Sigma_A$, as required.

Every sentence listed in $\Sigma_A$ is either true or false. Without having to determine the actual truth or falsehood of every sentence $\varphi \in \Sigma_A$, it can be stated that $\mathcal{T}_A \subseteq \Sigma_A$. Thus, all true sentences of $L_A$ (i.e. the elements $\mathcal{T}_A$) are elements of $\Sigma_A$, implying that the construction of all subsets of $\Sigma_A$ guarantees the construction of $\mathcal{T}_A$.

If $\Sigma_A$ were considered a finite set, the proof of the theorem would be immediate. It is reasonable to assume that set $\Sigma_A$ is infinite. This implies that a one-to-one correspondence can be established between $\Sigma_A$ and $\mathbb{N}$. Such a bijection enables the construction of the power set of $\Sigma_A$, $\mathcal{P}(\Sigma_A)$, simply by replicating the construction of the power set $\mathcal{P}(\mathbb{N})$ used in the proofs of the denumerability of $\mathcal{P}(\mathbb{N})$ (Theorem 4.1). The requirement is to maintain the relationship between each element of $\Sigma_A$ and the corresponding natural number $n \in \mathbb{N}$ throughout the construction of $\mathcal{P}(\Sigma_A)$. In this way, the construction of all the subsets of $\mathbb{N}$ can be matched by the parallel construction of all the subsets of $\Sigma_A$. Since $\mathcal{T}_A \subseteq \Sigma_A$, the above construction includes the construction of $\mathcal{T}_A$, implying that $\mathcal{T}_A$ is effectively enumerable. □

*ii)* Having proved that $\mathcal{T}_A$ is an effectively enumerable set, the proof of the second part of the theorem is straightforward. The construction of $\mathcal{P}(\Sigma_A)$ makes all the subsets of $\Sigma_A$ effectively enumerable. Consequently, a direct implication of the corollary of Craig's Re-axiomatization Theorem previously described is the existence of a complete collection of p.r. axiomatized theories whose theorems are exactly the members of each subset of $\Sigma_A$. Therefore, one of these theories will be $T_A$, i.e. a p.r. axiomatized theory whose theorems are exactly all the elements of $\mathcal{T}_A$, meaning all the truths of $L_A$. The construction of $\mathcal{T}_A$, as part of the construction of $\mathcal{P}(\Sigma_A)$, ensures both that the theory $T_A$ actually exists and that it is unique. Furthermore, a number of properties of $T_A$ can also be derived.

*a)* $T_A$ is sound. Since the theorems of $T_A$ are all and only all the truths of $L_A$, no theorem of $T_A$ can be false. Consequently, $T_A$ is sound. Furthermore, accepting the Law of Excluded Middle [22,52] implies that, if $T_A$ is sound, then $T_A$ is also consistent (the negation of any sentence proved by $T_A$ is false, hence it cannot be a theorem of $T_A$).

*b)* $T_A$ is negation-complete. It is apparent that $T_A$ decides every sentence of $L_A$. When a sentence $\varphi$ of $L_A$ is true, then $\varphi \in \mathcal{T}_A$ and, accordingly, it is a theorem of $T_A$, $T_A \vdash \varphi$. If the sentence $\varphi$ is false, its negation $\neg\varphi$ is true, hence $\neg\varphi \in \mathcal{T}_A$ and $T_A \vdash \neg\varphi$.

*c)* $T_A$ is decidable. It is a well-established fact that any consistent, axiomatized, negation-complete formal theory is decidable [59]. Since $T_A$ satisfies all these conditions, then $T_A$ is decidable. Furthermore, any consistent, recursively axiomatized, negation-complete formal theory is also recursively decidable [59]. According to Craig's Theorem, any axiomatized theory can be given a p.r. re-axiomatization; therefore, the fact that $T_A$ is p.r. axiomatized suggests that it is recursively decidable.

*d)* $T_A$ is able to prove its own consistency. The confirmation of this statement uses some properties and definitions applicable to any p.r. axiomatized theory, so defined when the numerical properties of being the g.n. of an axiom or a wff/sentence of such a theory are all p.r. The same applies to the numerical properties of being the super g.n. of a proof properly constructed from axioms of the theory [59]. The assumption that $T_A$ is p.r. adequate (i.e. it can capture every p.r. function as a function, hence also capturing every p.r. property and relation [59]) allows the replication of some of the elements originally used by Gödel in his own incompleteness proofs [31,59]. Three p.r. properties/relations are used.

- $Sent(n)$, a p.r. property which holds when $n$ is the g.n. of a sentence (i.e. a closed wff) of $L_A$. The p.r. adequacy of $T_A$ ensures the existence of a $\Sigma_1$-wff that canonically captures *Sent*. Following standard convention, this wff can be symbolised as $\mathsf{Sent(x)}$.

- $Prf(m,n)$, which holds when $m$ is the super g.n. of a sequence of wffs that is a $T_A$-proof of a sentence with g.n. $n$. Again, the p.r. adequacy of $T_A$ guarantees the existence of a $\Sigma_1$-wff that canonically captures *Prf*. This wff can be symbolised as $\mathsf{Prf(x,y)}$.

- Define $Prov(n) := \exists v\ Prf(v,n)$, which holds when some number $v$ codes for a $T_A$-proof of the sentence with g.n. $n$. Once more, the p.r. adequacy of $T_A$ implies there will be a $\Sigma_1$-wff that canonically captures *Prov*, i.e. $\mathsf{Prov(x)} := \exists \mathsf{v}\ \mathsf{Prf(v,x)}$, a (canonical) provability predicate for $T_A$ [59].

Since inconsistent classical theories prove every sentence [59], these three predicates are used to construct a new sentence capable of conveying the consistency of $T_A$: $\mathsf{Con} := \exists \mathsf{x}(\mathsf{Sent(x)} \wedge \neg \mathsf{Prov(x)})$. $\mathsf{Con}$ means that there is at least one $\mathcal{L}_A$-sentence that $T_A$ cannot prove. Accordingly, $T_A$ will be consistent if and only if the sentence $\mathsf{Con}$ is true. If Con is a true sentence of $L_A$, then it will follow that $\mathsf{Con} \in \mathcal{T}_A$ and, consequently, $T_A \vdash \mathsf{Con}$. Since $T_A$ has already been proven to be consistent, the final implication is that $T_A$ (provided it is p.r. adequate) is also able to prove its own consistency. □

The proof of Theorem 7.1 (Arithmetic Completeness) does not provide a precise identity of theory $T_A$, despite determining its existence. It could be conjectured that $T_A$ and PA are actually the same theory. Since PA is p.r. adequate [59], this would allow the elimination of the assumption made for $T_A$ in the last statement of Theorem 7.1. Then, the complete set of goals of Hilbert's Programme would be achieved. Furthermore, since Theorem 7.1 proves that $T_A$ is sound (as well as negation-complete), it could be said that an additional goal has been reached, beyond the formalist objective of Hilbert's Programme.

The validity of Theorem 7.1 has two major consequences. Firstly, it rules out the veracity of any other claims of arithmetic incompleteness, beyond those refuted in this article (Sections 7.1, 7.2, 7.3 and 7.5). Secondly, it firmly establishes the conviction that, however difficult it might be, all reasonable arithmetic statements can be proved true or false. Such a conviction applies equally to all $\Pi_1$-sentences, despite current preconceptions [59]. This conclusion should prompt a review of mathematical problems currently considered unsolvable (such as the previously mentioned Hilbert's 10th problem).

Arithmetic completeness has far reaching implications, philosophical as well as mathematical, which exceed the scope of this article. However, one area of interest that could be cited here is the potential development of artificial intelligence (AI). For decades, Gödelian incompleteness has been considered an unsummountable obstacle to the eventual construction of mechanical devices which could be considered 'intelligent', i.e. in some way able to replicate the scope of the human brain [59]. These perceived limitations have led to the conviction that some features of the physical world, beyond the classical domain, are required for the generation of consciousness [50,53,54]. It is clear that Theorem 7.1, and the refutations presented in the preceding sections, can provide a different scenario with new, and possibly raised expectations, of what is technologically achievable in this compelling area of research.

## 8. Concluding Remarks

As introduced in Section 2, the concept of nondenumerability plays a crucial role in Cantor's treatment of infinity. For this reason alone, the refutations presented in Section 3 seriously undermine the deepest foundations of the whole Cantorian landscape. Particular attention has been given to Cantor's diagonalization argument because of the widespread use by logicians of diagonal constructions [20,59], as discussed in Section 7. As a whole, the refutations of Section 3 have underlined a systemic failure in the correct implementation of proofs by contradiction that could extend beyond Cantorian set theory. Consequently, a workable strategy was needed to identify and eradicate this type of mistake. A simple and practical principle to achieve this goal was presented in Section 5, based on the formulation of Principle 5.2 (of Conceivable Proof). In particular, the use of self-referential statements should not be considered sound, and their deployment (as in the well-known Russell's Paradox [49,56]) might be seen as the cause of many of the difficulties that troubled the early development of set theory. Conjecture 5.3 (of Logical Imperfection) acknowledges the impossibility of preventing inconceivable statements by axiomatic means alone.

Theorem 4.1 (Denumerability of the Power Set), as well as corroborating the refutations of Section 3, implies that all infinite sets are countable. Although Theorem 4.1 cannot, by itself, rule out the existence of transfinite numbers, the denumerabilities of $\mathcal{P}(\mathbb{N})$ and higher power sets make it difficult to think of any other mathematical fact(s) capable of substantiating transfinite theory. It is important to remark here that Cantor's work with point sets [23], which predated his diagonalization argument, did not provide any solid evidence of nondenumerability beyond the earlier proof of 1874 [6], refuted in Section 3.5. Furthermore, as discussed in Section 7.4, the elementary proof of Goodstein's Theorem [55] undermines the whole structure of current transfinite theory.

One particularly interesting consequence of the denumerability of $\mathbb{R}$ revisits Cantor's first proof of nondenumerability published in 1874 [6]. Cantor used his result to deduce that the density of transcendental numbers in the real number line is infinitely greater than that of algebraic numbers, whose countable nature he managed to verify [6]. Cantor's claim on the relative density of algebraic and transcendental numbers was made on the basis of $\mathbb{R}$ being an uncountable set. Therefore, without any other independent results supporting this assertion, there is no choice but to declare this issue an open question.

Section 2 cited the results published by Gödel and Cohen in relation to the independence of CH (GCH) from ZF [17,18,32,33]. The countable nature of $\mathbb{R}$ and $\mathcal{P}(\mathbb{N})$ (Theorems 4.1 and 4.4, and Corollary 4.3) provides a plausible explanation for Gödel's and Cohen's claims. Given that CH makes the initial assumption that $\boldsymbol{c} > \aleph_0$, it can be deduced that the opposing statements of the truth or falsehood of CH are both false. Since the ZF axioms are designed to account for transfinite numbers [56], they are also compromised. The refutations and proofs presented here imply that no result involving the standard formulation of CH and ZF can be meaningful, particularly when related to issues of consistency.

The pre-eminence that Hilbert gave to CH [5,40] could be better understood with an evaluation of the impact that denumerability has on many areas of mathematics. Although Sections 5, 6 and 7 address some of the most immediate set-theoretical and logical implications of the denumerability of both the set of the real numbers $\mathbb{R}$, and the power set of $\mathbb{N}$, $\mathcal{P}(\mathbb{N})$, a complete audit of all the possible consequences will be required. Special attention should be given to the analytical implications brought about by the countable nature of the real number line [19].

Set theory can be largely simplified once inconceivable proofs and transfinite numbers are disregarded. Conjecture 6.1 (of Countable Infinity) endorses this principle. Once this simplification has been accepted, the interpretation of the axiom of infinity enshrined by Theorem 3.6 (of Actual Countable Infinity) could be sufficient to account for the construction of all infinite sets, given their countable nature. The concept of relative cardinality (Definitions 6.3 and 6.5) could be used to compare the size of infinite sets in a way not previously envisaged.

Diagonal arguments are widely used in mathematical logic. Therefore, their refutation as a method of proof (Section 7) has significant implications. The intimate link between nondenumerability and incompleteness, and their parallel refutation, has been further highlighted by the formulation of Theorem 7.1 (Arithmetic Completeness), with a proof which is a direct consequence of the denumerability of $\mathcal{P}(\mathbb{N})$ (Theorem 4.1). The completeness of first-order arithmetic thus deduced is poised to have major repercussions throughout mathematics and science at large.

A concluding remark has to be made in view of the primary purpose of the axiomatization of mathematical theory. As endorsed by the title of this article, "Addressing Mathematical Inconsistency", the ultimate objective should be the provision of a system of mathematics known to be consistent. In this sense, the appearance of any internal conflicts (which often manifest themselves in the form of seemingly intractable paradoxes [35,39,49,56]) should be interpreted, not necessarily as the inevitable limitations of mathematical scope, but as a call for a re-examination of the foundational pillars of that system. If the set of axioms and definitions that sit at the very heart of mathematical theory are internally consistent, they will be much more likely to always be perceived that way. This was the idealistic view embodied by Hilbert's Programme [41]. The fact that the concepts of denumerability and completeness, as defended in this article (Sections 4 and 7), have already been corroborated by a previous report of an elementary proof of Goodstein's Theorem [55], should be perceived as tantalizing evidence that, after all, Hilbert's goal can be achieved.

---

[19] A report on this topic is under preparation.

JUAN A PEREZ. BERKSHIRE, UK.
*E-mail address*: jap717@juanperezmaths.com

## Appendix: Commentary on a Recent Article

The integrity of Cantor's diagonalization argument has been questioned in an article recently published on arXiv [a]. Its author (A.L.S.) uses the idea of reordering a table containing all the real numbers (or just the rational numbers) within the unit interval $[0,1)$, with the purpose of rendering a rational diagonal number and a rational antidiagonal counterpart and implying a contradiction in Cantor's argument.

The above article correctly acknowledges the present author (J.A.P.) as the originator of the idea, which was part of an early draft of this current paper [c]. The draft was sent to A.L.S. and others in February/March 2009, together with the draft of another article [b]. However, a trivial error was noted in the idea by one of the reviewers (K.S.), hence it was not included in the submitted version of this current paper [c].

The invalidation applies to the idea used in [a], as much as it did to the original concept. Simply explained, an antidiagonal number (either before or after any reordering of the original table) cannot be part of a diagonal enumeration. Consequently, assuming that such an enumeration accounts for all the elements of the set under consideration, the process of reordering excludes the above antidiagonal and makes the list incomplete. Hence, no contradiction of Cantor's rationale can be achieved in this way, contrary to the claim made in [a].

An explicit way to illustrate this conclusion is to construct a table of all the rational numbers in the interval $[0,1)$ and, after selecting a suitable antidiagonal number, attempt the required reordering¶. To build a complete (countable) list of all the rational numbers within the interval $[0,1)$, $\mathbb{Q}_{01}$, take all the fractions $q=a/b$, $\forall a,b\in\mathbb{N}$, $a<b \wedge (a,b)=1$, i.e. $a$ and $b$ are coprime. Including all the possible values for $a$, as $b$ runs through all the natural numbers, the list will start with:

$$\text{(A.1)} \qquad \mathbb{Q}_{01} = \{0,\, 1/2,\, 1/3,\, 2/3,\, 1/4,\, 3/4,\, 1/5,\, 2/5,\, 3/5,\, 4/5,\, 1/6,\, 5/6,\, 1/7,\, 2/7,\, 3/7,\, 4/7,\, \cdots\},$$

where the first 16 members of the table are listed. Display the elements of $\mathbb{Q}_{01}$ using their binary representations and, for simplicity, when two could be used (as it is the case, for example, of $q=1/2=0.1000\cdots$ or $0.0111\cdots$), select only the representation with '0' as terminating digit.

Without any prior reordering, the binary representations of the rational numbers that appear in (A.1) will be as follows:

¶ A similar construction for all the real numbers in the unit interval $[0,1)$ could be considered, and the same conclusion would be reached.

$$
\begin{array}{rcl}
0 & = & 0\,.\,\mathbf{0}\,0\,0\,0\,0\,0\,0\,0\,0\,0\,0\,0\,0\,0\,0\,0\cdots \\
1/2 & = & 0\,.\,1\,\mathbf{0}\,0\,0\,0\,0\,0\,0\,0\,0\,0\,0\,0\,0\,0\,0\cdots \\
1/3 & = & 0\,.\,0\,1\,\mathbf{0}\,1\,0\,1\,0\,1\,0\,1\,0\,1\,0\,1\,0\,1\cdots \\
2/3 & = & 0\,.\,1\,0\,1\,\mathbf{0}\,1\,0\,1\,0\,1\,0\,1\,0\,1\,0\,1\,0\cdots \\
1/4 & = & 0\,.\,0\,1\,0\,0\,\mathbf{0}\,0\,0\,0\,0\,0\,0\,1\,0\,0\,0\,0\cdots \\
3/4 & = & 0\,.\,1\,1\,0\,0\,0\,\mathbf{0}\,0\,0\,0\,0\,0\,1\,0\,0\,0\,0\cdots \\
1/5 & = & 0\,.\,0\,0\,1\,1\,0\,0\,\mathbf{1}\,1\,0\,0\,1\,1\,0\,0\,1\,1\cdots \\
2/5 & = & 0\,.\,0\,1\,1\,0\,0\,1\,1\,\mathbf{0}\,0\,1\,1\,0\,0\,1\,1\,0\cdots \\
3/5 & = & 0\,.\,1\,0\,0\,1\,1\,0\,0\,1\,\mathbf{1}\,0\,0\,1\,1\,0\,0\,1\cdots \\
4/5 & = & 0\,.\,1\,1\,0\,0\,1\,1\,0\,0\,1\,\mathbf{1}\,0\,0\,1\,1\,0\,0\cdots \\
1/6 & = & 0\,.\,0\,0\,1\,0\,1\,0\,1\,0\,1\,0\,\mathbf{1}\,0\,1\,0\,1\,0\cdots \\
5/6 & = & 0\,.\,1\,1\,0\,1\,0\,1\,0\,1\,0\,1\,0\,\mathbf{1}\,0\,1\,0\,1\cdots \\
1/7 & = & 0\,.\,0\,0\,1\,0\,0\,1\,0\,0\,1\,0\,0\,1\,\mathbf{0}\,0\,1\,0\cdots \\
2/7 & = & 0\,.\,0\,1\,0\,0\,1\,0\,0\,1\,0\,0\,1\,0\,0\,\mathbf{1}\,0\,0\cdots \\
3/7 & = & 0\,.\,0\,1\,1\,0\,1\,1\,0\,1\,1\,0\,1\,1\,0\,1\,\mathbf{1}\,0\cdots \\
4/7 & = & 0\,.\,1\,0\,0\,1\,0\,0\,1\,0\,0\,1\,0\,0\,1\,0\,0\,\mathbf{1}\cdots
\end{array}
\tag{A.2}
$$

By examination of (A.2), it is apparent that the current diagonal number does not seem to have a periodic binary representation. Therefore, to achieve this objective, a reordering is required. Select as antidiagonal number $q_{AD} = 2/3 = 0.\overline{10} = 0.10101010\cdots$, which will need an order to be displayed such that the diagonal number is $q_D = 1/3 = 0.\overline{01} = 0.01010101\cdots$. It should be obvious that any possible reordering of the table (A.2) (and the subsequent elements of $\mathbb{Q}_{01}$) that renders the diagonal sequence associated with $q_D = 1/3$, will exclude $q_{AD} = 2/3$. It should be equally obvious that, whatever choice for $q_D$ ($q_{AD}$) is made, the resulting antidiagonal number will not appear in the diagonal enumeration.

However, the above analysis is initiated with the assumption that, if the set tabled is countably infinite, a diagonal enumeration should always be able to account for all the elements of the set. As Section 3.2.2 of the current paper [*c*] explained, such an assumption is erroneous. In fact, an infinitely large number of denumerable sets could be conceived such that their diagonal cover $Dc$ (as defined in Section 3.2.2 of this article [*c*]) is $Dc < 1$. This is the case for the set $\mathbb{Q}_{01}$ used in the preceding example and with the first 16 elements displayed in table (A.2). Assume that a reordering of (A.2) has been executed that, as mentioned above, renders the antidiagonal number $q_{AD} = 2/3 = 0.\overline{10} = 0.10101010\cdots$. That reordering will exclude $q_{AD} = 2/3$, which will no longer appear in any visible list. However, consider three more elements of $\mathbb{Q}_{01}$: $q_1 = 1/6 = 0.0\overline{01} = 0.\mathbf{0}0101010\cdots$, that differs from 2/3 in the *first* digit after the decimal point, $q_2 = 11/12 = 0.11\overline{10} = 0.1\mathbf{1}101010\cdots$, that differs from 2/3 in the *second* digit, and finally $q_3 = 5/12 = 0.01\overline{10} = 0.\mathbf{01}101010\cdots$, that differs in the *first two* digits. To avoid exclusion, $q_1$ and $q_2$ have to appear in the first and second positions, respectively, of the new ordering. But then there will be no positions left for $q_3$, which will inevitably be excluded alongside $q_{AD}$. In the same way as $q_3$, there will be a myriad of other elements of $\mathbb{Q}_{01}$ that will also be excluded.

The above example is representative of any other conceivable enumeration of the elements of $\mathbb{Q}_{01}$. Consequently, the reordering of any list of rational numbers cannot be used as a reason to question Cantor's diagonalization argument.